\documentclass[10pt]{amsart}
\usepackage{amsmath}
\usepackage{amsxtra}
\usepackage{amscd}
\usepackage{amsthm}
\usepackage{amsfonts}
\usepackage{amssymb}
\usepackage{eucal}

\usepackage{latexsym}

\newtheorem{cor}[subsection]{Corollary}
\newtheorem{lem}[subsection]{Lemma}
\newtheorem{prop}[subsection]{Proposition}
\newtheorem{propconstr}[subsection]{Proposition-Construction}

\newtheorem{thm}[subsection]{Theorem}

\newtheorem{lemconstr}[subsection]{Lemma-Construction}
\theoremstyle{remark}

\theoremstyle{definition}

\numberwithin{equation}{section}

\newcommand{\propconstrref}[1]{Proposition-Construction~\ref{#1}}
\newcommand{\lemconstrref}[1]{Lemma-Construction~\ref{#1}}
\newcommand{\thmref}[1]{Theorem~\ref{#1}}

\newcommand{\secref}[1]{Sect.~\ref{#1}}
\newcommand{\lemref}[1]{Lemma~\ref{#1}}
\newcommand{\propref}[1]{Proposition~\ref{#1}}
\newcommand{\corref}[1]{Corollary~\ref{#1}}

\newcommand{\nc}{\newcommand}
\nc{\renc}{\renewcommand}
\nc{\ssec}{\subsection}
\nc{\sssec}{\subsubsection}
\nc{\on}{\operatorname}

\nc\ol{\overline}
\nc\ul{\underline}
\nc\wt{\widetilde}
\nc\tboxtimes{\wt{\boxtimes}}
\nc{\wh}{\widehat}
\nc{\mc}{\mathcal}

\nc{\CM}{{\mathcal M}}
\nc{\CN}{{\mathcal N}}
\nc{\CF}{{\mathcal F}}
\nc{\D}{{\mathcal D}}
\nc{\CQ}{{\mathcal Q}}
\nc{\CY}{{\mathcal Y}}
\nc{\CX}{{\mathcal X}}
\nc{\CG}{{\mathcal G}}
\nc{\CE}{{\mathcal E}}
\nc{\CC}{{\mathcal C}}
\nc{\CO}{{\mathcal O}}
\renc{\CC}{{\mathcal C}}
\nc{\CT}{{\mathcal T}}
\nc{\CK}{{\mathcal K}}
\nc{\CS}{{\mathcal S}}
\nc{\CH}{{\mathcal H}}
\nc{\CU}{{\mathcal U}}
\nc{\CV}{{\mathcal V}}
\nc{\CA}{{\mathcal A}}
\nc{\CB}{{\mathcal B}}
\nc{\CW}{{\mathcal W}}
\nc{\CL}{{\mathcal L}}
\nc{\CP}{{\mathcal P}}
\nc{\CI}{{\mathcal I}}
\nc{\CJ}{{\mathcal J}}
\nc{\CR}{{\mathcal R}}
\nc{\CZ}{{\mathcal Z}}

\nc{\BA}{{\mathbb{A}}}
\nc{\BC}{{\mathbb{C}}}
\nc{\BG}{{\mathbb{G}}}
\nc{\BM}{{\mathbb{M}}}
\nc{\BN}{{\mathbb{N}}}
\nc{\BP}{{\mathbb{P}}}
\nc{\BR}{{\mathbb{R}}}
\nc{\BZ}{{\mathbb{Z}}}
\nc{\BV}{{\mathbb{V}}}
\nc{\BW}{{\mathbb{W}}}
\nc{\BS}{{\mathbb{S}}}
\nc{\BD}{{\mathbb{D}}}
\nc{\BQ}{{\mathbb{Q}}}
\nc{\BL}{{\mathbb{L}}}
\renc{\BW}{{\mathbb{W}}}

\nc{\fa}{{\mathfrak{a}}}
\nc{\fb}{{\mathfrak{b}}}
\nc{\fg}{{\mathfrak{g}}}
\nc{\fgl}{{\mathfrak{gl}}}
\nc{\fh}{{\mathfrak{h}}}
\nc{\fj}{{\mathfrak{j}}}
\nc{\fm}{{\mathfrak{m}}}
\nc{\fl}{{\mathfrak{l}}}
\nc{\fn}{{\mathfrak{n}}}
\nc{\fu}{{\mathfrak{u}}}
\nc{\fp}{{\mathfrak{p}}}
\nc{\ff}{{\mathfrak{f}}}
\nc{\fd}{{\mathfrak{d}}}
\nc{\fr}{{\mathfrak{r}}}
\nc{\fs}{{\mathfrak{s}}}
\nc{\fsl}{{\mathfrak{sl}}}

\nc{\hsl}{{\widehat{\mathfrak{sl}}}}
\nc{\hgl}{{\widehat{\mathfrak{gl}}}}
\nc{\hg}{{\widehat{\mathfrak{g}}}}
\nc{\hb}{{\widehat{\mathfrak{b}}}}
\nc{\hn}{{\widehat{\mathfrak{n}}}}

\nc{\fA}{{\mathfrak{A}}}
\nc{\fB}{{\mathfrak{B}}}
\nc{\fO}{{\mathfrak{O}}}
\nc{\fD}{{\mathfrak{D}}}
\nc{\fE}{{\mathfrak{E}}}
\nc{\fF}{{\mathfrak{F}}}
\nc{\fG}{{\mathfrak{G}}}
\nc{\fK}{{\mathfrak{K}}}
\nc{\fL}{{\mathfrak{L}}}
\nc{\fC}{{\mathfrak{C}}}
\nc{\fM}{{\mathfrak{M}}}
\nc{\fN}{{\mathfrak{N}}}
\nc{\fH}{{\mathfrak{H}}}
\nc{\fP}{{\mathfrak{P}}}
\nc{\fU}{{\mathfrak{U}}}
\nc{\fV}{{\mathfrak{V}}}
\nc{\fZ}{{\mathfrak{Z}}}
\nc{\fz}{{\mathfrak{z}}}

\nc{\bc}{{\mathbf{c}}}
\nc{\bd}{{\mathbf{d}}}
\nc{\bh}{{\mathbf{h}}}
\nc{\be}{{\mathbf{e}}}
\nc{\ba}{{\mathbf{a}}}
\nc{\bj}{{\mathbf{j}}}
\nc{\bn}{{\mathbf{n}}}
\nc{\bp}{{\mathbf{p}}}
\nc{\bg}{{\mathbf{g}}}
\nc{\bq}{{\mathbf{q}}}
\nc{\bs}{{\mathbf{s}}}
\nc{\bu}{{\mathbf{u}}}
\nc{\bv}{{\mathbf{v}}}
\nc{\bx}{{\mathbf{x}}}
\nc{\by}{{\mathbf{y}}}
\nc{\bb}{{\mathbf{b}}}
\nc{\bw}{{\mathbf{w}}}
\nc{\bA}{{\mathbf{A}}}
\nc{\bK}{{\mathbf{K}}}
\nc{\bB}{{\mathbf{B}}}
\nc{\bC}{{\mathbf{C}}}
\nc{\bD}{{\mathbf{D}}}
\nc{\bH}{{\mathbf{H}}}
\nc{\bM}{{\mathbf{M}}}
\nc{\bN}{{\mathbf{N}}}
\nc{\bV}{{\mathbf{V}}}
\nc{\bW}{{\mathbf{W}}}
\nc{\bL}{{\mathbf{L}}}
\nc{\bU}{{\mathbf{U}}}
\nc{\bX}{{\mathbf{X}}}
\nc{\bI}{{\mathbf{I}}}
\nc{\bZ}{{\mathbf{Z}}}
\nc{\bS}{{\mathbf{S}}}

\nc{\sA}{{\mathsf{A}}}
\nc{\sB}{{\mathsf{B}}}
\nc{\sC}{{\mathsf{C}}}
\nc{\sD}{{\mathsf{D}}}
\nc{\sF}{{\mathsf{F}}}
\nc{\sH}{{\mathsf{H}}}
\nc{\sE}{{\mathsf{E}}}
\nc{\sG}{{\mathsf{G}}}
\nc{\sK}{{\mathsf{K}}}
\nc{\sM}{{\mathsf{M}}}
\nc{\sO}{{\mathsf{O}}}
\nc{\sQ}{{\mathsf{Q}}}
\nc{\sP}{{\mathsf{P}}}
\nc{\sV}{{\mathsf{V}}}
\nc{\sZ}{{\mathsf{Z}}}
\nc{\sfp}{{\mathsf{p}}}
\nc{\sr}{{\mathsf{r}}}
\nc{\sg}{{\mathsf{g}}}
\nc{\sh}{{\mathsf{h}}}
\nc{\sk}{{\mathsf{k}}}
\nc{\ssf}{{\mathsf{f}}}
\nc{\sv}{{\mathsf{v}}}
\nc{\ssh}{{\mathsf{h}}}
\nc{\sse}{{\mathsf{e}}}
\nc{\sfb}{{\mathsf{b}}}
\nc{\sfc}{{\mathsf{c}}}
\nc{\sd}{{\mathsf{d}}}

\nc{\Av}{\on{Av}}
\nc{\act}{\on{act}}
\nc{\Hom}{\on{Hom}}
\nc{\Ext}{\on{Ext}}
\nc{\Tor}{\on{Tor}}
\nc{\End}{\on{End}}
\nc{\Lie}{\on{Lie}}
\nc{\Loc}{\on{Loc}}
\nc{\IC}{\on{IC}}
\nc{\Aut}{\on{Aut}}
\nc{\rk}{\on{rk}}
\nc{\Sh}{\on{Sh}}
\nc{\Perv}{\on{Perv}}
\nc{\pos}{{\on{pos}}}
\nc{\Conv}{\on{Conv}}
\nc{\Sph}{\on{Sph}}
\nc{\Sym}{\on{Sym}}
\nc{\Rep}{\on{Rep}}
\nc{\RepH}{{\mc R}ep(H)}
\nc{\Fun}{\on{Fun}}
\nc{\Id}{\on{Id}}
\nc{\id}{\on{id}}
\renc{\mod}{\on{--mod}}

\nc{\oG}{\overset{\circ}{G}{}}
\nc{\oGB}{{\overset{\circ}{G/B}{}}}
\nc{\oGN}{{\overset{\circ}{G/N}{}}}
\nc{\uBC}{\underline{\BC}}

\nc{\crit}{{\on{crit}}}
\nc{\reg}{{\on{reg}}}
\nc{\nilp}{{\on{nilp}}}
\nc{\ord}{\on{ord}}
\nc{\nil}{\wt{\on{reg}}}
\nc{\mb}{\mathbf}
\nc{\ren}{{\on{ren}}}
\nc{\res}{\on{res}}
\nc{\RS}{{\on{RS}}}
\nc{\Dist}{\on{Dist}}
\nc{\semiinf}{{\frac{\infty}{2}}}
\nc{\semiinfi}{{\frac{\infty}{2}+i}}
\nc{\semiinfb}{{\frac{\infty}{2}+\bullet}}
\nc{\torsemiinf}{{\overset{\semiinf}\otimes}}
\nc{\Hitch}{\on{Hitch}}

\nc{\hl}{\overset{\leftarrow}h}
\nc{\hr}{\overset{\rightarrow}h}
\nc\Dh{\widehat{\D}}
\nc{\Gr}{\on{Gr}}
\nc{\Grb}{\ol{\Gr}{}_G}
\nc{\Fl}{\on{Fl}}
\nc{\Flt}{\wt{\Fl}{}}
\nc{\Pic}{\on{Pic}}
\nc{\Bun}{\on{Bun}}

\nc{\bDR}{\mathbf {DR}}
\nc{\uV}{\underline{V}}
\nc{\arrowtimes}{\overset{\to}\otimes}
\nc{\hattimes}{\widehat\otimes}
\nc{\larrowtimes}{\overset{\leftarrow}\otimes}
\nc{\shriektimes}{\overset{!}\otimes}
\nc{\startimes}{\overset{*}\otimes}
\nc{\sCliff}{\mathsf {Cliff}}
\nc{\sSpin}{\mathsf {Spin}}

\nc{\one}{{\mathbf{1}}}

\nc\Spec{\on{Spec}}
\nc{\Pro}{\on{Pro}}
\nc{\QCoh}{\on{QCoh}}
\nc{\uHom}{\underline{\on{Hom}}}
\nc{\RHom}{\on{RHom}}
\nc{\uRHom}{\underline{\on{RHom}}}
\nc{\CHom}{{\mathcal Hom}}
\nc{\uCHom}{\underline{{\mathcal Hom}}}
\nc{\uCRHom}{\underline{{\mathcal R}{\mathcal Hom}}}

\nc{\cg}{{\check \fg}}
\nc{\Op}{\on{Op}}
\nc{\nOp}{\on{Op}^{\nilp}_{\cg}}
\nc{\nMOp}{\on{MOp}^{\nilp}_{\cg}}
\nc{\rOp}{\on{Op}^{\reg}_{\cg}}

\nc{\tg}{\wt{\check \fg}}
\nc{\cn}{\check \fn}
\nc{\tn}{\wt{\cn}}
\nc{\cG}{{\check G}}
\nc{\cB}{\check B}
\nc{\cT}{\check T}
\nc{\cH}{\check H}
\nc{\cb}{\check \fb}
\nc{\cN}{\check N}
\nc{\MOp}{\on{MOp}}
\nc{\tN}{\wt{\CN}_{\cG}}
\nc{\dIsom}{{\mathsf{Isom}}_{\Op}}
\nc{\disom}{{\mathsf{isom}}_{\Op}}
\nc{\Kdv}{{\mathsf{Isom}}_{\Op^\reg}}
\nc{\kdv}{{\mathsf{isom}}_{\Op^\reg}}
\nc{\Isom}{{\mathsf{Isom}}}
\nc{\isom}{{\mathsf{isom}}}

\nc{\wcosta}{j_{\wt{w},*}}
\nc{\wsta}{j_{\wt{w},!}}
\nc{\wcost}{j_{w,*}}
\nc{\wst}{j_{w,!}}

\nc{\epsi}{{\mathbf e}^\psi}
\nc{\epsip}{{\mathbf e}^{\psi'}}

\nc{\Ppi}{{\mathbf \Pi}}

\nc{\hCO}{{\hat{\CO}}}
\nc{\hCK}{{\hat{\CK}}}

\nc{\CPreg}{\CP_{G,\on{Op}^\reg}}
\nc{\CPBreg}{\CP_{B,\on{Op}^\reg}}
\nc{\CPnilp}{\CP_{G,\on{Op}^\nilp}}
\nc{\CPBnilp}{\CP_{B,\on{Op}^\nilp}}
\nc{\CPla}{\CP_{G,\on{Op}_{\cla}}}
\nc{\CPBla}{\CP_{B,\on{Op}_{\cla}}}

\nc{\Cat}{\hg_\crit\mod^{I,m}_\nilp}
\nc{\Catf}{{}^{fl}\hg_\crit\mod^{I,m}_\nilp}
\nc{\DCat}{D^b(\hg_\crit\mod_\nilp)^{I^0}}
\nc{\DCatf}{{}^{fl} D^b(\hg_\crit\mod_\nilp)^{I^0}}
\nc{\Catr}{\hg_\crit\mod^{I,m}_\reg}
\nc{\Catrf}{{}^{fl}\hg_\crit\mod^{I,m}_\reg}
\nc{\DCatr}{D^b(\hg_\crit\mod_\reg)^{I^0}}
\nc{\DCatrf}{{}^{fl} D^b(\hg_\crit\mod_\reg)^{I^0}}

\nc{\ch}{\mbox{ch}}
\nc{\Z}{{\mathbb Z}}
\nc{\C}{{\mathbb C}}
\nc{\pone}{{\mathbb C}{\mathbb P}^1}
\nc{\pa}{\partial}
\nc{\F}{{\mathcal F}}
\nc{\arr}{\rightarrow}
\nc{\larr}{\longrightarrow}
\nc{\al}{\alpha}
\nc{\ri}{\rangle}
\nc{\lef}{\langle}
\nc{\W}{{\mathcal W}}
\nc{\la}{\lambda}
\nc{\ep}{\epsilon}
\nc{\su}{\widehat{{\mathfrak s}{\mathfrak l}}_2}
\nc{\sw}{{\mathfrak s}{\mathfrak l}}
\nc{\g}{{\mathfrak g}}
\nc{\h}{{\mathfrak h}}
\nc{\n}{{\mathfrak n}}
\nc{\N}{\widehat{\n}}
\nc{\De}{\Delta}
\nc{\gt}{\widetilde{\g}}
\nc{\Ga}{\Gamma}
\nc{\z}{{\mathfrak Z}}
\nc{\La}{\Lambda}
\nc{\cri}{_{\kappa_c}}
\nc{\kk}{h^\vee}
\nc{\sun}{\widehat{\sw}_N}
\nc{\si}{\sigma}
\nc{\el}{\ell}
\nc{\bi}{\bibitem}
\nc{\om}{\omega}
\nc{\ds}{\displaystyle}
\nc{\dzz}{\frac{dz}{z}}
\nc{\Res}{\on{Res}}
\nc{\Cal}{\mathcal}
\nc{\ot}{\otimes}
\nc{\R}{{\mc R}}
\nc{\yy}{{\mc Y}}
\nc{\ga}{\gamma}

\nc{\us}{\underset}
\nc{\opl}{\oplus}
\nc{\beq}{\begin{equation}}
\nc{\Fq}{{\mathbb F}_q}
\nc{\Mq}{{\mathcal M}}
\nc{\lan}{\langle}
\nc{\ran}{\rangle}

\nc{\Vect}{\on{Vect}}
\nc{\ghat}{\wh\fg}
\nc{\T}{\mc T}
\nc{\Tloc}{\T^\g_{\on{loc}}}
\nc{\vac}{|0\ran}
\nc{\Wick}{{\mb :}}
\nc{\delz}{\partial_z}
\nc{\K}{{\cali K}}
\nc{\cali}{\mathcal}
\nc{\li}{\mathfrak l}
\nc{\lt}{\widetilde{\li}}
\nc{\astar}{a^*}
\nc{\cA}{{\mc A}}
\nc{\ka}{\kappa}

\nc{\OO}{{\mc O}}
\nc{\AutO}{\on{Aut}\OO}
\nc{\DerO}{\on{Der}\OO}
\nc{\DerpO}{\on{Der}_+\OO}
\nc{\mf}{\mathfrak}
\nc{\V}{{\mc V}}
\nc{\hh}{\wh{\h}}

\nc{\pp}{{\mathfrak p}}
\nc{\mm}{{\mathfrak m}}
\nc{\rr}{{\mathfrak r}}
\nc{\ket}{\rangle}
\nc{\zz}{{\mathfrak z}}
\nc{\gr}{\on{gr}}
\nc{\Spe}{\on{Spec}}
\nc{\rv}{\crho}
\nc{\can}{\on{can}}
\nc{\Db}{{\mathbb D}}
\nc{\ww}{w}

\nc{\RR}{\on{R}}
\nc{\PPi}{{\mathbf \Pi}}
\nc{\M}{{\mathbb M}}
\nc{\Mv}{{\mathbb M}^\vee}
\nc{\VV}{{\mathbb V}}
\nc{\bsl}{\backslash}

\nc{\bchi}{{\mathbf {\chi}}}
\nc{\anch}{{\mathbf {anch}}}

\nc{\cla}{{\check{\la}}}
\nc{\cmu}{{\check{\mu}}}
\nc{\crho}{{\check{\rho}}}
\nc{\com}{{\check{\omega}}}
\nc{\DD}{{\mc D}}
\nc{\E}{{\mc E}}
\nc{\Ll}{{\mc L}}

\nc{\ConnX}{\on{Conn}_{\cH}(\omega^{\rho}_X)}
\nc{\ConnD}{\on{Conn}_{\cH}(\omega_{\D}^{\rho})}
\nc{\ConnDt}{\on{Conn}_{\cH}(\omega_{\D^\times}^{\rho})}

\nc{\Hecke}{{\on{Hecke}}}

\nc{\cLambda}{{\check\Lambda}}
\nc{\cnu}{{\check\nu}}
\nc{\ceta}{{\check\eta}}

\nc{\Ind}{\on{Ind}}

\nc{\CTop}{{\mathcal Top}}

\nc{\ppart}{(\!(t)\!)}

\nc{\qu}{/\!/}

\nc{\gen}{\on{gen}}

\nc{\alg}{\on{alg}}
\nc{\geom}{\on{geom}}
\nc{\Jet}{\on{Jets}}
\nc{\aut}{\on{aut}}
\nc{\Der}{\on{Der}}

\nc{\uW}{\underline{W}}
\nc{\uU}{\underline{U}}

\nc{\aff}{\on{aff}}
\nc{\dirsum}{\oplus}

\begin{document}

\title{Geometric realizations of Wakimoto modules at the critical level}

\author{Edward Frenkel}\thanks{The research of E.F. was supported by
the DARPA grant HR0011-04-1-0031 and by the NSF grant DMS-0303529.}

\address{Department of Mathematics, University of California,
  Berkeley, CA 94720, USA}

\email{frenkel@math.berkeley.edu}

\author{Dennis Gaitsgory}

\address{Department of Mathematics, Harvard University,
Cambridge, MA 02138, USA}

\email{gaitsgory@math.harvard.edu}

\date{March 2006}

\begin{abstract}

We study the Wakimoto modules over the affine Kac-Moody algebras at
the critical level from the point of view of the equivalences of
categories proposed in our previous works, relating categories of
representations and certain categories of sheaves. In particular, we
describe explicitly geometric realizations of the Wakimoto modules as
Hecke eigen-D-modules on the affine Grassmannian and as quasi-coherent
sheaves on the flag variety of the Langlands dual group.

\end{abstract}

\maketitle

\section*{Introduction}

Wakimoto modules, introduced in \cite{Wak,FF:usp,FF:si}, have many
applications in representation theory of affine Kac-Moody algebras. In
our previous papers \cite{FG1}--\cite{FG4} we have undertaken a study
of representations of affine Kac-Moody algebras at the critical level
in the framework of the local geometric Langlands correspondence.
Wakimoto modules play an important role in it. In this paper we
elucidate further the geometric meaning of Wakimoto modules from the
point of view of the equivalences of categories proposed in
\cite{FG2,FG4}. These equivalences relate categories of
representations of the affine Kac-Moody algebras at the critical level
to certain categories of D-modules and quasi-coherent sheaves.

\ssec{}

The first step in establishing a connection between representation
theory of an affine Kac-Moody algebra at the critical level $\hg_\crit$
and geometric Langlands correspondence is the description of the
center $\fZ_\fg$ of the completed universal enveloping algebra of
$\hg_\crit$. According to a theorem of \cite{FF,F:wak}, $\fZ_\fg$ is
isomorphic to the algebra of functions on the space $\Op_{\cg}({\mc
D}^\times)$ of $\cg$-opers on the punctured disc, where $\cg$ is the
Langlands dual Lie algebra to $\fg$. This isomorphism was proved in
\cite{FF,F:wak} algebraically, in the framework of representation
theory of $\hg_\crit$, and in particular using Wakimoto
modules. It may be reformulated as an isomorphism
\begin{equation}    \label{first}
\on{map}_{\alg}: \Spec(\fz_{\fg,X}) \overset{\sim}\longrightarrow
\Op_{\cg,X},
\end{equation}
defined for an arbitrary smooth algebraic curve $X$. Here
$\fz_{\fg,X}$ is the center of the chiral algebra on $X$ associated to
$\hg_\crit$ and $\Op_{\cg,X}$ is the scheme of jets of $\cg$-opers on
$X$.

We recall from \cite{BD} (see also \cite{FG2}, Sect. 1) that a
$\cg$-oper on $X$ is a $\cG$-bundle on $X$ equipped with a reduction
to a Borel subgroup $\cB$ and a connection satisfying a certain
transversality condition. This means that $\Op_{\cg,X}$ carries a
tautological $\cG$-bundle equipped with a $\cB$-reduction and a
connection along $X$. Defining a morphism $\Spec(\fz_{\fg,X}) \to
\Op_{\cg,X}$ is equivalent to defining the pull-backs of these data to
$\Spec(\fz_{\fg,X})$.

Can we construct these data on $\Spec(\fz_{\fg,X})$ in a natural way?
This question was addressed by A. Beilinson and V. Drinfeld in
\cite{BD}. They constructed these data using the affine Grassmannian
$\Gr_G = G\ppart/G[[t]]$ and the geometric Satake correspondence which
identifies the category of $G[[t]]$-equivariant D-modules on $\Gr_G$
with the category $\Rep(\cG)$ of representations of $\cG$, see
\cite{MV}. Thus, they obtained a map
$$
\on{map}_{\geom}: \Spec(\fz_{\fg,X}) \to \Op_{\cg,X}.
$$

We recall its definition in \secref{def geom} below. We note that the
construction of the map $\on{map}_{\geom}$ in \cite{BD} relies on the
existence of the isomorphism \eqref{first}. However, {\em a
priori} it is not clear whether $\on{map}_{\geom}$ coincides with
$\on{map}_{\alg}$, nor whether $\on{map}_{\geom}$ is an
isomorphism. Beilinson and Drinfeld proved this in \cite{BD} by
showing that both of these maps are compatible with actions of
certain Lie algebroids and that this property essentially
characterizes these maps uniquely.

In Sects. \ref{comparison}--\ref{geom constr} of this paper we give a
different proof of the fact that the maps $\on{map}_{\alg}$ and
$\on{map}_{\geom}$ coincide. This proof uses Wakimoto modules in
an essential way, in particular, their behavior under the
"Harish-Chandra convolution" functors which was described in
\cite{FG2}. Using this additional structure, we will see the emergence
in representation theoretic context not only of the geometric data of
opers mentioned above (the "birth of opers", as Beilinson and Drinfeld
had put it), but also the geometric data of Miura opers which
parametrize Wakimoto modules.

We remark that our proof of the coincidence of $\on{map}_{\alg}$ and
$\on{map}_{\geom}$ also relies on the results of \cite{FF,F:wak}, in
which the existence of \eqref{first} was proved, and so it does not
give us an alternative proof of the isomorphism
\eqref{first}. However, it helps us understand better the geometric
meaning of this isomorphism and the role of Wakimoto modules in
it. This is the first main result of this paper.

\ssec{} Next, we analyze the Wakimoto modules from the point of view
of the equivalences of categories that appeared in our approach to the
local geometric Langlands correspondence in \cite{FG2,FG4}. 

Let us denote by $\hg_\crit\mod_\reg$ the category of discrete
$\hg$-modules at the critical level on which the center $\fZ_\fg$ acts
through its quotient $\fZ^\reg_\fg = \on{Fun}(\Op_{\cG}^\reg)$, where
$\Op_{\cG}^\reg \subset \Op_{\cG}({\mc D}^\times)$ is the space of
$\cg$-opers on the formal disc ${\mc D}$. Let
$\hg_\crit\mod_\reg^{I^0}$ be its subcategory of $I^0$-equivariant
$\hg$-modules, where $I^0$ is the radical of the Iwahori subgroup $I
\subset G\ppart$. The algebra $\on{Fun}(\Op_{\cG}^\reg)$ acts on the
category $\hg_\crit\mod_\reg^{I^0}$ in a natural way, so we may think
of $\hg_\crit\mod_\reg^{I^0}$ as "fibered" over the space
$\Op_{\cG}^\reg$. In this Introduction, in order to simplify our
notation, we will restrict ourselves to a "fiber" of this category,
denoted by $\hg_\crit\mod_\chi^{I^0}$, over a particular $\cg$-oper
$\chi$. The objects of $\hg_\crit\mod_\chi^{I^0}$ are
$I^0$-equivariant $\hg$-modules at the critical level on which the center
$\fZ_\fg$ acts through the character determined by $\chi$. (In the
main body of the paper we will work over the base $\Op_{\cG}^\reg$.)

We have constructed in \cite{FG4} an equivalence of categories
\begin{equation}    \label{Gamma Hecke}
\Gamma^{\Hecke}: \on{D}(\Gr_G)^{\Hecke}_\crit\mod^{I^0}
\overset{\sim}\longrightarrow \hg_\crit\mod_\chi^{I^0}
\end{equation}
(this equivalence is canonically defined for each trivialization of
the flat $\cG$-bundle underlying the oper $\chi$, which we will assume
fixed in what follows). Here $\on{D}(\Gr_G)^{\Hecke}_\crit\mod^{I^0}$
is the category of $I^0$-equivariant critically twisted right
D-modules on the affine Grassmannian $\Gr_G$ which satisfy the Hecke
eigensheaf property (see \cite{FG4}, Sect. 1.1, and \secref{def of cat}
below for the precise definition).

In addition, there is another equivalence of categories
\begin{equation}    \label{second eq}
D\bigl(\QCoh((\Fl^{\cG})^{DG})\bigr) \simeq
D\bigl(\on{D}(\Gr_G)^{\Hecke}_\crit\mod\bigr)^{I^0}.
\end{equation}
Here $\Fl^{\cG}$ is the flag variety of $\cG$, and
$D\bigl(\QCoh((\Fl^{\cG})^{DG})\bigr)$ is the derived category of
complexes of quasi-coherent sheaves over the DG-scheme
$$(\Fl^{\cG})^{DG}: =
\Spec\Bigl(\Sym_{\CO_{\Fl^{\cG}}}(\Omega^1(\Fl^{\cG})[1])\Bigr).$$
This DG-scheme can be realized as the derived Cartesian product
$\wt\cg\underset{\cg}\times \on{pt}$, where $\on{pt}\to \cg$
corresponds to the point $0\in \cg$, and $\wt\cg = \{ (x,\check\fb)|x
\in \check\fb \subset \cg \}$ is Grothendieck's alteration.
The equivalence \eqref{second eq} follows from the results of
\cite{ABG}, albeit in a somewhat indirect way.

Combining \eqref{Gamma Hecke} and \eqref{second eq} 
we obtain an equivalence 
\begin{equation} \label{third eq}
D\bigl(\QCoh((\Fl^{\cG})^{DG})\bigr) \simeq D(\hg_\crit\mod_\chi^{I^0}).
\end{equation}
The existence of such an equivalence is a corollary of the Main
Conjecture 6.1.1 of \cite{FG2} (see the Introduction to
\cite{FG4} for more details).

\medskip

We have a natural direct image functor
\begin{equation} \label{DG direct image}
D\bigl(\QCoh(\Fl^{\cG})\bigr) \to
D\bigl(\QCoh((\Fl^{\cG})^{DG})\bigr),
\end{equation}
and the second main objective of this paper is to describe
explicitly its composition with \eqref{third eq}, which is
a functor
\begin{equation} \label{fourth eq}
\sG:D\bigl(\QCoh(\Fl^{\cG})\bigr) \to D(\hg_\crit\mod_\chi^{I^0}).
\end{equation}
(Note that unlike \eqref{third eq}, the functor $\sG$
is not an equivalence.)

The functor $\sG$ of \eqref{fourth eq} turns out to be closely related
to Wakimoto modules, as we shall presently explain. This relationship
confirms the basic property of the equivalence \eqref{third eq} (and
the more general equivalence of the Main Conjecture 6.1.1 of
\cite{FG2}), conjectured in \cite{FG2}, Sect. 6.1.

\ssec{}

Recall that $\Fl^{\cG}$ has a stratification by Schubert cells
$\Fl^{\cG}_w$, where $w$ runs over the Weyl group $W$ of $\cG$.
In \cite{FG2}, Sect. 3.6 we have explained that $\Fl^{\cG}_w$ 
is the parameter space for Wakimoto modules at the critical level
(having a fixed central character $\chi$), and of "highest weight"
$w(\rho)-\rho$. 

More precisely, for every $w$ we have a functor $_w\BW$ from the
category $\QCoh(\Fl^{\cG})_w$ of quasi-coherent sheaves on
$\Fl^{\cG}$, which are {\it set-theoretically} supported on
$\Fl^{\cG}_w$, to the category $\hg_\crit\mod_\chi^{I^0}$.  By
definition, $_w\BW$ is a kind of semi-infinite induction functor that
is embodied by the Wakimoto module construction.

The second main result of this paper, \thmref{full ident},
asserts that, up to
a certain twist, the above functor $_w\BW$ is isomorphic
to the restriction of the functor $\sG$ of \eqref{fourth eq} to
$\QCoh(\Fl^{\cG})_w\subset \QCoh(\Fl^{\cG})$.
In other words, the functor $\sG$ "glues" the functors
$_w\BW$, which are defined for each $w\in W$ separately,
into a single functor. 

\medskip

This fact has a number of interesting representation-theoretic
implications. For example, let $\CT_1,\CT_2$ be two quasi-coherent
sheaves on $\Fl^{\cG}$, supported set-theoretically on two
different Schubert cells $\Fl^{\cG}_{w_1}$ and $\Fl^{\cG}_{w_2}$.
Suppose that $\CT_1\to \CT_2$ is a morphism between them.

Since $\sG$ is a functor, we obtain a homomorphism
$_{w_1}\BW(\CT_1)\to {}_{w_2}\BW(\CT_2)$ of $\hg_\crit$-modules.
The existence of such a homomorphism is not obvious
from the point of view of Wakimoto modules. However, examples
of such homomorphisms have already existed:

If we take $\CT_1$ to be the structure sheaf on the big cell
$\Fl^{\cG}_1\subset \Fl^{\cG}$, and $\CT_2$ to be the 
quasi-coherent sheaf, underlying the D-module
of distributions on a Schubert cell of codimension $1$,
the resulting morphism between the corresponding
Wakimoto modules is the "screening operator" of \cite{FF}.

As an application of the above results, we use
the Cousin-Grothendieck resolution of the structure sheaf of
$\Fl^{\cG}$ to construct a resolution of the
vacuum module $\BV_\crit$ at the critical level in terms of the
Wakimoto modules, corresponding to distributions along the Schubert
cells. The existence of such a resolution was conjectured in \cite{FF}.

\ssec{}

Let us now explain the structure of the part of the paper that
analyzes Wakimoto modules and the functor $\sG$.

In \secref{from to} we study the functor
$$\sE:D\bigl(\QCoh(\Fl^{\cG})\bigr) \to
D\bigl(\on{D}(\Gr_G)^{\Hecke}_\crit\mod\bigr)^{I^0},$$ obtained by
composing the functor \eqref{DG direct image} and the equivalence
\eqref{second eq}. The composition of $\sE$ and the functor
$\Gamma^\Hecke$ of \eqref{Gamma Hecke} is the functor $\sG$ of
\secref{fourth eq}.

The main idea is that the functor $\sE$ is fixed essentially uniquely
by the condition that it respects the action of the group $\cG$ on
both categories, where on the LHS this action comes from the
$\cG$-action on $\Fl^\cG$, and on the RHS this action is as in
\cite{FG4}, Sect. 2.1.

\medskip

In \secref{statement full ident} we collect some basic facts about
Wakimoto modules, as well as some material from \cite{FG4}, and state
our main result, \thmref{full ident}.

\medskip

Sections \ref{principal Wakimoto}, \ref{other Wakimoto} and \ref{more
general} are devoted to the proof of \thmref{full ident}. In
\secref{principal Wakimoto} we treat a particular case of \thmref{full
ident}, namely the one of $w=1$ and the quasi-coherent sheaf
$\CO_{\Fl^{\cG}_1}$, which corresponds to the vacuum Wakimoto module
$\BW_{\crit,0}$. From this particular case, in \secref{other Wakimoto}
we derive the statement of \thmref{full ident} for the quasi-coherent
sheaf $\CO_{\Fl^{\cG}_w}$ for any $w$.

In \secref{more general} we derive the general case of \thmref{full
ident} up to a certain twisting functor that acts from the category
$\QCoh(\Fl^{\cG})_w$ to itself.  To determine this twist, we need to
analyze the action of the {\it renormalized enveloping algebra}
$U^{\ren,\reg}(\hg_\crit)$ introduced in \cite{BD}.

In \secref{action} we review the definition of
$U^{\ren,\reg}(\hg_\crit)$ and make a digression describing the
behavior of certain natural Lie algebroids on the schemes
$\Spec(\fZ^\reg_\fg)$ and $\Op_\cg^\reg$ under the isomorphism
$\on{map}_{\geom}$.

In \secref{renorm Wak} we analyze the interaction of the above
algebroids and the Wakimoto functor, which will allow us to finish the
proof of \thmref{full ident}.

\medskip

Finally, in \secref{resolution} we give another construction of some
Wakimoto modules at the critical level by a procedure of
renormalization.  In particular, the Wakimoto modules appearing in the
Cousin-Grothendieck resolution mentioned above can be obtained in this
way.

\section{Comparison of two morphisms}    \label{comparison}

Let $\fg$ be a simple finite-dimensional Lie algebra, and $G$ the
connected simply-connected algebraic group with the Lie algebra
$\fg$. We shall fix a Borel subgroup $B\subset G$. Denote by $\cG$ the
Langlands dual group of $G$, and let $\cg$ be its Lie algebra. The
group $\cG$ comes equipped with a Borel subgroup $\cB\subset \cG$.

Now let $\hg_\crit$ be the affine Kac-Moody algebra associated
to the critical inner product $\kappa_{\crit}$ and
$\hg_\crit\mod$ the category of discrete $\hg_\crit$-modules (see
\cite{FG2}). Its objects are $\hg_\crit$-modules in which every vector
is annihilated by the Lie subalgebra $\fg \otimes t^n\BC[[t]]$ for
sufficiently large $n$. Let
$$
\BV_\crit :=\on{Ind}^{\hg_\crit}_{\fg[[t]] \oplus \BC}(\BC)
$$
denote the vacuum module, which is an object of the category
$\hg_\crit\mod$. Denote by $\fZ_\fg$ the topological commutative
algebra that is the center of $\hg_\crit\mod$. Let $\fZ^\reg_\fg$
denote its "regular" quotient, i.e., the quotient modulo the
annihilator of $\BV_\crit$. It is known (see \cite{FF,F:wak}) that
\begin{equation}    \label{isom}
\fZ^\reg_\fg = \on{End}_{\hg_\crit}(\BV_\crit) =
(\BV_\crit)^{G[[t]]}.
\end{equation}

Our goal in this section is to reprove a fundamental result of
Beilinson and Drinfeld \cite{BD}, which compares two morphisms from
$\Spec(\fZ^\reg_\fg)$ to the scheme $\Op_{\cg}^\reg$ of opers on the
formal disc $\on{Spec}(\C[[t]])$: one is the Feigin-Frenkel
isomorphism \cite{FF,F:wak}, and the other is a morphism constructed
by geometric means in \cite{BD}.

\ssec{Recollections}

Let $X$ be a smooth complex curve. For a point $x\in X$ we denote by
$\hCO_x$ (resp., $\hCK_x$) the corresponding completed local ring
(resp., field), and write $\D_x=\Spec(\hCO_x)$,
$\D^\times_x=\Spec(\hCK_x)$.

To $\fg$ and an invariant inner product $\kappa$ on $\fg$ one
associates a chiral algebra $\CA_{\fg,\kappa,X}$ on $X$ (see
\cite{CHA}). Denote the chiral algebra corresponding to the critical
inner product $\kappa_{\crit}$ by $\CA_{\fg,\crit,X}$. Consider
the category of chiral $\CA_{\fg,\crit,X}$-modules, supported at a 
given point $x\in X$.

Let us choose a formal coordinate $t$ at $x$. This choice defines a
(tautological) equivalence between the above category and
$\hg_{\crit}\mod$, preserving the forgetful functor to the category of
vector spaces. In particular, the vacuum object, i.e., the fiber
$\CA_{\fg,\crit,x}$ of $\CA_{\fg,\crit,X}$ at $x$, corresponds under
this equivalence to $\BV_\crit$.  For this reason, when we view
$\CA_{\fg,\crit,X}$ as a chiral module over itself we will sometimes
denote it by $\BV_{\crit,X}$.

\ssec{A digression}   \label{digr}

Set $\D=\Spec(\BC[[t]])$ and $\D^\times=\Spec(\BC\ppart)$.  Let
$\Aut(\D)$ be the group scheme of automorphisms of $\BC[[t]]$ that
preserve the maximal ideal. Let $\Der(\D)$ be the Lie algebra of all
derivations of $\BC[[t]]$; it has a basis formed by elements
$L_i=-t^{i+1}\partial_t$, $i \geq -1$. Note that $\Lie(\Aut(\D))$
is a subalgebra of codimension $1$ in $\Der(\D)$; the quotient is
spanned by the image of $L_{-1}$.

\medskip

We recall a general construction assigning to $\Aut(\D)$-modules
(resp., $(\Der(\D),\Aut(\D))$-modules) quasi-coherent $\CO$-modules
(resp., D-modules) on a smooth curve $X$ (see, e.g., \cite{FB} for
more details).

Let $\on{Coord}_X$ denote the $\Aut(\D)$-torsor over $X$, whose fiber
over $x\in X$, denoted by $\on{Coord}_x$, is the scheme of continuous
isomorphisms between $\hCO_x$ and $\BC[[t]]$, preserving the maximal
ideal (equivalently, ``formal coordinates'' at $x$).  Note that the
tautological $\Aut(\D)$-action on $\on{Coord}_X$ extends to a
$\Der(\D)$-action.

\medskip

Let $\sV$ be a representation of $\Aut(\D)$. We form the associated
$\CO$-module on $X$, denoted $\sV_X$, by setting
$$\sV_X := \on{Coord}_X \overset{\Aut(\D)}\times \sV.$$
We will call such $\CO$-modules {\it natural}.
In other words, the fiber of $\sV_X$ at any $x$ is identified with 
$\sV$ for every choice of a formal coordinate $t$ near $x$.
The $\CO$-module $\sV_X$ carries a natural action of the Lie
algebra of vector fields on $X$ by Lie derivatives, that we 
we will denote by $\xi,\sv\mapsto \on{Lie}_\xi(\sv)$. 

For example, for an integer $n$, consider the character
of $\Aut(\D)$ given by $\Aut(\D)\twoheadrightarrow \BG_m
\overset{n}\to \BG_m$, where the first arrow corresponds to
the action of $\Aut(\D)$ on the cotangent space to $\D$ at
the origin. The corresponding natural $\CO$-module
identifies with $\omega_X^n$, where $\omega_X$ is
the canonical line bundle on $X$ and the superscript
"$n$" denotes the $n$th tensor power.

\medskip

If $\sV$ is a Harish-Chandra $(\Der(\D),\Aut(\D))$-module, then
$\sV_X$ acquires a natural left D-module structure. Let $\sv_X$ be a
local section of $\sV_X$ near $x$, and let us choose a formal
coordinate at $x$, thereby identifying $\sV_x$ with $\sV$; let $\sv\in
\sV$ be the value of $\sv_X$ at $x$, and let $\xi$ be a vector field
defined in a neighborhood of $x$. We have the following relation:
\begin{equation} \label{Lie der}
\bigl(\on{Lie}_\xi(\sv_X)-\xi\cdot \sv_X\bigr)_x=\xi(\sv),
\end{equation}
where $\xi\cdot \sv_X$ refers to the left D-module structure on 
$\sV_X$, and $\xi(\sv)$ to the $\Der(D)$-action on $\sV$.

\medskip

The chiral algebra $\CA_{\fg,\crit,X}$ itself is natural in the above
sense (see \cite{FB}, Ch. 19), and most chiral modules over it that we
will consider in this paper will also be natural. Such modules are the
same as $\Aut(\D)$-equivariant $\hg_\crit$-modules. Thus, statements
concerning such modules $\sV_X$ (for all $X$ simultaneously) are
equivalent to $\Aut(\D)$-equivariant statements concerning $\sV$.

\ssec{}

Let $\fz_{\fg,X}$ be the center of $\CA_{\fg,\crit,X}$, regarded as a
commutative D-algebra. A choice of a coordinate as above identifies
the fiber $\fz_{\fg,x}$ at $x$ with $\fZ^\reg_\fg$. The
topological commutative algebra $\wh\fz_{\fg,x}$, corresponding to
$\fz_{\fg,X}$ and $x\in X$ (see \cite{CHA} Sect. 3.6.18), identifies
with $\fZ_\fg$. 

\medskip

Let $\Op_{\cg,X}$ be the D-scheme of $\cg$-opers on $X$ introduced in
\cite{BD} (in this paper we follow the notation of \cite{FG2}). Its
fiber over $x\in X$ is the scheme $\Op_{\cg}(\D_x)$ of
$\cg$-opers on the disc $\D_x$ at $x$. Again, a choice of a coordinate
identifies $\Op_{\cg}(\D_x)$ with $\Op^\reg_\cg$, the scheme of
$\cg$-opers on $\D$, as subschemes in $\Op_\cg(\D^\times_x) \simeq
\Op_\cg(\D^\times)$.

According to the results of \cite{FF,F:wak}, there
exists a canonical isomorphism of D-schemes over $X$:
\begin{equation} \label{FF isomorphism}
\on{map}_{\alg}:\Spec(\fz_{\fg,X}) \overset{\sim}\longrightarrow
\Op_{\cg,X}.
\end{equation}
In particular, we have isomorphisms
$$
\Spec(\fz_{\fg,x}) \simeq \Op_{\cg}(\D_x), \qquad \Spec(\fZ^\reg_\fg)
\simeq \Op^\reg_\cg.
$$
In \secref{property} we will recall an important property of the
isomorphism \eqref{FF isomorphism} that fixes it uniquely.

\medskip

On the other hand, in \cite{BD} a different map
$$
\on{map}_{\geom}: \Spec(\fz_{\fg,X}) \to \Op_{\cg,X}
$$
was constructed, using the affine Grassmannian and the geometric
Satake equivalence. We recall its definition below. We will then prove
the following result:

\begin{thm} \label{two maps coincide}
The morphisms $\on{map}_{\alg}$ and $\on{map}_{\geom}$ coincide.
\end{thm}

This theorem was proved in \cite{BD} by showing that these maps are
compatible with the action of the Lie algebroids
$\isom_{\cG,\fZ_\fg^\reg}$ and $\isom_{\cG,\Op^\reg_\cg}$ (see
\cite{FG2}, Sect. 4), defined on their left and right hand sides,
respectively, and that a map satisfying this property is essentially
unique.

Here we give a different proof in which we use the Wakimoto modules at
the critical level and a key property of $\on{map}_{\alg}$ established in
\cite{FF,F:wak} which we mentioned above.

\ssec{Definition of $\on{map}_{\geom}$}    \label{def geom}

Let us recall the definition of the map $\on{map}_{\geom}$ from
\cite{BD}. Let $\Gr_G = G\ppart/G[[t]]$ be the affine
Grassmannian. This is a strict ind-scheme; in particular
it makes sense to consider the category $\on{D}(\Gr_{G})\mod$
of right D-modules on $\Gr_G$.

We have a natural action of $G[[t]]$ on $\Gr_G$
(given by left multiplication), such that the orbit of every 
finite-dimensional subscheme of $G[[t]]$ is still
finite-dimensional. This insures that the category of
$G[[t]]$-equivariant right D-modules on $\Gr_G$ is
well-defined; we will denote this category by 
$\on{D}(\Gr_{G})\mod^{G[[t]]}$.

The category $\on{D}(\Gr_{G})\mod^{G[[t]]}$ is known to be
semi-simple. Its irreducible objects can be described as follows.  For
$\cla \in \check{\Lambda}^+$ (here, and in the rest of this paper,
$\check\Lambda^+$ denotes the set of dominant coweights of $G$), let
$\Gr^{\cla}_G\subset \Gr_G$ be the $G[[t]]$-orbit of the point
$\cla(t)\in G\ppart$, and let $\Grb^\cla$ be its closure. The
irreducible objects of $\on{D}(\Gr_{G})\mod^{G[[t]]}$ are the
irreducible right D-modules $\IC_{\Grb^\cla}$ corresponding to the
strata $\Grb^\cla$.

By the geometric Satake equivalence (see \cite{MV}), the category
$\on{D}(\Gr_G)\mod^{G[[t]]}$ has a natural structure of tensor
category under convolution, and as such it is equivalent to
$\Rep(\cG)$.  We will denote by $V\mapsto \CF_{V}$ the corresponding
tensor functor $\Rep(\cG) \to \on{D}(\Gr_{G})\mod^{G[[t]]}$. For
$V=V^\cla$, the irreducible $\cG$-representation of highest weight
$\cla$, the corresponding object $\CF_{V^\cla}$ is by definition
isomorphic to $\IC_{\Grb^\cla}$.  For example, if $V=\BC$ is the
trivial representation, the corresponding D-module $\CF_{V}$ is the
$\delta$-function D-module $\delta_{1,\Gr_G}$ at the unit point of
$\Gr_G$.

\medskip

Recall now that over $\Gr_G$ there exists a canonical line bundle
$\CL_{\crit}$; it is $G[[t]]$-equivariant. The action of 
$\g\ppart$ on $\Gr_G$ lifts to an action of
$\ghat_\crit$ on local sections of $\CL_{\crit}$ (with the central
element of $\ghat_\crit$ mapping to the identity). 

Thus, we can consider the category $\on{D}(\Gr_{G})_\crit\mod$ of
$\CL_{\crit}$-twisted right D-modules on $\Gr_G$, and its
$G[[t]]$-equivariant counterpart, denoted
$\on{D}(\Gr_{G})_\crit\mod^{G[[t]]}$. The functor $\CF\mapsto
\CF\otimes \CL_{\crit}$ defines an equivalence between the two
categories (but this equivalence, of course, does not commute with the
functor of global sections). By a slight abuse of notation, for $V\in
\Rep(\cG)$ we will denote by $\CF_V$ also the corresponding
irreducible object of $\on{D}(\Gr_{G})_\crit\mod^{G[[t]]}$.

\medskip

We have the global sections functor
$$\Gamma:\on{D}(\Gr_{G})_\crit\mod\to \hg_\crit\mod.$$ The main result
of \cite{FG1} asserts that this functor is exact and faithful.

In \cite{BD} it is shown (using the results of \cite{FF,F:wak}) that
for $V \in \Rep(\cG)$ there exists a locally free
$(\Der(\D),\Aut(\D))$-equivariant $\fZ^\reg_\fg$-module $\CV_\fZ$,
such that
\begin{equation} \label{BD constr loc}
\Gamma(\Gr_G,\CF_V)\simeq \CV_{\fZ}\underset{\fZ^\reg_\fg}\otimes
\BV_{\crit}.
\end{equation}
Moreover, the assignment $V\mapsto \CV_{\fZ}$ extends to a tensor
functor from $\Rep(\cG)$ to the category of locally free sheaves
$\fZ^\reg_\fg$. By the Tannakian formalism, it defines a $\cG$-torsor
$\CP_{\cG,\fZ}$ on $\Spec(\fZ^\reg_\fg)$, equivariant with respect to
the pair $(\Der(\D),\Aut(\D))$, such that $$\CV_{\fZ} = \CP_{\cG,\fZ}
\overset{\cG}\times V.$$

\medskip

We will now consider the relative versions of the above objects over
the curve $X$. Let $\Gr_{G,X}$ be the global version 
of the affine Grassmannian: this is an ind-scheme over $X$ whose fiber at 
$x \in X$ is $\Gr_{G,x} = G(\hCK_x)/G(\hCO_x)$. Globally, we have:
$$\Gr_{G,X}\simeq \on{Coord}_X\overset{\Aut(\D)}\times \Gr_G,$$ where
$\on{Coord}_X$ is as in \secref{digr}. In addition, $\Gr_{G,X}$ is
endowed with a connection along $X$. By the same construction, the
line bundle $\CL_{\crit}$ gives rise to a line bundle $\CL_{\crit,X}$
on $\Gr_{G,X}$, and the connection on $\Gr_{G,X}$ lifts to
$\CL_{\crit,X}$.

Let $\Jet(G)_X$ denote the group D-scheme of jets of sections of the
group scheme $X \times G$ over $X$; the fiber $\Jet(G)_x$
of $\Jet(G)_X$ at $x \in X$ is by definition $G(\hCO_x)$. 

We will consider the category $\on{D}(\Gr_{G,X})_\crit\mod$ (resp.,
$\on{D}(\Gr_{G,X})_\crit\mod^{\Jet(G)_X}$) of critically twisted
(resp., $\Jet(G)_X$-equivariant) right D-modules on $\Gr_{G,X}$.
Since each object of the form $\CF_V\in
\on{D}(\Gr_{G})_\crit\mod^{G[[t]]}$ for $V\in \Rep(\cG)$ is
$\Aut(\D)$-equivariant, it gives rise to a well-defined object of
$\on{D}(\Gr_{G,X})_\crit\mod^{\Jet(G)_X}$, which we will denote by
$\CF_{V,X}$.

For any object $\CF_X\in \on{D}(\Gr_{G,X})_\crit\mod$, the direct
image of $\CF_X$, considered as a quasi-coherent sheaf, onto $X$, is
naturally a right D-module. Moreover, it has a natural structure of
chiral $\CA_{\fg,\crit,X}$-module.  By a slight abuse of notation we
denote it by $\Gamma(\Gr_{G,X},\CF_X)$.

\medskip

{}From \eqref{BD constr loc} we obtain that, globally over $X$, there
exist $\fz_{\fg,X}$-modules $\CV_{\fZ,X}$, endowed with a connection
along $X$, such that
\begin{equation} \label{BD constr}
\Gamma(\Gr_{G,X},\CF_{V,X})\simeq
\CV_{\fZ,X}\underset{\fz_{\fg,X}}\otimes \BV_{\crit,X},
\end{equation}
and a $\cG$-torsor $\CP_{\cG,\fZ,X}$ on
$\Spec(\fz_{\fg,X})$ with a connection along $X$ such
that
$$\CV_{\fZ,X} = \CP_{\cG,\fZ,X} \overset{\cG}\times V.$$

Furthermore, this $\cG$-torsor is endowed with a reduction to
the Borel subgroup $\cB \subset \cG$, as we shall presently explain.

\medskip

Consider the object $\CF_{V^\cla,X}\in
\on{D}(\Gr_{G,X})_\crit\mod^{\Jet(G)_X}$, corresponding to a dominant
coweight $\cla$ of $G$. By the semi-simplicity of the category
$\on{D}(\Gr_{G})_\crit\mod^{G[[t]]}$, this D-module equals the $0$-th
cohomology of the *-extension of the {\it constant} critically twisted
right D-module on the corresponding $\Jet(G)_X$-orbit
$\Gr^\cla_{G,X}\subset \Gr_{G,X}$, the latter being by definition the
line bundle $\CL_{\crit,X}|_{\Gr^\cla_{G,X}}\otimes
\Omega^{\on{top}}_{\Gr^\cla_{G,X}}$.

Recall now, that according to \cite{BD}, there exists a canonical
isomorphism
\begin{equation} \label{twist cancel}
\CL_{\crit,X}|_{\Gr^\cla_{G,X}}\otimes
\Omega^{\on{top}}_{\Gr^\cla_{G,X}}\simeq 
p^*(\omega_X^{\langle \rho,\cla\rangle})|_{\Gr^\cla_{G,X}},
\end{equation}
where $p: \Gr_{G,X} \to X$ is the canonical projection, and $\omega_X$
is the line bundle of $1$-forms on $X$.  Hence we obtain a map
$\omega_X^{\langle\rho,\cla\rangle} \to
\Gamma(\Gr_{G,X},\CF_{V,X})^{\Jet(G)_X}$.  Using the isomorphism
\eqref{BD constr} and the fact that
$$(\BV_{\crit,X})^{\Jet(G)_X} = \fz_{\fg,X},$$
which follows from \eqref{isom}, we obtain a map
\begin{equation} \label{birth of opers}
\kappa^{\cla}: \fz_{\fg,X}\underset{\CO_X}\otimes 
\omega_X^{\langle \rho,\cla\rangle} \to \CV^\cla_{\fZ,X}.
\end{equation}

Note that $\fz_{\fg,X}\underset{\CO_X}\otimes \omega_X^{\langle
\rho,\cla\rangle}$ is a plain line bundle over $\Spec(\fz_{\fg,X})$,
i.e., it has no connection along $X$.

It is easy to see that the system of maps $\kappa^{\cla}$ satisfies
the Pl\"ucker relations (for the definition, see, e.g., \cite{FG4},
Sect. 4.1) and therefore defines a reduction of
$\CP_{\cG,\fZ,X}$ to $\cB$. Moreover, it is shown in \cite{BD}
that this reduction to $\cB$ satisfies the oper condition relative to
the connection along $X$ on $\CP_{\cG,\fZ,X}$.

\medskip

This defines the desired morphism
$\on{map}_{\geom}:\Spec(\fz_{\fg,X})\to \Op_{\cg,X}$.

\ssec{The defining property of $\on{map}_{\alg}$}    \label{property}

Let us now recall from \cite{FF,F:wak} a property of the morphism
$\on{map}_{\alg}$ that fixes it uniquely (we will follow
the notation of \cite{FG2}). Let $\fH_{\crit,X}$ be the commutative
D-algebra on $X$, introduced in \cite{FG2}, Sect. 10.6. By definition,
$\Spec(\fH_{\crit,X})$ is the D-scheme of induction parameters for
Wakimoto modules. (We will review what this means below.) 

According to the results of \cite{FF,F:wak} (see also \cite{FG2},
Sect. 10.8), we have a canonical map of D-algebras
$$\varphi:\fz_{\fg,X}\to \fH_{\crit,X}.$$

\medskip

Let  $\cH =\cB/[\cB,\cB]$ be the Cartan group of $\cG$. Let 
$\omega^{\rho}_X$ be the $\cH$-torsor over $X$, induced from
the $\BG_m$-torsor, corresponding to the line bundle $\omega_X$
via the co-character $\rho:\BG_m\to \cH$.

Consider the D-scheme $\ConnX$, classifying
connections on the $\cH$-torsor $\omega^{\rho}_X$.

\medskip

Both $\Spec(\fH_{\crit,X})$ and $\ConnX$ are naturally torsors with
respect to the D-scheme classifying $\fh^*=\check\fh$-valued
one-forms on $X$. According to \cite{CHA}, 2.8.17, we have a canonical
isomorphism
\begin{equation} \label{Miura ident alg}
\on{map}^M_{\alg}:\Spec(\fH_{\crit,X})\simeq \ConnX,
\end{equation}
respecting the torsor structure.

Finally, recall from \cite{F:wak}, Sect. 10.3 (see also \cite{FG2},
Sect. 3.3) that we have a canonical map $\on{MT}: \ConnX\to
\Op_{\cg,X}$, called the Miura transformation.

The defining property of $\on{map}_{\alg}$ proved in \cite{FF,F:wak}
is that the diagram
\begin{equation} \label{basic diagram}
\CD
\Spec(\fH_{\crit,X})   @>{\tau\circ \on{map}^M_{\alg}}>>  \ConnX \\
@V{\varphi}VV   @V{\on{MT}}VV   \\
\Spec(\fz_{\fg,X}) @>{\on{map}_{\alg}}>>  \Op_{\cg,X}
\endCD
\end{equation}
is commutative, where $\tau$ denotes the automorphism of $\ConnX$,
induced by the automorphism of $\cH$, given by $\cla\mapsto
-w_0(\cla)$. (Note that this automorphism is well-defined on $\ConnX$,
since $-w_0(\rho)=\rho$.) Moreover, this property determines the
isomorphism $\on{map}_{\alg}$ uniquely.

\medskip

Our strategy of the proof of \thmref{two maps coincide} will be as
follows: we will construct another map
\begin{equation} \label{Miura ident geom}
\on{map}^M_{\geom}:\Spec(\fH_{\crit,X})\simeq \ConnX,
\end{equation}
in the spirit of the above construction of $\on{map}_{\geom}$, for
which the diagram
\begin{equation} \label{geom diagram}
\CD
\Spec(\fH_{\crit,X})   @>{\on{map}^M_{\geom}}>>  \ConnX \\
@V{\varphi}VV   @V{\on{MT}}VV   \\
\Spec(\fz_{\fg,X}) @>{\on{map}_{\geom}}>>  \Op_{\cg,X}
\endCD
\end{equation}
is manifestly commutative. In addition, we will see that the maps
$\on{map}^M_{\geom}$ and $\tau\circ \on{map}^M_{\alg}$ coincide,
thereby implying \thmref{two maps coincide}.

\section{A geometric construction of Miura opers}    \label{geom constr}

The construction of the map $\on{map}^M_{\geom}$ utilizes the Wakimoto
modules introduced in \cite{Wak,FF:usp,FF:si,F:wak}. In this paper we
will mostly follow the notation of \cite{FG2}, Part III.

\ssec{}

Let $\wh\fH_{\crit,x}$ be the topological commutative algebra,
corresponding to the commutative chiral algebra $\fH_{\crit,X}$ and
the point $x\in X$. There exists a canonically defined topological
commutative algebra $\wh\fH_\crit$, acted on by $\Aut(\D)$, such that
every choice of a coordinate $t$ at $x$ defines an isomorphism
$\wh\fH_{\crit,x}\simeq \wh\fH_\crit$.

In fact, $\wh\fH_\crit$ is the completed universal enveloping algebra
of a certain canonical central (and in fact commutative) extension 
$\wh\fh_\crit$ of $\fh\ppart$. In what follows for $\mu\in \fh^*$,
we will denote by $\pi_{\crit,\mu}$ the $\wh\fH_\crit$-module
$$\on{Ind}^{\wh\fh_\crit}_{\fh[[t]]}\bigl(\BC^\mu\bigr),$$ where
$\BC^\mu$ is the $1$-dimensional representation of $\fh[[t]]$,
corresponding to the character $\fh[[t]]\twoheadrightarrow
\fh\overset{\mu}\to \BC$.

Since $\wh\fH_\crit$ is commutative, $\pi_{\crit,\mu}$ is in fact a
quotient algebra, which we will also denote by
$\fH_\crit^{\RS,\mu}$. For $\mu=0$, we have
$\fH_\crit^{\RS,0}=\fH^\reg_\crit$, which identifies with the fiber
$\fH_{\crit,x}$ at $x$. The corresponding global chiral
$\fH_{\crit,X}$-module, denoted $\pi_{\crit,X}$, is, by definition,
the vacuum module $\fH_{\crit,X}$.

\ssec{Wakimoto modules}    \label{wak functor}

In \cite{FG2}, Sect. 11.3 we defined a functor $\BW^{w_0}$ from
the category of chiral $\fH_{\crit,X}$-modules to that of chiral
$\CA_{\fg,\crit,X}$-modules. In the present paper we will only
consider Wakimoto modules "of type $w_0$", so we will omit the
superscript $w_0$ from the notation.

In particular, we obtain a functor 
$$\BW:\wh\fH_\crit\mod\to \hg_\crit\mod.$$

Set $\BW_{\crit,\mu}:=\BW(\pi_{\crit,w_0(\mu)})$. We remark that in
\cite{F:wak} this module was denoted by $W_{\mu,\ka_c}$, and in
\cite{FG2} by $\BW^{w_0}_{\crit,\mu}$. For $\mu=0$ we will denote by
$\BW_{\crit,X}$ the chiral $\CA_{\fg,\crit,X}$-module
$\BW(\pi^0_{\crit,X})$.

\medskip

Assume now that $\mu$ is integral. Then by \cite{FG2}, Sect. 11, the
Wakimoto module $\BW_{\crit,\mu}$ is integrable with respect to the
Iwahori subgroup $I\subset G[[t]]$, i.e., it is naturally an object of
$\hg_\crit\mod^I$, where the latter denotes the category of
$I$-integrable representations of $\hg_\crit$.

\ssec{}

The starting point of our construction of the map $\on{map}^M_{\geom}$
is the following observation.

\medskip

Let $\Fl^{\aff}_G = G\ppart/I$ be the affine flag scheme of $G$. Let
$\on{D}(\Fl^{\aff}_G)_\crit\mod$ be the category of critically twisted
right D-modules on $\Fl_G^{\aff}$ and
$\on{D}(\Fl_G^{\aff})_\kappa\mod^I$ the corresponding $I$-equivariant
category. Recall that to $\CM\in \hg_\crit\mod^{I}$ and $\CF\in
\on{D}(\Gr_G)_\crit\mod^I$ we can associate their convolution
$$\CF\star \CM\in \on{D}(\hg_\crit\mod)^I$$ (see \cite{FG2},
Sect. 22.5 for precise definitions). In particular, for each integral
coweight $\cla$ of $G$ we have an object $j_{\cla,*}$ of
$\on{D}(\Fl_G^{\aff})_\crit\mod^I$ (see \cite{FG2}, Sect. 12.1). It is
the *-extension of the twisted right D-module on $\Fl_G^{\aff}$
corresponding to the constant sheaf on the $I$--orbit $I \cdot \cla(t)
\subset \Fl_G^{\aff}$. The following proposition is a generalization
of \cite{FG2}, Corollary 13.4.2:

\begin{prop}   \label{twist of Wakimoto, local}
For $\cla\in \cLambda^+$, there exists a canonical isomorphism
$$j_{\cla,*}\star \BW_{\crit,\mu}\simeq \BW_{\crit,\mu}
\otimes \fl^{\cla}_{\mu+\rho},$$
respecting the action of $\Aut(\D)$, where $\fl^{\cla}_{\mu+\rho}$
is a $1$-dimensional $\Aut(\D)$-module, corresponding to
the character $$\Aut(\D)\twoheadrightarrow \BG_m\overset
{\langle \mu+\rho,\cla\rangle}\longrightarrow \BG_m.$$
\end{prop}

\begin{proof}

By \cite{FG2}, Corollary 13.4.1, there exists a {\it non-canonical}
isomorphism
$$j_{\cla,*}\star \BW_{\crit,\mu}\simeq \BW_{\crit,\mu}.$$ Hence, by
Proposition 13.1.2 of {\it loc. cit.}, there exists an
$\Aut(\D)$-equivariant line bundle $\CL^\cla_{\fH,\mu}$ over
$\Spec(\fH^{\RS,w_0(\cmu)}_\crit)$ and an isomorphism
$$j_{\cla,*}\star \BW_{\crit,\mu}\simeq \BW_{\crit,\mu}
\underset{\fH_{\crit}^{\RS,w_0(\mu)}}\otimes \CL^\cla_{\fH,\mu}.$$

Consider the embedding $\BG_m\hookrightarrow \Aut(\D)$
given by loop rotations $c\in \BG_m\mapsto (t\mapsto c\cdot t)$.
Since the weights of the resulting action of $\BG$ on
$\fH_{\crit}^{\RS,w_0(\mu)}$ are non-positive, and $\BC\subset
\fH_{\crit}^{\RS,w_0(\mu)}$ is the weight zero subspace, the line
bundle $\CL^\cla_{\fH,\mu}$ is of the form
$\fH_{\crit}^{\RS,w_0(\mu)}\otimes \fl^{\cla}_{\mu+\rho}$,
where $\fl^{\cla}_{\mu+\rho}$ is a $1$-dimensional representation
of $\Aut(\D)$. 

The group $\Aut(\D)$ acts on $\fl^{\cla}_{\mu+\rho}$ via its projection
onto $\BG_m$. We have to show that the corresponding character
of $\BG_m$ is given by the integer $\langle \mu+\rho,\cla\rangle$.
I.e., we have to compare the natural actions of $L_0 = - t
\pa_t\in \Der(\D)$ on $j_{\cla,*}\star \BW_{\crit,\mu}$ and 
$\BW_{\crit,\mu}$, and show that their difference equals the
above integer.

\medskip

Let $\kappa_\hslash$ be a one-parameter deformation of $\kappa_\crit$
away from the critical level, and consider the Wakimoto modules
$\BW_{\hslash,\mu_\hslash}:=\BW_{\kappa_\hslash,\mu_\hslash}$, where
$\mu_\hslash$ is an $\hslash$-family of weights such that
$\mu_\hslash\,\on{mod}\, \hslash=\mu$.  Each of these modules is
endowed with an action of $L_0$ via the Segal-Sugawara
construction. When we shift this action by $\frac{1}{2}\cdot
(\kappa_\hslash-\kappa_\crit)^{-1}(\mu_\hslash,\mu_\hslash+2\cdot\rho)$,
this limiting action at $\kappa_\hslash=\kappa_\crit$ equals the
natural one on $\BW_{\crit,\mu}$.

By \cite{FG2}, Proposition 12.5.1, 
$$\wt{j}_{\kappa_\hslash,\cla}\underset{I^0}\star
\BW_{\hslash,\mu_\hslash}\simeq
\BW_{\hslash,\mu_\hslash-\nu_\hslash},$$ where $\nu_\hslash$ is the
weight such that $\langle \nu_\hslash,\cmu \rangle =
(\kappa_\hslash-\kappa_\crit)(\cla,\cmu)$. Here $I^0$ is the
pro-unipotent radical of $I$ and we use the convolution functor,
denoted by $\underset{I^0}\star$, between $I^0$-equivariant twisted
D-modules on $\wt{\Fl}_G = G\ppart/I^0$ and $I^0$-equivariant
$\hg_{\ka_\hslash}$-modules (note that since the weights on
$\BW_{\hslash,\mu_\hslash}$ are non-integral for generic $\hslash$, it
is not $I$-equivariant, but only $I^0$-equivariant). The corresponding
$I^0$-equivariant D-module $\wt{j}_{\kappa_\hslash,\cla}$ on
$\wt{\Fl}_G$ was defined in \cite{FG2}, Sect. 12.1, where the
connection between the convolution functors
$\wt{j}_{\kappa_\hslash,\cla}\underset{I^0}\star M$ and $j_{\cla,*}
\star M$ is also discussed.

We obtain that in the limit $\kappa_\hslash \to \kappa_\crit$ the resulting
$L_0$ action on $j_{\cla,*}\star \BW_{\crit,\mu}$ differs from that on
$\BW_{\crit,\mu}$ by the limit of
$$\frac{1}{2}\cdot
(\kappa_\hslash-\kappa_\crit)^{-1}(\mu_\hslash-\nu_\hslash,
\mu_\hslash-\nu_\hslash+2\cdot\rho)-\frac{1}{2}\cdot
(\kappa_\hslash-\kappa_\crit)^{-1}(\mu_\hslash,\mu_\hslash+2\cdot\rho),$$
which equals $- \langle \mu+\rho,\cla\rangle$.

\end{proof}

\ssec{Definition of $\on{map}'{}^M_{\geom}$}

Let us consider the chiral $\CA_{\fg,\crit,X}$-module
$\BW_{\crit,X}$; it is obtained from $\BW_{\crit,0}$ by the
procedure described in \secref{digr}. 

In a similar fashion, we can consider a global version 
of $\Fl_G^{\aff}$, and of the twisted D-modules $j_{\cla,*}$;
we will denote the latter by $j_{\cla,*,X}$. 

Consider now the chiral $\CA_{\fg,\crit,X}$-module 
$j_{\cla,*,X}\star \BW_{\crit,X}$. From Proposition
13.1.2 of \cite{FG2}, we obtain that there exists a line bundle
$\CL_{\fH,X}^\cla$ over $\Spec(\fH_{\crit,X})$, endowed with a
connection along $X$, and an isomorphism:

\begin{equation} \label{twist of Wakimoto, global}
j_{\cla,*,X}\star \BW_{\crit,X} \simeq
\CL_{\fH,X}^\cla\underset{\fH_{\crit,X}}\otimes \BW_{\crit,X}.
\end{equation}

The data $\cla\mapsto \CL_{\fH,X}^\cla$ give rise a $\cH$-torsor on
$\Spec(\fH_{\crit,X})$, once we define isomorphisms
\begin{equation}    \label{cH torsor}
\CL_{\fH,X}^\cla\underset{\fH_{\crit,X}}\otimes
\CL_{\fH,X}^\cmu\simeq \CL^{\cla+\cmu}_{\fH,X}
\end{equation}
compatible with the triple tensor products. These are defined by the
requirement that the diagram
$$ \CD j_{\cmu,*,X} \star (j_{\cla,*,X}\star \BW_{\crit,X})
@>>> j_{\cmu,*,X} \star
(\CL_{\fH,X}^\cla\underset{\fH_{\crit,X}}\otimes \BW_{\crit,X}) \\
@VVV @V{\sim}VV \\ j_{\cla+\cmu,*,X}\star \BW_{\crit,X} & &
\CL_{\fH,X}^\cla\underset{\fH_{\crit,X}}\otimes
(j_{\cmu,*,X} \star \BW_{\crit,X}) \\ @VVV @VVV \\
\CL_{\fH,X}^{\cla+\cmu} \underset{\fH_{\crit,X}}\otimes \BW_{\crit,X}
@>>> \CL_{\fH,X}^\cla\underset{\fH_{\crit,X}}\otimes \CL_{\fH,X}^\cmu
\underset{\fH_{\crit,X}}\otimes \BW_{\crit,X} \endCD
$$
be commutative.

\medskip

{}From \propref{twist of Wakimoto, local}, we obtain that as a plain
$\fH_{\crit,X}$-module (i.e., disregarding the connection along $X$),
$\CL_{\fH,X}^\cla$ is isomorphic to $\omega_X^{\langle
\rho,\cla\rangle} \underset{\CO_X}\otimes \fH_{\crit,X} \otimes
\fl^\cla$, where $\fl^\cla$ is a (constant) line.

Thus, we obtain that the line bundles 
$\omega_X^{\langle \rho,\cla\rangle}
\underset{\CO_X}\otimes \fH_{\crit,X}$
on $\Spec(\fH_{\crit,X})$ acquire connections along $X$. These
connections are compatible with the isomorphisms \eqref{cH torsor}.
These data give rise to a morphism of D-schemes
$$\on{map}'{}^M_{\geom}:\Spec(\fH_{\crit,X})\to \ConnX.$$

\medskip

We will establish the following:

\begin{prop}  \label{first comp}
The map $\on{map}'{}^M_{\geom}$ coincides with the map 
$\tau\circ \on{map}^M_{\alg}$.
\end{prop}

A proof of this proposition will be given at the end of this
section. We proceed with the proof of 
\thmref{two maps coincide}.

\ssec{Definition of the map $\on{map}^M_{\geom}$}  
\label{constr geom Miura}

Let $\MOp_{\cG,X}$ be the (non-affine) D-scheme of Miura opers,
as defined in \cite{F:wak} (see also \cite{FG2}, Sect. 3). By definition,
it classifies opers on $X$, endowed with an additional data of
reduction of the corresponding $\cG$-bundle to the subgroup 
$\cB^-\subset \cG$, compatible with the connection. (Here
$\cB^-\subset \cG$ is an opposite Borel subgroup that we
choose once and for all.)

Let $\MOp_{\cG,\gen,X}\subset \MOp_{\cG,X}$ be the open subscheme,
classifying generic Miura opers, as defined in \cite{F:wak},
Sect. 10.3 (see also \cite{FG2}, Sect. 3.3). The genericity condition
is that the above reduction to $\cB^-$ is at all points of $X$ in the
generic position with respect to the reduction to $\cB$, given by the
oper structure.

\medskip

By {\it loc. cit.} there exists a canonical isomorphism of D-schemes
\begin{equation} \label{Miura as Con}
\MOp_{\cG,\gen,X}\simeq \ConnX.
\end{equation}

\medskip

Thus, constructing a morphism  $\on{map}^M_{\geom}:\Spec(\fH_{\crit,X})\to 
\ConnX$ is equivalent to constructing a morphism
\begin{equation} \label{map geom M}
\Spec(\fH_{\crit,X})\to \MOp_{\cg,\gen,X}.
\end{equation}

By definition, to define a map as in \secref{map geom M}, 
we need to construct a $\cG$-bundle
$\CP_{\cG,\fH,X}$ on $\Spec(\fH_{\crit,X})$, equipped with a
connection along $X$, endowed with an oper structure, and a reduction
to $\cB^-$, which is compatible with connection, and which is in
generic relative position with respect to the reduction to $\cB$,
given by the oper structure. 

We define the $\cG$-bundle $\CP_{\cG,\fH,X}$ on
$\Spec(\fH_{\crit,X})$ with its oper structure as the pull-back from
$\CP_{\cG,\fZ,X}$ on $\Spec(\fz_{\fg,X})$ via the map $\varphi$. To 
define a reduction to $\cB^-$ on $\CP_{\cG,\fH,X}$ we proceed as follows:

\medskip

We have the canonical embedding of chiral $\hg_\crit$-modules
$$\phi:\BV_{\crit}\to \BW_{\crit,0}.$$ For $\cla\in \cLambda^+$
consider the $I$-orbit $\Gr^{\cla,I}_G = I \cdot \cla(t)$ in
$\Gr_G$. Let $j_{\cla,\Gr_G,*}$ be the *-extension of the right
D-module on this orbit corresponding to the constant sheaf. Observe
that $\Gr^{\cla,I}_G$ is open and dense in the $G[[t]]$-orbit
$\Gr^\cla_G = G[[t]] \cdot \cla(t)$, so that we have a map
$\CF_{V^\cla} \to j_{\cla,\Gr_G,*}$. Furthermore, under the projection
$\Fl_G^{\aff} \to \Gr_G$ the $I$-orbit $I \cdot \cla(t) \subset
\Fl_G^{\aff}$ is mapped to $\Gr^{\cla,I}_G$ one-to-one. Hence, we have
an isomorphism
$$
j_{\cla,*}\star \BV_{\crit} \simeq j_{\cla,\Gr_G,*}
\underset{G[[t]]}\star \BV_{\crit},
$$
where $\underset{G[[t]]}\star$ denotes the convolution functor between
the category of twisted D-modules on $\Gr_G$ and the category of
$G[[t]]$-equivariant $\hg_\crit$-modules (to which $\BV_{\crit}$
belongs). Consider the composition:
\begin{equation}  \label{imp comp} 
\Gamma(\Gr_{G},\CF_{V})\simeq \CF_{V^\cla}
\underset{G[[t]]}\star \BV_{\crit}\to j_{\cla,*}\star
\BV_{\crit} \overset{j_{\cla,*}\star \phi}\longrightarrow
j_{\cla,*}\star \BW_{\crit}.
\end{equation}

Hence, using \eqref{BD constr loc} and \propref{twist of Wakimoto, local},
we obtain a map
\begin{equation} \label{bos lambda, local}
\CV^\cla_{\fZ}\underset{\fZ^\reg_\fg}\otimes \BV_{\crit}\to
\CL_{\fH}^\cla\otimes \BW_{\crit,0},
\end{equation}
where $\CL^\cla_{\fH}:=\CL_{\fH,0}^\cla\simeq \fH^\reg_\crit\otimes
\fl_\rho^\cla$ is an in the proof of \propref{twist of Wakimoto,
local}

However, by \cite{FG2}, Proposition 10.7.1, the map
$$\BV_{\crit}\underset{\fZ^\reg_\fg}\otimes \fH^\reg_{\crit}\to
\BW_{\crit,0}$$ induces an isomorphism of the spaces of
$G[[t]]$-invariants
$$\fH^\reg_{\crit}\to (\BW_{\crit,0})^{G[[t]]}.$$

Therefore, the map \eqref{bos lambda, local} gives rise to a map of
$\fH_{\crit,X}$-modules,  compatible with the connection:
\begin{equation} \label{des Miura, local}
\kappa^{-,\cla}_{\fH}: \CV^\cla_{\fZ}\underset{\fZ^\reg_\fg}\otimes
\fH^\reg_{\crit}\to \CL_\fH^\cla.
\end{equation}

\medskip

The same construction can be performed over the base $X$, and we 
obtain a map


\begin{equation} \label{bos lambda}
\CV^\cla_{\fZ,X}\underset{\fz_{\fg,X}}\otimes \BV_{\crit,X}\to
\CL_{\fH,X}^\cla\underset{\fH_{\crit,X}}\otimes \BW_{\crit,X},
\end{equation}
and hence a map
\begin{equation} \label{des Miura}
\kappa^{-,\cla}_{\fH}: \CV^\cla_{\fZ,X}\underset{\fz_{\fg,X}}\otimes
\fH_{\crit,X}\to \CL_{\fH,X}^\cla.
\end{equation}

The maps $\kappa^{-,\cla}_{\fH}$ for $\cla\in \cLambda^+$ satisfy the
Pl\"ucker equations since the diagrams
$$
\CD
\CF_{V^\cla}\star \CF_{V^\cmu} @>>> j_{\cla,\Gr_G,*}
\underset{G[[t]]}\star \CF_{V^\cmu}
@>{\sim}>> j_{\cla,*} \star \CF_{V^\cmu} \\
@VVV  & & @VVV   \\
\CF_{V^{\cla+\cmu}} @>>> j_{\cla+\cmu,\Gr_G,*}
@>{\sim}>>  j_{\cla,*} \star j_{\cmu,\Gr_G,*}
\endCD
$$
are clearly commutative. Therefore the system of maps
$\{\kappa^{-,\cla}_{\fH}\}$ defines a reduction of
$\CP_{\cG,\fH,X}$ to $\cB^-$. By construction, this reduction
is horizontal with respect to the connection.

In order to show that this $\cB^-$-reduction is in generic relative
position with the oper $\cB$-reduction it suffices to prove the
following:

\begin{prop}   \label{crucial comp}
The composed arrow
$$\omega_X^{\langle \rho,\cla\rangle}\underset{\CO_X}\otimes
\fH_{\crit,X} \overset{\kappa^{\cla}}\to \CV^\cla_{\fZ,X}
\underset{\fz_{\fg,X}}\otimes
\fH_{\crit,X}\overset{\kappa^{-,\cla}_{\fH}}\to \CL_{\fH,X}^\cla \simeq
\omega_X^{\langle \rho,\cla\rangle}\underset{\CO_X}\otimes
\fH_{\crit,X} \otimes \fl^\cla$$ is an isomorphism, which is induced
by a trivialization of the $\cH$-torsor $\{\cla\mapsto \fl^\cla\}$.
\end{prop}

Let us assume this proposition and finish the proof of 
\thmref{two maps coincide}.  

\medskip

First, \propref{crucial comp} implies that the data
$(\CP_{\cG,\fH,X}, \kappa^{\cla},\kappa^{-,\cla}_{\fH})$ define
a generic Miura oper over $\Spec(\fH_{\crit,X})$. This gives rise to
a map $\Spec(\fH_{\crit,X})\to \MOp_{\cg,\gen,X}$, which we compose
with the identification \eqref{Miura as Con} to produce the
sought-after $\on{map}^M_{\geom}$,

The diagram \eqref{geom diagram} is commutative,
since the map $\on{MT}:\ConnX\to \Op_{\cg,X}$ is by definition
the composition of \eqref{Miura as Con} and the tautological
projection $\MOp_{\cg,\gen,X}\to \Op_{\cg,X}$.

\medskip

In view of \propref{first comp}, it remains to show that the map
$\on{map}^M_{\geom}$, constructed above coincides with
$\on{map}'{}^M_{\geom}$.

\medskip

To see that, we recall from \cite{F:wak}, Sect. 10.3 (see also
\cite{FG2}, Sect. 3.3) that the map $\MOp_{\cg,\gen,X}
\overset{\sim}\longrightarrow \ConnX$ of \eqref{Miura as Con}
is defined as follows:

Given a Miura oper, the genericity assumption implies that
the $\cH$-bundle, induced from the $\cB^-$-bundle, is
isomorphic to $\omega_X^\rho$, and hence, the latter
acquires a connection.

Therefore, the connection along $X$ on every line
bundle $\omega_X^{\langle \rho,\cla\rangle}\underset{\CO_X}\otimes 
\fH_{\crit,X}$, corresponding to $\on{map}^M_{\geom}$,
equals the one arising from the composed isomorphism of
\propref{crucial comp}. The latter equals, by definition, to
the connection on $\omega_X^{\langle \rho,\cla\rangle}\underset{\CO_X}
\otimes \fH_{\crit,X}$, corresponding to the map $\on{map}'{}^M_{\geom}$. 

This completes the proof of \thmref{two maps coincide}.

\ssec{Proof of \propref{crucial comp}}

It is enough to show that for a fixed point $x\in X$ the composition
\begin{equation}   \label{crucial comp, local}
\omega_x^{\langle \rho,\cla\rangle}\to \Gamma(\Gr_G,\CF_{V^\cla})\to
j_{\cla,*} \star \BW_{\crit,0}\simeq 
\omega_x^{\langle \rho,\cla\rangle}\otimes \fl^\cla\otimes \BW_{\crit,0}
\end{equation}
is an isomorphism onto the sub-space
$\omega_x^{\langle \rho,\cla\rangle}\otimes \fl^\cla\subset \BW_{\crit,0}$,
corresponding to the generating vector of $\BW_{\crit,0}$.

In fact, we claim that it is enough to show that the composition in
\eqref{crucial comp, local} is non-zero. Indeed, according to
\cite{BD}, Proposition 8.1.5, the image of $\omega_x^{\langle
\rho,\cla\rangle}$ in $\Gamma(\Gr_G,\CF_{V^\cla})$ equals the subspace
on which the operator $L_0$ acts with the eigenvalue $- \langle
\rho,\cla\rangle$ (all other eigenvalues of $L_0$ being strictly
greater).  Similarly, the generating vector of $\BW_{\crit,0}$ spans
the subspace corresponding to the zero eigenvalue of $L_0$.

\medskip

To show the non-vanishing of \eqref{crucial comp, local} we proceed as
follows.  Let us apply to the two sides of the morphism
$$\Gamma(\Gr_G,\CF_{V^\cla})\to j_{\cla,*} \star
\BW_{\crit,0}$$ the semi-infinite cohomology functor
$H^\semiinf\bigl(\fn\ppart,\fn[[t]], ?\otimes \Psi_0\bigr)$ (see
\cite{FG2}, Sect. 18).

As in \cite{FG2}, Sect. 18.3, we have a commutative diagram
$$ \CD H^\semiinf\bigl(\fn\ppart,\fn[[t]],
\Gamma(\Gr_G,\CF_{V^\cla})\otimes \Psi_0\bigr) @>>>
H^\semiinf\bigl(\fn\ppart,\fn[[t]], (j_{\cla,*} \star
\BW_{\crit,0})\otimes \Psi_0\bigr) \\ @V{\sim}VV @V{\sim}VV \\
H^\semiinf\bigl(\fn\ppart,\fn[[t]], \BV_\crit\otimes \Psi_\cla\bigr)
@>{\phi}>> H^\semiinf\bigl(\fn\ppart,\fn[[t]], \BW_{\crit,0}\otimes
\Psi_\cla\bigr).  \endCD
$$ It is easy to see that under the left vertical map the image of
$$\omega_x^{\langle \rho,\cla\rangle}\to
\Gamma(\Gr_G,\CF_{V^\cla})^{G[[t]]}\to
H^\semiinf\bigl(\fn\ppart,\fn[[t]], \Gamma(\Gr_G,\CF_{V^\cla})\otimes
\Psi_0\bigr)$$ is mapped to the one-dimensional vector space spanned
by the image of the canonical generator of $\BV_\crit$ in
$$\BV_{\crit}^{G[[t]]}\to H^\semiinf\bigl(\fn\ppart,\fn[[t]],
\BV_\crit\otimes \Psi_\cla\bigr).$$ Under the bottom horizontal map
the latter goes to the image of the the canonical generator of
$\BW_{\crit,0}$ in
$$\fH^\reg_\crit \simeq H^\semiinf\bigl(\fn\ppart,\fn[[t]],
\BW_{\crit,0}\otimes \Psi_\cla\bigr),$$ and, in particular, it is
non-zero.

\ssec{Proof of \propref{first comp}}

We claim that it is enough to show that the map
$$\tau\circ \on{map}'{}^M_{\geom}:\Spec(\fH_{\crit,X})\to \ConnX$$
respects the torsor structure with respect to the D-scheme of
$\fh^*$-valued one-forms on $X$. 

Indeed, the morphism $\on{map}^M_{\alg}$ has this property and is an
isomorphism, by definition.  Therefore the difference of the two maps
$\tau\circ \on{map}'{}^M_{\geom}-\on{map}^M_{\alg}$ can be regarded as
an $\fh^*$-valued one-form, canonically attached to the curve $X$. In
particular, this one-form would be invariant under Lie derivatives by
all vector fields, which implies that it is equal to zero.

\medskip

To verify the required property of the map $\on{map}'{}^M_{\geom}$
with respect to the torsor structure we need to check the following:

Consider the action of $\Der(\D)$ on
$j_{\cla,*} \star \BW_{\crit,0}\simeq 
\fl^\cla_{\rho}\otimes \BW_{\crit,0}$
at some point $x\in X$, given by the formula \eqref{Lie der}.
This action preserves the subspace of $\fg[[t]]$-invariant vectors. 

Recall that $L_{-1}\in \Der(\D)$ denotes the element $-\partial_t$.
Let $\bw$ be an element from the $1$-dimensional vector space
$\fl^\cla_{\rho}\subset j_{\cla,*} \star \BW_{\crit,0}$.

We need to show that
$$L_{-1}\cdot \bw \, \on{mod}\, \BC\cdot \bw=(w_0(\cla)\otimes
t^{-1})\cdot \bw,$$ where we identify $(\fh\otimes t^{-1}\cdot
\BC[[t]])/\BC[[t]]$ with a sub-space of $\wh\fh_\crit/\fh[[t]]\subset
\fH^\reg_\fh$. The latter is a straightforward calculation, performed
below.

\medskip

We may realize the action of $\on{Der}(\D)$ on $j_{\cla,*} \star
\BW_{\crit,0}$ as the limit $\kappa_\hslash\to \kappa_\crit$ of its
actions on
$$\wt{j}_{\kappa_\hslash,\cla,*}\underset{I^0}\star \BW_{\hslash,0}
\simeq \BW_{\hslash,-\nu_\hslash},$$ given by the Segal-Sugawara
construction, where $\nu_\hslash$ is as in the proof of \propref{twist
of Wakimoto, local}.

Let us note also that in order to show that two elements of
$\fH^{\reg}_{\crit}\simeq \BW_{\crit,0}^{G[[t]]}$ are equal, it is
sufficient to analyze their images in
$$\fH^\reg_{\crit}\simeq H^\semiinf\bigl(\fn\ppart,\fn[[t]],
\BW_{\crit,0}\bigr).$$ Let us regard
$H^\semiinf\bigl(\fn\ppart,\fn[[t]], \BW_{\hslash,-\nu_\hslash}\bigr)$
as acted on by the chiral algebra $\fH'_{\kappa_\hslash}$ in the
notation of \cite{FG2}, Sect. 10.2. Let $\wh\fh'_\hslash$ be the
corresponding central extension of $\fh\ppart$.

We have a canonical isomorphism
$$H^\semiinf\bigl(\fn\ppart,\fn[[t]], \BW_{\hslash,-\nu_\hslash}\bigr)
\simeq \pi'_{\hslash,-\nu_\hslash}:=
\on{Ind}^{\wh\fh'_\hslash}_{\fh[[t]]}(\BC^{-\nu_\hslash}),$$
compatible with the $\on{Der}(\D)$-action.

Note that for $\kappa_\hslash=\kappa_\crit$ the action of
$\fH^\reg_\crit$ on $H^\semiinf\bigl(\fn\ppart,\fn[[t]],
\BW_{\crit,0}\bigr),$ coming from its action on $\BW_{\crit,0}$,
coincides with the one given by the isomorphism of chiral algebras
$\fH'_{\crit,X}\simeq \fH_{\crit,X}$ of \cite{FG2}, 10.6, composed
with the automorphism, induced by $\tau$.

\medskip

To summarize, we need to compute the action of $L_{-1}$, given by the
Segal-Sugawara construction, on the $\wh\fh'_\hslash$-module
$\pi'_{\hslash,-\nu_\hslash}$ in the limit $\kappa_\hslash\to
\kappa_\crit$.

If $\bv_\hslash$ denotes the canonical generator in
$\pi_{\hslash,-\nu_\hslash}$, we need to show that
$$L_{-1}\cdot \bv_\hslash \,\on{mod}\, \BC\cdot \bv_\hslash = -
(\cla\otimes t^{-1})\cdot \bv_\hslash \, \on{mod}\, \hslash.$$

However, this coincides with the formula for the action of the
operator $L_{-1}$ on the Heisenberg algebra $\wh\fh'_\hslash$ (see
\cite{F:wak}, Sect. 5.5). This completes the proof of \propref{first
comp}.

\section{From $\OO$-modules on $\Fl^\cG$ to D-modules on $\Gr_G$}
\label{from to}

\ssec{}
Let $\Fl^\cG$ be the flag variety of the group $\cG$, thought of as
the quotient $\cG/\cB^-$. Consider the category
$\QCoh\bigl(\Fl^{\cG}\bigr)$ of quasi-coherent sheaves on $\Fl^\cG$,
and the corresponding derived category
$D\bigl(\QCoh\bigl(\Fl^{\cG}\bigr)\bigr)$.

Consider the category $\on{D}(\Gr_G)^{\Hecke}_\crit\mod$ of {\it Hecke
eigensheaves} on $\Gr_G$, introduced in \cite{FG4}, Sect. 2.1 (the
definition will be recalled below).

In this section we will study the functor
$$\sE:D^+\bigl(\QCoh\bigl(\Fl^{\cG}\bigr)\bigr)\to
D^+\bigl(\on{D}(\Gr_G)^{\Hecke}_\crit\mod\bigl),$$ obtained by
composing the equivalence \eqref{second eq} of \cite{ABG} and the
direct image functor \eqref{DG direct image}, as was explained in the
Introduction.  The functor $\sE$ will be left-exact.

We should remark, however, that in this paper we will not formally
rely on the results of \cite{ABG}. We will construct {\it a} functor
$\sE$ "from scratch", with {\it loc. cit.}  serving as a guide.

We remark that the contents of this section have a significant
intersection with Sect. 3.2.13 of \cite{ABBGM} and Sect. 3 of
\cite{FG4}.

\ssec{}    \label{def of cat}

Let us first recall the definition of the category
$\on{D}(\Gr_G)^{\Hecke}_\crit\mod$.\footnote{We remark that in
\cite{FG2} we used the category $\on{D}(\Gr_G)^{\Hecke}_\crit\mod$
when the group $G$ was of the adjoint type. It is easy to see,
however, that the two categories are equivalent (see \cite{AG}).}

Its objects are the data of
$(\CF,\alpha_V,\,\, \forall\,\, V\in \Rep(\cG))$, where $\CF$ is an
object of $\on{D}(\Gr_G)_\crit\mod$, and $\alpha_V, V\in \Rep(\cG))$,
are isomorphisms of D-modules
$$\al_V: \CF\underset{G[[t]]}\star \CF_V \overset{\sim}\longrightarrow
\uV\otimes \CF,$$ 
(here $\uV$ denotes the vector space underlying the representation $V$)
such that the following two conditions are satisfied:

\begin{itemize}

\item
If $V$ is the trivial representations $\BC$, then the morphism
$\alpha_V$ is the identity map.

\item
For $V,W\in \Rep(\cG)$ and $U:=V\otimes W$, the diagram
$$
\CD
(\CF\underset{G[[t]]}\star \CF_V)\underset{G[[t]]}\star \CF_W  @>{\sim}>>  
\CF\underset{G[[t]]}\star \CF_U \\
@V{\alpha_V\underset{G[[t]]}\star \on{id}_{\CF_W}}VV  @V{\alpha_U}VV \\
(\uV\otimes \CF)\underset{G[[t]]}\star \CF_W  & &
\uU\otimes \CF \\
@V{\sim}VV   @V{\sim}VV  \\
\uV\otimes (\CF\underset{G[[t]]}\star \CF_W) 
@>{\on{id}_{\uV}\otimes \alpha_W}>>  
\uV \otimes \uW \otimes\CF
\endCD
$$ is commutative.

\end{itemize}

Morphisms in this category between $(\CF,\alpha_V)$ and
$(\CF',\alpha'_V)$ are maps of D-modules $\phi:\CF\to \CF'$, and such
that
$$(\on{id}_{\uV}\otimes \phi)\circ \alpha_V=
\alpha'_V\circ (\phi\star \on{id}_{\CF_V}).$$

Note that this category carries a canonical action of the group
$\cG$. Indeed, for $\bg\in \cG$ and an object $(\CF,\alpha_V)\in
\on{D}(\Gr_G)^{\Hecke}_\crit\mod$ we define a new object as
$(\CF,\bg|_{\uV}\cdot \alpha_V)$, where $\bg|_{\uV}$ denotes the
automorphism by means of which $\bg$ acts on the vector space
underlying $V$.

\ssec{}     \label{gen framework of functor}

Consider now the following general situation. Let $\CC$ be an abelian
category equipped with an action of $\cG$; let $\on{act}^*:\CC\to
\QCoh(\cG)\otimes \CC$ be the corresponding action functor (we refer
the reader to \cite{FG2}, Sect. 20.1 where the corresponding notions
are discussed in detail).

We will describe a general framework, in which one constructs a functor
$$\sF: D^+(\QCoh(\Fl^\cG))\to D^+(\CC),$$ compatible with the $\cG$-actions.

\medskip

Suppose we are given a collection of $\cG$-equivariant objects
$\CM_\cla\in \CC,\,\,\cla\in \cLambda^+$, and a collection of
$\cG$-equivariant morphisms
$$\beta^{\cmu,\cla}:V^\cmu\otimes \CM_\cla\to \CM_{\cla+\cmu},$$
defined for $\cmu\in \cLambda^+$,
where the LHS is endowed with a $\cG$-equivariant structure
via the diagonal action. We will assume that for $\cmu_1,\cmu_2 \in
\cla\in \cLambda^+$ the diagram
\begin{equation}   \label{la mu and mu}
\CD
V^{\cmu_1}\otimes V^{\cmu_2}\otimes\CM_\cla 
@>{\on{id}_{V^{\cmu_1}}\otimes \beta^{\cmu_2,\cla}}>> 
V^{\cmu_1}\otimes \CM_{\cla+\cmu_2}  \\
@VVV     @V{\beta^{\cmu+1,\cla+\cmu_2}}VV \\
V^{\cmu_1+\cmu_2}\otimes\CM_\cla @>{\beta^{\cmu_1+\cmu_2,\cla}}>>
\CM_{\cla+\cmu_1+\cmu_2}
\endCD
\end{equation}
is commutative.

\medskip

One example of the above situation arises when $\CC=\QCoh(\Fl^\cG)$
and $\CM_\cla$ is taken to be $\CL^\cla_{\Fl^\cG}$, i.e., the
$\cG$-equivariant line bundle on $\Fl^\cG$, attached to the weight
$\cla$.  Our normalization is such that
$\Gamma(\Fl^\cG,\CL^\cla_{\Fl^\cG})=V^\cla$.

We shall now perform the following procedure in this example that will
allow us to expresses a large class of objects of $\QCoh(\Fl^\cG)$ in
terms of the line bundles $\CL^\cla_{\Fl^\cG}$:

\medskip

Let $\CT$ be a quasi-coherent sheaf on $\Fl^\cG$ obtained as the
direct image of a quasi-coherent sheaf on an affine locally closed
subscheme of $\Fl^\cG$.

For $\cmu\in \cLambda^+$ consider the direct sums
$$\CT^+:=\underset{\cla\in \cLambda^+}\oplus \, 
\Gamma(\Fl^\cG,\CT\otimes \CL^{-\cla}_{\Fl^\cG})\otimes
\CL^\cla_{\Fl^\cG}$$
and
$$\CT^{++}:=\underset{\cmu,\cla\in \cLambda^+}\oplus \, 
\uV^\cmu\otimes \Gamma(\Fl^\cG,\CT\otimes \CL^{-\cla-\cmu}_{\Fl^\cG})\otimes
\CL^\cla_{\Fl^\cG}.$$

Note that we have a canonical map 
\begin{equation} \label{T quot of dir sum}
\CT^+\to \CT.
\end{equation}

In addition, we have two maps
\begin{equation} \label{two maps to co-eq}
\CT^{++}\rightrightarrows \CT^+.
\end{equation}
The first map comes from the canonical map 
\begin{equation} \label{beta basic}
\uV^\cmu \otimes
\CL^\cla_{\Fl^\cG}\to \CL^{\cla+\cmu}_{\Fl^\cG}, 
\end{equation}
tensored with the identity on each $\Gamma(\Fl^\cG,\CT\otimes
\CL^{-\cla-\cmu}_{\Fl^\cG})$.  The second map comes from the map
$$\uV^\cmu\otimes \Gamma(\Fl^\cG,\CT\otimes \CL^{-\cla-\cmu}_{\Fl^\cG})\to
\Gamma(\Fl^\cG,\CT\otimes \CL^{-\cla}_{\Fl^\cG}),$$
tensored with the identity map on $\CL^\cla_{\Fl^\cG}$.

The proof of the following result is straightforward.

\begin{lem}
The map \eqref{T quot of dir sum} identifies $\CT$ with the
co-equalizer of the map \eqref{two maps to co-eq}.
\end{lem}

We will now perform the same construction in a category $\CC$ equipped
with the above structures. For $\CT$ as above consider two objects of
$\CC$ defined as:
$$\sF(\CT)^+:=\underset{\cla\in \cLambda^+}\oplus \, 
\Gamma(\Fl^\cG,\CT\otimes \CL^{-\cla}_{\Fl^\cG})\otimes \CM_\cla$$
and
$$\sF(\CT)^{++}:=\underset{\cmu,\cla\in \cLambda^+}\oplus \, 
\uV^\cmu\otimes \Gamma(\Fl^\cG,\CT\otimes \CL^{-\cla-\cmu}_{\Fl^\cG})\otimes
\CM_\cla.$$

The maps $\beta^{\cmu,\cla}$ and \eqref{beta basic} give rise to two 
morphisms:
\begin{equation}  \label{two arr gen}
\sF(\CT)^{++}\rightrightarrows
\sF(\CT)^+.
\end{equation}

Define $\sF(\CT)$ as the co-equalizer of the above maps. 
It is clear that if $\CT$ is $\cG$ (or $\cB$)-equivariant, then 
so is $\sF(\CT)$. Set 
\begin{equation}  \label{M w0}
\CM_{w_0}:=\sF(\BC_{w_0}), 
\end{equation}
where $\BC_{w_0}$
is the sky-scraper at the $\cB$-invariant point $w_0\in \Fl^\cG$. By the
above, $\CM_{w_0}$ is $\cB$-equivariant.

\medskip

To extend the functor $\sF$ to arbitrary quasi-coherent sheaves on $\Fl^\cG$,
we need to make a digression and discuss the following general construction.

\ssec{Coherent convolutions}    \label{coherent convolutions}

Let $\CC$ be as above, $\CM\in \CC$ an object equivariant with respect
to $\cB\subset \cG$, and let $\CT$ be a quasi-coherent sheaf on
$\Fl^\cG$. Let us identify $\Fl^\cG$ with the quotient $\cG/\cB$, by
the action of $\cG$ on the point $w_0\in \Fl^\cG:=\cG/\cB^-$
stabilized by $\cB$. Let $\wt{\CT}$ be the $\cB$-equivariant
quasi-coherent sheaf on $\cG$, corresponding to $\CT$.

Consider the object $\wt\CT \underset{\Fun(\cG)}\otimes
\on{act}^*(\CM) \in \QCoh(\cG)\otimes \CC$. Regarded as an object of
$\CC$, it carries an action of $\cB$ by automorphisms, compatible with
the action of the latter on $\Fun(\cG)$ by right translations.

\medskip

We define the convolution $\CT\ast \CM\in \CC$ as
\begin{equation} \label{defn conv}
\CT\ast \CM = \Bigl(\wt{\CT}\underset{\Fun(\cG)}\otimes
\on{act}^*(\CM)\Bigr)^{\cB}.
\end{equation}

Keeping $\CM\in \CC$ fixed, we can regard $\CT\mapsto \CT\ast \CM$ as
a functor $\QCoh(\cG)\to \CC$. This functor is evidently left exact, and
we will denote by $\CT\overset{R}\ast \CM$ its right derived functor.

It is easy to see that if $\CT$ is isomorphic to the direct image of a
quasi-coherent sheaf on an affine locally closed subscheme of
$\Fl^\cG$, then we have an isomorphism $\CT\overset{R}\ast \CM\simeq
\CT\ast \CM$.


\medskip

Let us consider some examples of the above situation. Take first
$\CC=\QCoh(\Fl^\cG)$ with the natural action of $\cG$, and let
$\CM$ be $\BC_{w_0}$--the skyscraper at the the point $w_0\in \Fl^\cG$.
Since this point is fixed by $\cB$, the object $\BC_{w_0}$ is
$\cB$-equivariant.

It is clear that in this case the functor $\CT\mapsto \CT \ast
\BC_{w_0}$ is tautologically isomorphic to the identity functor.

\medskip

Suppose now that $\sF$ is a left-exact functor $\QCoh(\Fl^\cG)\to \CC$,
compatible with the actions of $\cG$ in the evident sense.
Set $\CM_{w_0}:=\sF(\BC_{w_0})$. This object is $\cB$-equivariant,
by transport of structure.

We then obtain a functorial isomorphism:
\begin{equation} \label{funct as conv}
\sF(\CT)\simeq \CT\ast \CM_{w_0}.
\end{equation}

\medskip

Finally, let us consider the case $\CC=\on{D}(\Gr_G)^{\Hecke}_\crit\mod$,
and let us write down the convolution functors more explicitly.

Namely, for an object $\CM:=(\CF,\alpha_V)$ to be $\cB$-equivariant means
that $\CF$, as an object of $\on{D}(\Gr_G)_\crit\mod$, is endowed with an 
algebraic action of $\cB$ by automorphisms, and the morphisms
$$\alpha_V:\CF\underset{G[[t]]}\star \CF_V\simeq \uV\otimes \CF$$
respect the $\cB$-action, where on the LHS the action comes by
functoriality from the $\cB$-action on $\CF$, and on the RHS it is
diagonal with respect to the natural $\cB$-action on $\uV$.

Then the object of $\on{D}(\Gr_G)_\crit\mod$, underlying $\CT\ast \CM$
is given by $\Bigl(\wt\CT\otimes \CF\Bigr)^{\cB}$. (It is also easy to
write down the remaining data of $\CT\ast \CM$ making it into an
object of $\on{D}(\Gr_G)^{\Hecke}_\crit\mod$.)

\ssec{}   \label{gen framework, contd}

Let us return to the framework of \secref{gen framework of functor}. 
Let $\CM_{w_0}$ be given by \eqref{M w0}. From
the construction we have:
\begin{lem} \label{coherent conv on affine}
For $\CT$, which is the direct image of a quasi-coherent sheaf on an
affine locally closed subscheme of $\Fl^\cG$, there exists a canonical
isomorphism
$$\sF(\CT)\simeq \CT\overset{R}\ast \CM_{w_0}.$$
\end{lem}

We now define the functor
$$\sF:D^+\left(\QCoh(\Fl^\cG)\right)\to D^+\left(\CC\right)$$
by the formula
$$\sF(\CT):=\CT\overset{R}\ast \CM_{w_0}.$$

\lemref{coherent conv on affine} implies that for $\CT$ being the
direct image of a quasi-coherent sheaf on an affine locally closed
subscheme of $\Fl^\cG$, this definition agrees with the one
of \secref{gen framework of functor}.

\medskip

Let $\CL^\cla_{\Fl^\cG_{w_0}}$ denote the line equal to the fiber of
the line bundle $\CL^\cla_{\Fl^\cG}$ at the point $w_0\in \Fl^\cG$. We
will regard it as equipped with an action of $\cB$ given by the
character $w_0(\cla)$.  

We have the following general assertion, valid in the set-up of
\secref{coherent convolutions}, whose proof is straightforward:
\begin{lem}   \label{coh conv equiv}
Let $\CM$ be a $\cG$-equivariant object of $\CC$. Then for $\CT\in
\QCoh(\Fl^\cG)$ and $\cmu\in \cLambda$ there exists a canonical
isomorphism
$$\CT\overset{R}\ast \Bigl(\CL^{\cmu}_{\Fl^\cG_{w_0}}\otimes \CM\Bigr)
\simeq R\Gamma(\Fl^\cG,\CT\otimes \CL^\cmu_{\Fl^\cG})\otimes \CM.$$
\end{lem}

By construction, for $\cla\in \cLambda^+$ we
have a map
\begin{equation} \label{vac to verm}
\CL^{-\cla}_{\Fl^\cG_{w_0}}\otimes \CM_\cla \to \CM_{w_0},
\end{equation}
of $\cB$-equivariant objects of
$\CC$. Hence, by functoriality, and using \lemref{coh conv equiv},
we obtain a map
\begin{equation} \label{gen A}
\CM_\cla\to \CL^\cla_{\Fl^\cG}\overset{R}\ast
\CM_{w_0}:=\sF(\CL^\cla_{\Fl^\cG}).
\end{equation}


We are going to show now that under some additional hypotheses the map
\eqref{gen A} is in fact an isomorphism for large enough $\cla \in
\cLambda^+$. We will make the following two additional assumptions:

\medskip

\noindent{(A)} The maps $\beta^{\cmu,\cla}:V^\cmu\otimes \CM_\cla\to
\CM_{\cla+\cmu}$ are surjective once $\cla$ is deep enough in the
dominant chamber. (In our main example, namely,
$\on{D}(\Gr_G)^{\Hecke}_\crit\mod$, the surjectivity will hold for all
$\cla$ that are regular.)

\medskip

\noindent{(B)} The objects $\CM_\cmu$ are Artinian, and for every
$\cB$-equivariant object $\CM$, which is Artinian as an object of
$\CC$, the higher cohomologies of $\CL^\cla_{\Fl^\cG}\overset{R}\ast
\CM$ vanish for all $\cla$ that are large enough (i.e., deep enough in
the dominant chamber).

We start with the following

\begin{lem}    \label{surj}
Assume that condition (A) holds. Then the map
\begin{equation} \label{dual to baby}
\CL^{-\cmu}_{\Fl^\cG_{w_0}}\otimes \CM_\cmu \to \CM_{w_0}
\end{equation}
is surjective for $\cmu$ deep inside the dominant chamber.
\end{lem}

\begin{proof}
By assumption (A), there exists $\cmu \in \cLambda^+$ such that the
maps $\beta^{\cnu,\cmu}$ are surjective for all $\cnu\in \cLambda^+$.
Let $\cmu'\in \cLambda^+$ be such that
$\cmu'-\cmu=:\cnu\in \cLambda^+$. Then we have a commutative diagram
$$
\CD
\CL^{-\cmu'}_{\Fl^\cG_{w_0}}\otimes \CM_{\cmu'} @>>> \CM_{w_0}  \\
@A{\beta^{\cnu,\cmu}}AA   @A{\text{\eqref{dual to baby}}}AA   \\
\CL^{-\cmu'}_{\Fl^\cG_{w_0}}\otimes \uV^\cnu \otimes \CM_{\cmu}
@>>> \CL^{-\cmu}_{\Fl^\cG_{w_0}}\otimes \CM_\cmu,
\endCD
$$
and hence the image of the upper horizontal arrow is contained
in that of \eqref{dual to baby}, since the left vertical arrow is
surjective.

For an arbitrary $\cmu'' \in \cLambda^+$ we can find $\cnu'\in
\cLambda^+$ such that $\cmu':=\cmu''+\cnu'$ satisfies $\cmu'-\cmu\in
\cLambda^+$. Then we have a commutative diagram
$$
\CD
\CL^{-\cmu'}_{\Fl^\cG_{w_0}}\otimes \uV^{\cnu'}\otimes \CM_{\cmu''} @>>> 
\CL^{-\cmu''}_{\Fl^\cG_{w_0}}\otimes \CM_{\cmu''} \\
@V{\beta^{\cnu',\cmu''}}VV   @VVV  \\
\CL^{-\cmu'}_{\Fl^\cG_{w_0}}\otimes \CM_{\cmu'} @>>> \CM_{w_0}.
\endCD
$$
Since the upper horizontal arrow in this diagram
(which comes from the map $\uV^{\cnu'}\otimes\CO_{\Fl^\cG}\to
\CL^{\cnu''}_{\Fl^\cG}$) is surjective, we obtain that the
image of $\CL^{-\cmu''}_{\Fl^\cG_{w_0}}\otimes \CM_{\cmu''}$
in $\CM_{w_0}$ is contained in that of 
$\CL^{-\cmu'}_{\Fl^\cG_{w_0}}\otimes \CM_{\cmu'}$, and the latter
is contained in the image of \eqref{dual to baby}, as we have seen
above.

Since $\CM_{w_0}$ is by definition the quotient of the direct sum
$\underset{\cla \in \cLambda^+}\oplus\, \CL^{-\cla}_{\Fl^\cG_{w_0}}\otimes
\CM_\cla$, whose components map to $\CM_{w_0}$ by means of
\eqref{dual to baby}, we obtain the statement of the lemma. 
\end{proof}

\begin{prop}  \label{verify A}
Assume that conditions (A) and (B) hold. Then the map \eqref{gen A}
is an isomorphism for all $\cla\in \cLambda^+$ that are large enough.
\end{prop}

\begin{proof}
By assumption (B) and \lemref{surj}, $\CM_{w_0}$ is Artinian. Applying
assumption (B) once more, we obtain that $\CL^\cla_{\Fl^\cG}
\overset{R}\ast\CM_{w_0}$ has no higher cohomologies for $\cla$ large
enough.

\medskip

Let us prove the surjectivity of the map 
\begin{equation} \label{la large}
\CM_\cla\to \CL^\cla_{\Fl^\cG}\ast\CM_{w_0},
\end{equation}
provided that $\cla$ is large.

Let $\CM'$ denote the kernel of this map. Let $\cla$ be large so that
$\CL^\cla_{\Fl^\cG}\overset{R}\ast \CM'$ has no higher
cohomologies. Assume also that $\cla-\cmu$ is dominant.

We obtain a short exact sequence
$$0\to \CL^\cla_{\Fl^\cG}\ast \CM'\to 
\CL^\cla_{\Fl^\cG}\ast (\CL^{-\cmu}_{\Fl^\cG_{w_0}}\otimes \CM_\cmu)
\twoheadrightarrow \CL^\cla_{\Fl^\cG}\ast\CM_{w_0}\to 0,$$
that fits into the commutative diagram
$$
\CD
\CL^\cla_{\Fl^\cG}\ast (\CL^{-\cmu}_{\Fl^\cG_{w_0}}\otimes \CM_\cmu)
@>>>  \CL^\cla_{\Fl^\cG}\ast\CM_{w_0} \\
@V{\text{\lemref{coh conv equiv}}}VV    @A{\text{\eqref{la large}}}AA
\\ \uV^{\cla-\cmu}\otimes \CM_\cmu  @>{\beta^{\cla-\cmu,\cmu}}>>
\CM_\cla,
\endCD
$$
implying surjectivity of the right vertical arrow for $\cla$ as above.

\medskip

Now we prove the injectivity of \eqref{la large}.  Let $\CM''$ denote
the kernel of the map \eqref{vac to verm}.  By (B), we can find a
finite collection of elements $(\cmu_i,\cnu_i)$, so that $\CM''$ is
contained in the image of $\underset{i}\oplus\, \uV^{\cnu_i}\otimes
\CL^{-\cmu_i-\cnu_i}_{\Fl^\cG_{w_0}}\otimes \CM_{\cmu_i}$ under the
difference of the two maps in \eqref{two arr gen}. Let $\CM'''$ denote
the pre-image of $\CM''$ in the above direct sum. Let us assume that
$\cla$ is such that $\CL^\cla_{\Fl^\cG}\overset{R}\ast \CM'''$ does
not have higher cohomologies. To prove the injectivity of \eqref{la
large}, it suffices to prove that the map
$$\CL^\cla_{\Fl^\cG}\ast \CM'''\to \CL^\cla_{\Fl^\cG}\ast
\CM''\hookrightarrow \CL^\cla_{\Fl^\cG}\ast
(\CL^{-\cla}_{\Fl^\cG_{w_0}}\otimes \CM_\cla) \simeq \CM_\cla$$ is
zero, but this follows from the commutativity of \eqref{la mu and mu}.

\end{proof}

\ssec{An explicit resolution}

Let us return to the case when $\CT\in \QCoh(\Fl^\cG)$ is the direct
image of a quasi-coherent sheaf on an affine subvariety of $\cG$. For
a non-negative integer $i$ consider the following object of $\CC$
$$\sF(\CT)^{+(i)}:=\underset{\cmu_i,...,\cmu_1,\cla\in
\cLambda^+}\oplus\, \Gamma(\Fl^\cG,\CT\otimes
\CL^{-\cla-\cmu_1-...-\cmu_i}_{\Fl^\cG})\otimes
\uV^{\cmu_1}\otimes...\otimes \uV^{\cmu_i}\otimes \CM_\cla.$$ We have
$i+1$ maps $\fd^j_i:\sF(\CT)^{+(i)}\to \sF(\CT)^{+(i-1)}$, where
$\fd^0_i$ comes from \eqref{beta basic}, $\fd^i_i$ comes from
$\beta^{\cmu_i,\cla}$, and $\fd^k_i$ with $0<k<i$ comes from the
natural map $\uV^{\cmu_{k}}\otimes \uV^{\cmu_{k+1}}\to
\uV^{\cmu_{k}+\cmu_{k+1}}$.

Let $\fd_i$ be the alternating sum of the $\fd^j_i$'s. The commutativity
of \eqref{la mu and mu} implies that $\fd_{i-1}\circ \fd_i=0$, i.e.,
$$...\to \sF(\CT)^{+(i+1)}\to \sF(\CT)^{+(i)}\to \sF(\CT)^{+(i-1)}\to...$$
is a complex. We will denote it by $\fC(\sF(\CT))^+$. 

\begin{prop}  \label{+ resolution}
Under the assumptions of \propref{verify A}, the natural map
$\fC(\sF(\CT))^+\to \sF(\CT)$ is a quasi-isomorphism.
\end{prop}

Before giving a proof let us consider the following construction. For
an element $\cla_0\in \cLambda^+$, let $\fC(\sF(\CT))^+_{\geq
\cla_0}$ be a sub-complex of $\fC(\sF(\CT))^+$ consisting of terms
$$\sF(\CT)^{+(i)}_{\geq
\cla_0}:=\underset{\cmu_1,...,\cmu_i,\cla;\cla-\cla_0\in
\cLambda^+}\oplus\, \Gamma(\Fl^\cG,\CT\otimes
\CL^{-\cla-\cmu_1-...-\cmu_i}_{\Fl^\cG})\otimes
\uV^{\cmu_1}\otimes...\otimes \uV^{\cmu_i}\otimes \CM_\cla.$$

The next assertion holds in the general
framework of \secref{gen framework of functor}:

\begin{lem} \label{restr lambda}
The embedding 
\begin{equation} \label{lambda emb} 
\fC(\sF(\CT))^+_{\geq \cla_0}\hookrightarrow 
\fC(\sF(\CT))^+
\end{equation}
is a quasi-isomorphism.
\end{lem}

\begin{proof}

As in \eqref{two arr gen}, we have two pairs of maps
$$\fC(\sF(\CT\otimes \CL^{-\cla_0}_{\Fl^\cG}))^+\otimes
\uV^{\cla_0}\rightrightarrows \fC(\sF(\CT))^+, \text{ and }
\fC(\sF(\CT\otimes \CL^{-\cla_0}_{\Fl^\cG}))^+_{\geq \cla_0}\otimes
\uV^{\cla_0}\rightrightarrows \fC(\sF(\CT))^+_{\geq \cla_0}$$ but it
is easy to see that they are tautologically pairwise homotopic.

Let $\CY\subset \Fl^\cG$ be an affine subvariety, which
$\CT$ was the direct image from. Let us choose a splitting,
denoted by $\fs$, of the surjection
$$\CO_{\CY}\otimes \uV^{\cla_0}\to \CL^{\cla_0}_{\Fl^\cG}|_{\CY}.$$
Then for each $\cmu$ we obtain a splitting
$$\Gamma(\Fl^\cG,\CT \otimes \CL^{-\cmu-\cla_0}_{\Fl^\cG})\otimes
\uV^{\cla_0} \overset{\fs}\leftrightarrows \Gamma(\Fl^\cG,\CT \otimes
\CL^{-\cmu}_{\Fl^\cG}).$$

We obtain a map
$$\fC(\sF(\CT))^+ \overset{\fs}\to \fC(\sF(\CT\otimes
\CL^{-\cla_0}_{\Fl^\cG}))^+\otimes \uV^{\cla_0}
\overset{\beta^{\cla_0,?}}\longrightarrow \fC(\sF(\CT))^+_{\geq
\cla_0},$$ which is easily seen to be the inverse to the embedding
\eqref{lambda emb}, up to homotopy.

\end{proof}

To prove \propref{+ resolution} we shall first consider the case when
$\CC=\QCoh(\Fl^\cG)$ and $\sF$ is the identity functor.

\begin{lem} \label{+ resolution coh}
\propref{+ resolution} holds for $\CC=\QCoh(\Fl^\cG)$ and
$\sF=\on{Id}$.
\end{lem}

\begin{proof}

In the case under consideration, we will use a shorthand notation
$\fC(\CT)^+$ (resp., $\fC(\CT)^+_{\geq \cla_0}$).

For a given $\CT$ it suffices to
show that for any $\cla_0\in \cLambda^+$ the complex of vector spaces
$\Gamma\Bigl(\Fl^\cG,\fC(\CT)^+\otimes
\CL^{\cla_0}_{\Fl^\cG}\Bigr)$ is quasi-isomorphic to
$\Gamma(\Fl^\cG,\CT\otimes \CL^{\cla_0}_{\Fl^\cG})$.

The above complex is tautologically isomorphic to the complex
$\Gamma\Bigl(\Fl^\cG,\fC(\CT\otimes
\CL^{\cla_0}_{\Fl^\cG})^+_{\geq \cla_0}\Bigr)$, and by 
\lemref{restr lambda}, the latter is quasi-isomorphic to
$\Gamma\Bigl(\Fl^\cG,\fC (\CT\otimes
\CL^{\cla_0}_{\Fl^\cG})^+\Bigr)$. Setting $\CT':= \CT\otimes
\CL^{\cla_0}_{\Fl^\cG}$, we need to verify the exactness of the
complex
\begin{multline*}
...\to\underset{\cmu_1,...,\cmu_i,\cmu_{i+1}\in \cLambda^+}\oplus\,
\Gamma(\Fl^\cG,\CT'\otimes
\CL^{-\cmu_1-...-\cmu_i-\cmu_{i+1}}_{\Fl^\cG})\otimes
\uV^{\cmu_1}\otimes...\otimes \uV^{\cmu_i}\otimes
\uV^{\cmu_{i+1}}\to... \\ ...\to \underset{\cmu_1,...,\cmu_i\in
\cLambda^+}\oplus\, \Gamma(\Fl^\cG,\CT'\otimes
\CL^{-\cmu_{1}-...-\cmu_i}_{\Fl^\cG})\otimes
\uV^{\cmu_{1}}\otimes...\otimes \uV^{\cmu_i}\to...\\
\to\underset{\cmu\in \cLambda^+}\oplus\, \Gamma(\Fl^\cG,\CT'\otimes
\CL^{-\cmu}_{\Fl^\cG})\otimes \uV^\cmu\to \Gamma(\Fl^\cG,\CT').
\end{multline*}
However, this complex admits an explicit homotopy 
\begin{multline*}
\sh^i: \underset{\cmu_1,...,,\cmu_i\in \cLambda^+}\oplus\,
\Gamma(\Fl^\cG,\CT'\otimes
\CL^{-\cmu_{1}-...-\cmu_i}_{\Fl^\cG})\otimes
\uV^{\cmu_{1}}\otimes...\otimes \uV^{\cmu_i}\to \\ \to
\underset{\cmu_1,...,\cmu_i,\cmu_{i+1}\in \cLambda^+}\oplus\,
\Gamma(\Fl^\cG,\CT'\otimes
\CL^{-\cla-\cmu_1-...-\cmu_i-\cmu_{i+1}}_{\Fl^\cG})\otimes
\uV^{\cmu_1}\otimes...\otimes \uV^{\cmu_i}\otimes \uV^{\cmu_{i+1}},
\end{multline*}
obtained by taking $\cmu_{i+1}=0$.

\end{proof}

Now we are ready to prove \propref{+ resolution}:

\begin{proof}

According to \lemref{restr lambda}, 
it suffices to show that $$\fC(\sF(\CT))^+_{\geq \cla_0}\to
\CT\overset{R}\ast \CW_{w_0}$$ is a quasi-isomorphism for any
$\cla_0\in \cLambda^+$. We take $\cla_0$ to be sufficiently large so
that the conclusion of \propref{verify A} holds.

Then we have:
$$\fC(\sF(\CT))^+_{\geq \cla_0}\simeq
\fC(\CT)^+_{\geq \cla_0}\overset{R}\ast \CW_{w_0},$$
and the assertion follows from the fact that
$\fC(\CT)^+_{\geq \cla_0}$ is quasi-isomorphic to
$\CT$, as we have seen above.

\end{proof}

\ssec{The functor}

We are now going to apply the above general discussion to the category
$\CC=\on{D}(\Gr_G)^{\Hecke}_\crit\mod$, and construct the sought-after
functor $$\sE:D^+\bigl(\QCoh\bigl(\Fl^{\cG}\bigr)\bigr)\to
D^+\bigl(\on{D}(\Gr_G)^{\Hecke}_\crit\mod\bigl).$$

\medskip

Let $R(\cG)$ denote the left-regular representation of $\cG$, and let
$\CF_{R(\cG)}$ be the corresponding object of
$\on{D}(\Gr_G)_\crit\mod$. However, it is easy to see that
$\CF_{R(\cG)}$ is naturally an object of
$\on{D}(\Gr_G)^{\Hecke}_\crit\mod$, and, moreover, as such it is
$\cG$-equivariant. Furthermore, the assignment
$$\CF\mapsto \CF\underset{G[[t]]}\star \CF_{R(\cG)}$$
defines an equivalence between $\on{D}(\Gr_G)_\crit\mod$ and the category 
of $\cG$-equivariant objects of $\on{D}(\Gr_G)^{\Hecke}_\crit\mod$. 

Following \cite{FG4}, we will denote the above functor also by
$\on{Ind}^\Hecke$.  When viewed as a functor
$\on{D}(\Gr_G)_\crit\mod\to \on{D}(\Gr_G)^\Hecke_\crit\mod$, it is the
left adjoint to the tautological forgetful functor.

\medskip

For an element $\cla\in \cLambda^+$ consider the $\cG$-equivariant
object of $\on{D}(\Gr_G)^{\Hecke}_\crit\mod$ equal to
$\on{Ind}^\Hecke(j_{\cla,\Gr_G,*})$.

\begin{lemconstr}   \label{la mu map}
For two elements $\cla,\cmu\in \cLambda^+$ there exists a canonical
map
\begin{equation} \label{j lambda and mu}
\uV^\cmu\otimes \on{Ind}^\Hecke(j_{\cla,\Gr_G,*})\to
\on{Ind}^\Hecke(j_{\cla+\cmu,\Gr_G,*}).
\end{equation}
\end{lemconstr}

\begin{proof}

We can rewrite
$$\uV^\cmu\otimes \on{Ind}^\Hecke(j_{\cla,\Gr_G,*})\simeq
j_{\cla,*}\star \CF_{V^\cmu}\underset{G[[t]]}\star \CF_{R(\cG)}$$ and
$$\on{Ind}^\Hecke(j_{\cla+\cmu,\Gr_G,*})\simeq
j_{\cla,*}\star j_{\cmu,\Gr_G,*}\underset{G[[t]]}\star \CF_{R(\cG)}.$$

The sought-after morphism comes from the map
$$\CF_{V^\cmu}\simeq \IC_{\Grb^\cmu}\to j_{\cmu,\Gr_G,*}.$$

\end{proof}

It is easy to see that the morphism of \lemconstrref{la mu map}
respects the natural $\cG$-equivariant structures on both sides,
and the resulting diagrams \eqref{la mu and mu} are
commutative.

Thus, we find ourselves in the framework of Sections \ref{gen
framework of functor} and \ref{gen framework, contd}. We set $\sE$ to
be the resulting functor $D^+\bigl(\QCoh\bigl(\Fl^{\cG}\bigr)\bigr)\to
D^+\bigl(\on{D}(\Gr_G)^{\Hecke}_\crit\mod\bigl)$. We let $\CW_{w_0}$
denote the object $\sE(\CW_0)\in \on{D}(\Gr_G)^{\Hecke}_\crit\mod$.

\ssec{}

By the results of the previous subsection, for every $\cla\in \cLambda^+$
we have a canonical map
\begin{equation} \label{ver of A}
\on{Ind}^\Hecke(j_{\cla,\Gr_G,*})\to \CL^\cla_{\Fl^\cG}\overset{R}\ast
\CW_{w_0}=:\sE(\CL^\cla_{\Fl^\cG}).
\end{equation}

We claim that assumptions (A) and (B) of \propref{verify A} hold
in our situation. Indeed,
assumption (A) follows from \cite{ABBGM}, Proposition 2.3.2.
The Artinian property of assumption (B) follows from 
\cite{ABBGM}, Corollary 1.3.10. To prove the cohomology vanishing
property, we can assume that the object in question is
irreducible. Hence, again by \cite{ABBGM}, Corollary 1.3.10,
it is of the form $\on{Ind}^\Hecke(\CF)$ for some irreducible 
$\CF\in  \on{D}(\Gr_G)_\crit\mod^I$, and the $\cB$-equivariant
structure is given by some character $\cla'\in \cLambda$.
The vanishing of higher cohomologies follows from the
Bott-Borel-Weil theorem for all $\cla$ with $\cla+\cla'\in \cLambda^+$.

\medskip

Thus, we obtain that the map \eqref{ver of A}
is an isomorphism for all $\cla$ that are deep enough in the
dominant chamber. We shall now establish a strengthening of
this result:

\begin{prop} \label{verify full A}
The map \eqref{ver of A} is an isomorphism for any
$\cla\in \cLambda^+$.
\end{prop}

This proposition implies that our functor $\sE$ coincides with the one
coming from the equivalence of \cite{ABG}. Before giving a proof we
need to make a couple of observations.

\medskip

Since the objects $j_{\cla,\Gr_G,*}$ are $I$-equivariant, by
construction, the functor $\sE$ factors as
$$D^+\bigl(\QCoh\bigl(\Fl^{\cG}\bigr)\bigr)\to
D^+\bigl(\on{D}(\Gr_G)^{\Hecke}_\crit\mod\bigl)^I\to
D^+\bigl(\on{D}(\Gr_G)^{\Hecke}_\crit\mod\bigl).$$ In particular, its
image belongs to $D^+(\on{D}(\Gr_G)^{\Hecke}_\crit\mod\bigl)^{I^0}$,
the latter being a full triangulated subcategory of
$D^+\bigl(\on{D}(\Gr_G)^{\Hecke}_\crit\mod\bigl)$.

\begin{prop}   \label{coh and I conv}
For $\cla\in \cLambda^+$ and $\CT\in
D^+\bigl(\QCoh\bigl(\Fl^{\cG}\bigr)\bigr)$ there exists a functorial
isomorphism
$$\sE(\CL^\cla_{\Fl^\cG}\otimes \CT)\simeq j_{\cla,*}\star \sE(\CT)\in 
D^+\bigl(\on{D}(\Gr_G)^{\Hecke}_\crit\mod\bigl)^I$$
\end{prop}

\begin{proof}
By the construction of the functor $\sE$, it suffices to show that
$j_{\cla,*}\star \CW_{w_0}$ is isomorphic to
$\CL^{\cla}_{\Fl^\cG_{w_0}}\otimes \CW_{w_0}$ as a $\cB$-equivariant
object of $\on{D}(\Gr_G)^{\Hecke}_\crit\mod^I$.

Since the functor $j_{\cla,*}\star ?$ is right exact, $j_{\cla,*}\star
\CW_{w_0}$ equals the co-equalizer of the maps
$$\underset{\cmu,\cla'}\oplus\, 
\uV^\cmu\otimes \CL^{-\cla'-\cmu}_{\Fl^\cG_{w_0}}\otimes 
j_{\cla,*}\star \on{Ind}^\Hecke(j_{\cla',\Gr_G,*}) \rightrightarrows 
\underset{\cla''}\oplus\, \CL^{-\cla''}_{\Fl^\cG_{w_0}}\otimes 
\on{Ind}^\Hecke(j_{\cla'',\Gr_G,*}).$$ 

Recall that $j_{\cla,*}\star \on{Ind}^\Hecke(j_{\cla',\Gr_G,*}) \simeq
\on{Ind}^\Hecke(j_{\cla+\cla',\Gr_G,*})$. Therefore, we obtain that 
the assertion of the proposition holds from \lemref{restr lambda}.

\end{proof}

We are now ready to prove \propref{verify full A}.

\begin{proof}

Let $\cla$ be as in the proposition, and let $\cla'$ be large
so that the map
$$\on{Ind}^\Hecke(j_{\cla+\cla',\Gr_G,*})\to
\sE(\CL^{\cla+\cla'}_{\Fl^\cG})$$ is an isomorphism. Let us apply the
functor $j_{\cla',*}\star ?$ to the two sides of \eqref{ver of A}. We
obtain a commutative diagram
$$ \CD j_{\cla',*}\star \on{Ind}^\Hecke(j_{\cla,\Gr_G,*})
@>{\text{\eqref{ver of A}}}>> j_{\cla',*}\star
\sE(\CL^{\cla}_{\Fl^\cG}) \\ @VVV @V{\text{\propref{coh and I
conv}}}VV \\ \on{Ind}^\Hecke(j_{\cla+\cla',\Gr_G,*})
@>{\text{\eqref{ver of A}}}>> \sE(\CL^{\cla+\cla'}_{\Fl^\cG}), \endCD
$$ and our assertion follows from the fact that the functor
$j_{\cla',*}\star ?$ is a self-equivalence of
$D^+(\on{D}(\Gr_G)^{\Hecke}_\crit\mod)^I$

\end{proof}

\section{Identification of Wakimoto modules}   \label{statement full ident}

The goal of this section is to formulate a result, \thmref{full
ident}, that will describe Wakimoto modules in terms of the
equivalence \eqref{Gamma Hecke} and the functor $\sE$.

\ssec{The category}  
We begin by spelling 
out the constructions of the previous section in a relative
situation, namely, over the scheme $\Spec(\fZ_\fg^\reg) =
\Op_{\cg}^\reg$.
Recall that over $\Spec(\fZ^\reg_\fg)$ there exists a canonical
$\cG$-torsor, denoted $\CP_{\cG,\fZ}$, that corresponds to the
tautological $\cG$-torsor $\CP_{\cG,\Op^\reg_\fg}$ over $\Op^\reg_\cg$
under the isomorphism $\Spec(\fZ^\reg_\fg)\to \Op^\reg_\cg$ induced by
$\on{map}_{\alg} = \on{map}_{\geom}$. We will denote by
$\CP_{\cB,\fZ}$ its reduction to $\cB$, corresponding to the oper
structure on $\CP_{\cG,\Op^\reg_\fg}$. For $V\in \Rep(\cG)$, recall
that $\CV_{\fZ}$ denotes the associated vector bundle over
$\Spec(\fZ^\reg_\fg)$.

\medskip

Following \cite{FG4}, we introduce the category
$\on{D}(\Gr_G)^{\Hecke_\fZ}_\crit\mod$ as follows. Its objects are the
data of $(\CF,\alpha_V,\,\, \forall\,\, V\in \Rep(\cG))$, where $\CF$
is an object of $\on{D}(\Gr_G)_\crit\mod$, endowed with an action of
the algebra $\fZ^\reg_\fg$ by endomorphisms, and $\alpha_V, V\in
\Rep(\cG))$, are isomorphisms of D-modules
$$\al_V: \CF\star \CF_V \overset{\sim}\longrightarrow
\CV_\fZ\underset{\fZ^\reg_\fg}\otimes \CF,$$ compatible with the
action of $\fZ^\reg_\fg$ on both sides, and such that the following
two conditions are satisfied:

\begin{itemize}

\item
If $V$ is the trivial representations $\BC$, then the morphism
$\alpha_V$ is the identity map.

\item
For $V,W\in \Rep(\cG)$ and $U:=V\otimes W$, the diagram
$$
\CD
(\CF\star \CF_V)\star \CF_W  @>{\sim}>>  \CF\star \CF_U \\
@V{\alpha_V\star \on{id}_{\CF_W}}VV  @V{\alpha_U}VV \\
(\CV_\fZ\underset{\fZ^\reg_\fg}\otimes \CF)\star \CF_W  & &
\CU_\fZ\underset{\fZ^\reg_\fg}\otimes \CF \\
@V{\sim}VV   @V{\sim}VV  \\
\CV_\fZ\underset{\fZ^\reg_\fg}\otimes (\CF\star \CF_W) 
@>{\on{id}_{\CV_\fZ}\otimes \alpha_W}>>  
\CV_\fZ\underset{\fZ^\reg_\fg} \otimes
\CW_\fZ\underset{\fZ^\reg_\fg} \otimes\CF
\endCD
$$ is commutative.

\end{itemize}

Morphisms in this category between $(\CF,\alpha_V)$ and
$(\CF',\alpha'_V)$ are maps of D-modules $\phi:\CF\to \CF'$ that are
compatible with the actions of $\fZ^\reg_\fg$ on both sides, and such
that
$$(\on{id}_{\CV_\fZ}\otimes \phi)\circ \alpha_V=
\alpha'_V\circ (\phi\star \on{id}_{\CF_V}).$$

\medskip

Note that the $\cG$-torsor $\CP_{\cG,\fZ}$ can be (non-canonically)
trivialized.  A choice of such a trivialization defines an
equivalence between $\on{D}(\Gr_G)^{\Hecke_\fZ}_\crit\mod$ and the
category consisting of objects of $\on{D}(\Gr_G)^{\Hecke}_\crit\mod$,
endowed with an action of the algebra $\fZ^\reg_\fg$.

\medskip

Let $R(\cG)^\fZ$ denote the direct image of $\fZ^\reg_\fg$ onto
$\on{pt}/\cG$, regarded as an object of $\Rep(\cG)$, endowed with a
commuting action of $\fZ^\reg_\fg$. Let $\CF_{R(\cG)^\fZ}$ be the
corresponding object of $\on{D}(\Gr_G)^{\Hecke_\fZ}_\crit\mod$.

As in the case of $\CF_{R(\cG)}$, one shows that $\CF_{R(\cG)^\fZ}$ is
naturally an object of
$\on{D}(\Gr_G)^{\Hecke_\fZ}_\crit\mod$. Moreover, the assignment
$$\CF\mapsto \CF\underset{G[[t]]}\star \CF_{R(\cG)^\fZ}$$ defines a
functor $\on{D}(\Gr_G)_\crit\mod\to
\on{D}(\Gr_G)^{\Hecke_\fZ}_\crit\mod$, denoted,
$\on{Ind}^{\Hecke_\fZ}$, which is the left adjoint to the tautological
forgetful functor.

\ssec{The functor to $\hg_\crit\mod_\reg$}

Let us now recall the definition of the functor
$\Gamma^{\Hecke_\fZ}:\on{D}(\Gr_G)^{\Hecke_\fZ}_\crit\mod\to
\hg_\crit\mod_\reg$, following \cite{FG2} and \cite{FG4}.

\medskip

Consider the groupoid
\begin{equation} \label{consider groupoid}
\Isom_{\cG,\fZ_\fg^\reg}:\Spec(\fZ^\reg_\fg)\underset{\on{pt}/\cG}\times
\Spec(\fZ^\reg_\fg),
\end{equation}
where the morphism $\Spec(\fZ^\reg_\fg)\to \on{pt}/\cG$ corresponds to
the $\cG$-bundle $\CP_{\cG,\fZ}$ on $\Spec(\fZ^\reg_\fg)$. Let ${{\bf
1}_{\Isom_{\cG,\fZ_\fg^\reg}}}$ denote the unit section
$\Spec(\fZ^\reg_\fg)\to \Isom_{\cG,\fZ_\fg^\reg}$. As shown in {\em
loc. cit.}, for every object $\CF^H\in
\on{D}(\Gr_G)^{\Hecke_\fZ}_\crit\mod$, the $\hg_\crit$-module
$\Gamma(\Gr_G,\CF^H)$ carries a natural action of the algebra
$\Fun(\Isom_{\cG,\fZ_\fg^\reg})$ by endomorphisms. We then define the
functor $\Gamma^{\Hecke_\fZ}$ by $$\CF^H\mapsto \Gamma(\Gr_G,\CF^H)
\underset{\Fun(\Isom_{\cG,\fZ_\fg^\reg}),{\bf
1}^*_{\Isom_{\cG,\fZ_\fg^\reg}}}\otimes \fZ^\reg_\fg.$$

For $\CF\in \on{D}(\Gr_G)_\crit\mod$ we have a natural isomorphism:
$$\Gamma^{\Hecke_\fZ}(\Gr_G,\on{Ind}^{\Hecke_\fZ}(\CF))\simeq
\Gamma(\Gr_G,\CF).$$

\medskip

In \cite{FG2} we conjectured that the functor $\Gamma^{\Hecke_\fZ}$ is
exact and defines an equivalence of categories
$\on{D}(\Gr_G)^{\Hecke_\fZ}_\crit\mod$ and $\hg_\crit\mod_\reg$.
In \cite{FG4}, Theorem 1.7, we proved that the corresponding functor
\begin{equation}    \label{equivalence}
\on{D}(\Gr_G)^{\Hecke_\fZ}_\crit\mod^{I^0}\to
\hg_\crit\mod_\reg^{I^0},
\end{equation}
between the $I^0$-equivariant categories, is indeed an equivalence of
categories.

\ssec{The functor $\sE^\fZ$}  \label{sE Z}

Let $\cG_{\fZ}$ be the group-scheme over $\Spec(\fZ^\reg_\fg)$ equal
to the adjoint twist of $\cG$ by means of $\CP_{\cG,\fZ}$. By
construction, $\on{D}(\Gr_G)^{\Hecke_\fZ}_\crit\mod$, viewed as a
category over $\Spec(\fZ^\reg_\fg)$ carries an action of $\cG_{\fZ}$
in the same way as $\on{D}(\Gr_G)^{\Hecke}_\crit\mod$ carried an
action of $\cG$.

\medskip 

Let $\cB_\fZ$ denote the group-subscheme of $\cG_\fZ$, corresponding
to the reduction $\CP_{\cB,\fZ}$ of $\CP_{\cG,\fZ}$ to $\cB$.
Let $\Fl^\cG_\fZ$ be the scheme over $\Spec(\fZ^\reg_\fg)$, classifying
reductions of $\CP_{\cG,\fZ}$ to $\cB^-\subset \cG$, i.e.,
$$\Fl^{\cG}_\fZ =  \CP_{\cG,\fZ}\overset{\cG}\times \Fl^{\cG}.$$
The group-scheme $\cG_\fZ$ acts naturally on $\Fl^\cG_\fZ$.

The $\cB_\fZ$-orbits on $\Fl^\cG_\fZ$ are in a natural bijection with
the elements of the Weyl group; for $w\in W$ we will denote by
$\Fl^\cG_{w,\fZ}\subset \Fl^\cG_\fZ$ the Schubert stratum,
corresponding to the coset $\cB\cdot w^{-1}\cdot \cB^-$.

For $w=1$, this is an open subscheme of $\Fl^\cG_\fZ$. For $w=w_0$,
the subscheme $\Fl^\cG_{w_0,\fZ}\subset \Fl^\cG_\fZ$ is the section of
the natural projection $$p:\Fl^\cG_\fZ\to \Spec(\fZ^\reg_\fg),$$
corresponding to the reduction $\CP_{\cB,\fZ}$ of $\CP_{\cG,\fZ}$.
The stabilizer of $\Fl^\cG_{w_0,\fZ}$ in $\cG_\fZ$ is, by definition,
$\cB_\fZ$.

For $\cla\in \cLambda^+$, let $\CL^\cla_{\Fl^\cG_\fZ}$ denote the
corresponding $\cG_\fZ$-equivariant line bundle on $\Fl^\cG_\fZ$,
normalized so that
$$p_*(\CL^\cla_{\Fl^\cG_\fZ})\simeq \CV^\cla_\fZ.$$

\medskip

As in the case of $\Fl^\cG$ and $\on{D}(\Gr_G)^{\Hecke}_\crit\mod$, we
construct a functor
\begin{equation}    \label{sEfZ}
\sE^\fZ:D^+\bigl(\QCoh(\Fl^\cG_\fZ)\bigr)\to
D^+\bigl(\on{D}(\Gr_G)^{\Hecke_\fZ}_\crit\mod\bigr).
\end{equation}
It is obtained as the convolution $\sE^\fZ(\CT) = \CT\overset{R}\ast
\CW^\fZ_{w_0},$ where $\CW^\fZ_{w_0}$ is a $\cB_\fZ$-equivariant
object of $\on{D}(\Gr_G)^{\Hecke_\fZ}_\crit\mod$ that can be recovered
as $\sE^\fZ(\CO_{\Fl^\cG_{w_0,\fZ}})$.

\medskip

Explicitly, for $\CT$ being the direct image of a quasi-coherent sheaf
on an affine locally closed subscheme of $\Fl^\cG_\fZ$ (in particular,
for $\CT=\CO_{\Fl^\cG_{w_0,\fZ}}$), the object $\sE^\fZ(\CT)$ is
defined as follows:

Set
$$\sE^\fZ(\CT)^+:=\underset{\cla\in
\cLambda^+}\oplus \, p_*(\CT\otimes
\CL^{-\cla-\cmu}_{\Fl^\cG_\fZ})\underset{\fZ^\reg_\fg}\otimes
\on{Ind}^{\Hecke_\fZ}(j_{\cla,\Gr_G,*})$$ 
and
$$\sE^\fZ(\CT)^{++}:=\underset{\cmu,\cla\in\cLambda^+}\oplus \, 
\CV^\cmu_\fZ\underset{\fZ^\reg_\fg}\otimes 
p_*(\CT\otimes \CL^{-\cla-\cmu}_{\Fl^\cG_\fZ})
\underset{\fZ^\reg_\fg}\otimes
\on{Ind}^{\Hecke_\fZ}(j_{\cla,\Gr_G,*})$$

Then $\sE^\fZ(\CT)$ is, by definition, the co-equalizer of the map
$$\sE^\fZ(\CT)^{++}\rightrightarrows
\sE^\fZ(\CT)^+.$$ 

We now set $\CW^\fZ_{w_0} = \sE^\fZ(\CO_{\Fl^\cG_{w_0,\fZ}})$. It
follows from the construction that $\CW^\fZ_{w_0}$ is a
$\cB_\fZ$-equivariant object of
$\on{D}(\Gr_G)^{\Hecke_\fZ}_\crit\mod$. Having defined
$\CW^\fZ_{w_0}$, we define the functor \eqref{sEfZ} as follows:

We have a functor
$$
\on{act}_\fZ^*: \on{D}(\Gr_G)^{\Hecke_\fZ}_\crit\mod \to
\QCoh(\cG_\fZ)
\underset{\fZ^\reg_\fg}\otimes \on{D}(\Gr_G)^{\Hecke_\fZ}_\crit\mod
$$
corresponding to the action of $\cG_\fZ$ on
$\on{D}(\Gr_G)^{\Hecke_\fZ}_\crit\mod$. Let $\CT$ be a quasi-coherent
sheaf on $\Fl^\cG_\fZ$ and $\wt{\CT}$ the corresponding
$\cB_{\fZ}$-equivariant quasi-coherent sheaf on $\cG_\fZ$. We
set
\begin{equation} \label{defn sEZ}
\sE^\fZ(\CT):=\CT\overset{R}\ast \CW^\fZ_{w_0}:= 
R\on{Inv}\left(\cB_\fZ,\wt\CT \underset{\Fun(\cG_\fZ)}\otimes
\on{act}_\fZ^*(\CW^\fZ_{w_0})\right),
\end{equation}
where $R\on{Inv}\left(\cB_\fZ,?\right)$ denotes the derived functor of
$\cB_\fZ$-invariants.

The object of $D^+\left(\on{D}(\Gr_G)_\crit\mod\right)$, underlying 
$\sE^\fZ(\CT)$, can be explicitly written as follows:
$$\sE^\fZ(\CT)=
R\on{Inv}\left(\cB_\fZ,\wt\CT\underset{\fZ^\reg_\fg}\otimes
\CW^\fZ_{w_0}\right).$$

\medskip

The next lemma follows from \propref{verify full A}:

\begin{lem} \label{cond A & B}
For $\cla\in \cLambda^+$ we have:
$$\sE^\fZ(\CL^{\cla}_{\Fl^\cG_\fZ})\simeq
\on{Ind}^{\Hecke_\fZ}(j_{\cla,\Gr_G,*}).$$
\end{lem}

Let now $\cla$ and $\cmu$ be two elements of $\cLambda^+$. As in
\lemconstrref{la mu map} we have a canonical map
\begin{equation} \label{la mu map Z}
\CV^\cmu_\fZ\underset{\fZ^\reg_\fg}\otimes
\on{Ind}^{\Hecke_\fZ}(j_{\cla,\Gr_G,*})\to
\on{Ind}^{\Hecke_\fZ}(j_{\cla+\cmu,\Gr_G,*}),
\end{equation}
and a commutative diagram:
$$ \CD \sE^\fZ(p^*(\CV^\cmu_\fZ)\otimes \CL^{\cla}_{\Fl^\cG_\fZ}) @>>>
\CV^\cmu_\fZ\underset{\fZ^\reg_\fg}\otimes
\on{Ind}^{\Hecke_\fZ}(j_{\cla,\Gr_G,*}) \\ @VVV @VVV \\
\sE^\fZ(\CL^{\cla+\cmu}_{\Fl^\cG_\fZ}) @>>>
\on{Ind}^{\Hecke_\fZ}(j_{\cla+\cmu,\Gr_G,*}).  \endCD
$$

Finally, as in the case of $\sE$, the functor $\sE^\fZ$ factors as
$$D^+\bigl(\QCoh\bigl(\Fl^{\cG}_\fZ)\bigr)\to
D^+(\on{D}(\Gr_G)^{\Hecke_\fZ}_\crit\mod)^I\to
D^+(\on{D}(\Gr_G)^{\Hecke_\fZ}_\crit\mod),$$
and as in \propref{coh and I conv}, we have:
\begin{equation} \label{lambda inv T}
j_{\cla,*}\star \sE^\fZ(\CT)\simeq
\sE^\fZ(\CL^\cla_{\Fl^\cG_\fZ}\otimes \CT).
\end{equation}

\ssec{The functor $\sG$}

Let us denote by $\sG$ the functor
$$\Gamma^{\Hecke_\fZ}\circ
\sE^\fZ:D^+\bigl(\QCoh\bigl(\Fl^{\cG}_\fZ)\bigr)\to
D^+(\hg_\crit\mod_\reg)^{I^0}.$$ The goal of the rest of this paper is
to describe the relationship between this functor and Wakimoto
modules.

\ssec{$\Fl^\cG_\fZ$ and Miura opers}

Recall the D-scheme $\MOp_{\cg,X}$ of Miura opers, and let
$\MOp^\reg_{\cg}$ be the corresponding scheme, attached to
the formal disc. By the construction of $\on{map}_{\geom}$, 
it identifies the $\cG$-torsors $\CP_{\cG,\fZ}$ over  
$\Spec(\fZ^\reg_\fg)$ and $\CP_{\cG,\Op^\reg_\cg}$
over $\Op^\reg_\cg$. Hence, we have an isomorphism
\begin{equation} \label{flags and Miura}
\CD
\Fl^\cG_\fZ   @>{\sim}>>  \MOp^\reg_\cg \\
@V{p}VV    @VVV  \\
\Spec(\fZ^\reg_\fg) @>{\on{map}_{\geom}}>>  \Op^\reg_\cg.
\endCD
\end{equation}

For every $w\in W$, let $\MOp^{w,\reg}_{\cg}\subset \MOp^\reg_{\cg}$
be the corresponding Schubert cell; under the isomorphism
\eqref{flags and Miura}, $\MOp^{w,\reg}_{\cg}$ goes over
to $\Fl^\cG_{w,\fZ}$. Let $\MOp^{w,\on{th},\reg}_{\cg}$ 
(resp., $\Fl^\cG_{w,\on{th},\fZ}$) denote the formal neighborhood
of $\MOp^{w,\reg}_{\cg}$ in $\MOp^\reg_{\cg}$
(resp., of $\Fl^\cG_{w,\fZ}$ in $\Fl^\cG_{\fZ}$).

\medskip

Recall the D-scheme $\ConnX$, and let $\ConnD$ (resp., $\ConnDt$) be
the scheme (resp., ind-scheme) of its sections over the formal (resp.,
punctured) disc around $x$. Recall the isomorphism \eqref{Miura as
Con} $\ConnX\simeq \MOp_{\cg,\gen,X}$. Taking the fibers at $x\in X$
we obtain an isomorphism
$$\ConnD\simeq \MOp^{1,\reg}_{\cg}.$$

{}From \cite{FG2}, Sect. 3.6, we obtain that the above isomorphism
generalizes to the following canonical isomorphism
\begin{equation} \label{Miura as Con, general}
\ConnDt\underset{\Op_{\cg}(\D^\times)}\times \Op^\reg_\cg \simeq
\underset{w\in W}\sqcup\, \MOp^{w,\on{th},\reg}_{\cg}.
\end{equation}

{}From \eqref{geom diagram}, we obtain a commutative diagram
of ind-schemes
$$
\CD
\Spec(\wh\fH_\crit)  @>{\on{map}^M_{\geom}}>> \ConnDt  \\
@V{\varphi}VV    @V{\on{MT}}VV  \\
\Spec(\fZ_\fg)  @>{\on{map}_{\geom}}>> \Op_\cg(\D^\times).
\endCD
$$

Combing this with \eqref{Miura as Con, general}, we obtain an 
isomorphism
\begin{equation}  \label{induction parameters}
\Spec(\wh\fH_{\crit})\underset{\Spec(\fZ_{\fg})}\times
\Spec(\fZ^\reg_\fg) \simeq \underset{w\in W}\sqcup\,
\Fl^\cG_{w,\on{th},\fZ},
\end{equation}
which will play a crucial role for the rest of this paper.

\ssec{Relation to Wakimoto modules}

Recall from \secref{wak functor} the functor $\BW:\wh\fH_\crit\mod\to
\hg_\crit\mod$.  In particular, we obtain a functor
$$\QCoh\Bigl(\Spec(\wh\fH_{\crit})\underset{\Spec(\fZ_{\fg})}\times
\Spec(\fZ^\reg_\fg)\Bigr)\to \hg_\crit\mod_\reg.$$

Using \eqref{induction parameters}, for each element $w$ of the Weyl
group, we obtain a functor
\begin{equation} \label{Wak w functor}
_w\BW:\QCoh(\Fl^\cG_{w,\on{th},\fZ})\to \hg_\crit\mod_\reg.
\end{equation}

The direct image functor identifies $\QCoh(\Fl^\cG_{w,\on{th},\fZ})$
with a full subcategory of $\QCoh(\Fl^\cG_{\fZ})$. Our main result
compares the functors $_w\BW$ and
$\sG|_{\QCoh(\Fl^\cG_{w,\on{th},\fZ})}$:

\begin{thm} \label{full ident}
For $\CT\in \QCoh(\Fl^\cG_{w,\on{th},\fZ})$ there exists a canonical
isomorphism:
$$_w\BW(\CT)\simeq \sG\Bigl(\CT
\underset{\CO_{\Fl^\cG_{w,\on{th},\fZ}}}\otimes
\CL_{\Fl^\cG_\fZ}^{\crho-w(\crho)}|_{\Fl^\cG_{w,\on{th},\fZ}}\Bigr).$$
\end{thm}

\ssec{Some particular cases}   \label{part cases}

For a weight $\mu\in \fh^*$, let $\BW_{\crit,\mu}$ be the corresponding
Wakimoto module (see \secref{wak functor}). By definition, it is
induced from the $\wh\fH_\crit$-module $\Fun(\fH_\crit^{\RS,w_0(\mu)})$.

Its support over $\Spec(\fZ_\fg)$ is 
contained in $\Spec(\fZ^{\RS,\varpi(-w_0(\mu)+\rho)}_\fg)$.  Let us take
$\mu=w(\rho)-\rho$; then $\Spec(\fZ^{\RS,\varpi(-w_0(\mu)+\rho)}_\fg)=
\Spec(\fZ^\nilp_\fg)$.

We introduce the module
$$\BW_{\crit,w(\rho)-\rho,\reg} :=
\BW_{\crit,w(\rho)-\rho}\underset{\fZ_\fg^\nilp}\otimes
\fZ^\reg_\fg.$$ This is the maximal quotient of
$\BW_{\crit,w(\rho)-\rho}$ that belongs to the category
$\hg_\crit\mod_\reg$.  In the particular case when $w=1$, the module
$\BW_{\crit,0}$ is itself supported over $\Spec(\fZ^\reg_\fg)$, and
so $\BW_{\crit,0}\to \BW_{\crit,1,\reg}$ is an isomorphism.

\medskip

We would like now to apply \thmref{full ident} and describe the above
modules in terms of objects of $\on{D}(\Gr_G)^{\Hecke_\fZ}_\crit\mod$.

\medskip 

For every $\lambda\in \fh^*$ we can consider the D-scheme
$\ConnX^{\RS,\lambda}$, whose restriction to $X-x$ is isomorphic to
$\ConnX$, and whose fiber at $x$ identifies with
$\ConnD^{\RS,\lambda}$, where the latter is as in \cite{FG2},
Sect. 3.5.  By definition, $\ConnX^{\RS,\lambda}$ is the scheme of
connections on the $\cH$-bundle $\omega_X^\rho$ that have a pole of
order $1$ at $x$ with residue $\lambda$.

\medskip

Let us take $\lambda=\rho-w(\rho)$. By \cite{FG2}, Sect. 3.6, the
morphism
$$\ConnX\to \MOp_{\cg,X},$$ given by \eqref{Miura as Con}, composed
with the tautological embedding $\MOp_{\cg,\gen,X}\hookrightarrow
\MOp_{\cg,X}$ extends to a map of D-schemes
\begin{equation} \label{RS extension}
\ConnX^{\RS,\rho-w(\rho)}\underset{\Op^{\nilp}_{\cg,X}}\times
\Op^\reg_{\cg,X} \to \MOp_{\cg,X},
\end{equation}
where $\Op^{\nilp}_{\cg,X}$ is the D-scheme of opers with a nilpotent
singularity at $x$.

The resulting map on the level of fibers fits into a commutative
diagram
\begin{equation} \label{Miura RS}
\CD \ConnD^{\RS,\rho-w(\rho)}\underset{\Op^\nilp_\cg}\times
\Op^\reg_\cg @>{\sim}>> \MOp_\cg^{w,\reg} @>>> \MOp^\reg_\cg \\ @VVV
@VVV \\ \ConnDt\underset{\Op_{\cg}(\D^\times)}\times \Op^\reg_\cg
@>{\text{\eqref{Miura as Con, general}}}>> \underset{w\in W}\sqcup\,
\MOp^{w,\on{th},\reg}_{\cg}.  \endCD
\end{equation}

The left portion of the upper horizontal arrow is an isomorphism, and
hence it identifies
$\ConnD^{\RS,\rho-w(\rho)}\underset{\Op^\nilp_\cg}\times \Op^\reg_\cg$
with the reduced scheme of a connected component of
$\ConnDt\underset{\Op_{\cg}(\D^\times)}\times \Op^\reg_\cg$.

\medskip

Combining this with the isomorphism
$\fH_\crit^{\RS,w_0(w(\rho)-\rho)}\simeq \ConnD^{\RS,\rho-w(\rho)}$,
induced by $\on{map}^M_{\geom}$, we obtain an isomorphism:
\begin{equation} \label{H RS}
\Spec\Bigl(\fH_\crit^{\RS,w_0(w(\rho)-\rho)}\underset{\fZ^\nilp_\fg}
\otimes \fZ^\reg_\fg\Bigr)\simeq \Fl^{\cG}_{w,\fZ},
\end{equation}
where regular functions on the (affine) scheme on the LHS is, by
definition, the $\wh\fH_\crit$-module that
$\BW_{\crit,w(\rho)-\rho,\reg}$ is induced from.  In other words,
$$\BW_{\crit,w(\rho)-\rho,\reg}\simeq{}_w\BW(\CO_{\Fl^{\cG}_{w,\fZ}}).$$

\medskip

Let us introduce a short-hand notation
$\CW^\fZ_w:=\sE^\fZ(\CO_{\Fl^{\cG}_{w,\fZ}})$. As we will see below,
for every $\cla\in \cLambda$ and $w\in W$, the restriction of the line
bundle $\CL^\cla_{\Fl^\cG_\fZ}$ to $\Fl^{\cG}_{w,\fZ}$ is constant and
isomorphic to $\fl^\cla_{w(\rho)}\otimes \CO_{\Fl^{\cG}_{w,\fZ}}$,
where $\fl^\cla_{w(\rho)}$ is as in \propref{twist of Wakimoto,
local}.  Hence, we obtain:

\begin{thm} \label{ident w Wak}
There exists a canonical isomorphism
$$\BW_{\crit,w(\rho)-\rho,\reg}\otimes \fl^{w(\crho)-\crho}_{w(\rho)}
\simeq \Gamma^{\Hecke_\fZ} (\Gr_G,\CW_w^\fZ).$$
\end{thm}

We should note that a particular case of \thmref{ident w Wak}, namely,
for $w=w_0$ has been established in \cite{FG4}, and it was a key
calculation on which the proof of the main result of {\it loc. cit.}
was based. The method of proof of \thmref{ident w Wak} for a general
$w$ presented below will be quite different.

\medskip

Taking $w=1$, we obtain the following description of the "main"
Wakimoto module:

\begin{thm} \label{ident prince Wak}
There exists a canonical isomorphism
$$\BW_{\crit,0} \simeq \Gamma^{\Hecke_\fZ} (\Gr_G,\CW_1^\fZ).$$
\end{thm}

Our strategy for the proof of \thmref{full ident} will be as follows.
First, we will prove  \thmref{ident prince Wak}; this will be done
in \secref{principal Wakimoto} by a rather explicit argument.
Then, in \secref{other Wakimoto}, we will prove \thmref{ident w Wak}.
In \secref{more general} we will prove a weakened version
of \thmref{full ident}: namely, we will show that
$$_w\BW(\CT)\simeq \sG\Bigl(\CL_w^{\on{twist}}
\underset{\Fun(\Fl^\cG_{w,\on{th},\fZ})}\otimes \CT\Bigr),$$ where
$\CL_w^{\on{twist}}$ is a certain {\it bi-module} over the topological
algebra $\Fun(\Fl^\cG_{w,\on{th},\fZ})$.

Finally, in \secref{renorm Wak}, we will show that the left and right
actions of $\Fun(\Fl^\cG_{w,\on{th},\fZ})$ on $\CL_w^{\on{twist}}$
coincide, and that it is an invertible sheaf canonically isomorphic
to $\CL_{\Fl^\cG_\fZ}^{\crho-w(\crho)}|_{\Fl^\cG_{w,\on{th},\fZ}}$.

\ssec{A BGG type resolution}

As an application of \thmref{full ident}, we construct a BGG type
resolution of the vacuum module $\BV_{\crit}$.  For an element $w\in
W$ let $\on{Dist}_w$ denote the $\CO$-module on $\Fl^{\cG}_\fZ$
underlying the left D-module of distributions on the Schubert cell
$\Fl^{\cG}_{w,\fZ}$.  This $\CO$-module can be naturally thought of as
an object of $\QCoh(\Fl^\cG_{w,\on{th},\fZ})$.

The Cousin-Grothendieck resolution of the structure sheaf of
$\Fl^{\cG}_\fZ$ by means of $\on{Dist}_w$, combined with 
\lemref{cond A & B}, yields the following:

\begin{cor}   \label{screening BGG}
There exists a right resolution $C^\bullet$ of $\BV_{\crit}$, whose
$k$th term $C^k$ is isomorphic to $$\bigoplus_{w\in W,\ell(w)=k} \,
{}_w\BW\bigl(\CL^{w(\crho)-\crho}_{\Fl^{\cG}_{\fZ}}
\underset{\CO_{\Fl^{\cG}_\fZ}}\otimes\on{Dist}_w\bigr).$$
\end{cor}

In \secref{resolution} a more explicit realization of the modules
$_w\BW\bigl(\on{Dist}_w\underset{\CO_{\Fl^{\cG}_\fZ}}\otimes
\CL^{w(\crho)-\crho}_{\Fl^{\cG}_{\fZ}}\bigr)$, involved in this
resolution, will be obtained. In addition, we will make 
contact with a conjecture of \cite{FF}.

\section{Geometric realization of the Wakimoto module $\BW_{\crit,0}$}
\label{principal Wakimoto}

As was mentioned above, the goal of this section is to prove
\thmref{ident prince Wak}, which is a particular case of
\thmref{full ident}. First, we will give a more explicit description
of the objects $\CW^\fZ_w$, which is valid for any $w\in W$.

\ssec{Explicit construction of $\CW^\fZ_w$}   \label{explicit prince}

By definition, the object $\CW^\fZ_w\in
\on{D}(\Gr_G)^{\Hecke_\fZ}_\crit$ is obtained by applying the functor
$\sE^\fZ$ to $\CO_{\Fl^\cG_{w,\fZ}}$.  Since $\Fl^\cG_{w,\fZ}$ is
affine, $\sE^\fZ(\CO_{\Fl^\cG_{w,\fZ}})$ can be explicitly described
as a quotient as in \secref{sE Z}.

Let us denote the restriction
$\CL^{-\cla-\cmu}_{\Fl^\cG_\fZ}|_{\Fl^\cG_{w,\fZ}}$ by
$\CL^{-\cla-\cmu}_{\Fl^\cG_{w,\fZ}}$. In what follows we will not
distinguish in the notation between quasi-coherent sheaves on this
scheme and their global sections, thought of as
$\Fun(\Fl^\cG_{w,\fZ})$-modules.

Set 
$$\CW^{\fZ+}_w:=\bigoplus_{\cla\in \cLambda^+} \,
j_{\cla,*}\star
\on{Ind}^{\Hecke_\fZ}(\delta_{1,\Gr_G})\underset{\fZ^\reg_\fg} \otimes
\CL^{-\cla}_{\Fl^\cG_{w,\fZ}},$$
and
$$\CW^{\fZ++}_w:=\bigoplus_{\cmu,\cla\in \cLambda^+} \,
j_{\cla,*}\star
\on{Ind}^{\Hecke_\fZ}(\delta_{1,\Gr_G})
\underset{G[[t]]}\star
\CF_{V^\cmu}
\underset{\fZ^\reg_\fg} \otimes\CL^{-\cla-\cmu}_{\Fl^\cG_{w,\fZ}}.$$

Using the maps 
\begin{equation}  \label{exch map}
j_{\cla,*}\star
\on{Ind}^{\Hecke_\fZ}(\delta_{1,\Gr_G})\underset{G[[t]]}\star
\CF_{V^\cmu}\to j_{\cla+\cmu,*}\star
\on{Ind}^{\Hecke_\fZ}(\delta_{1,\Gr_G})
\end{equation}
and 
\begin{equation} \label{kappa w}
\kappa^{-,\cmu,w}:\CV^\cmu_\fZ\simeq
\Gamma(\Fl^\cG_{\fZ},\CL^{\cmu}_{\Fl^\cG_\fZ})\to
\CL^{\cmu}_{\Fl^\cG_{w,\fZ}},
\end{equation}
we obtain two maps
$\CW^{\fZ++}_w\to \CW^{\fZ+}_w$, 
and $\CW_w^\fZ$ is, by definition, their co-equalizer.

\medskip

By construction, the algebra $\Fun\bigl(\Fl^{\cG}_{w,\fZ}\bigr)$ acts
on $\CW_w^\fZ$ by endomorphisms. Moreover, for $\cla\in \cLambda^+$ we
have a canonical map
\begin{equation} \label{lambda inv of Wak}
j_{\cmu,*}\star \CW_w^\fZ\to \CL^{\cmu}_{\Fl^{\cG}_{w,\fZ}}
\underset{\Fun\bigl(\Fl^{\cG}_{w,\fZ}\bigr)}\otimes \CW_w^\fZ.
\end{equation}
By \eqref{lambda inv T}, this map is in fact an isomorphism.

\ssec{The universal property}   \label{univ ppty}

Let $\CM$ be an object of $\hg_\crit\mod^I$, endowed with an action of
$\Fun\bigl(\Fl^{\cG}_{w,\fZ}\bigr)$, compatible with the action of
$\fZ^\reg_\fg$, and a system of morphisms
\begin{equation} \label{lambda inv M}
j_{\cmu,*}\star \CM\to \CL^{\cmu}_{\Fl^{\cG}_{w,\fZ}}
\underset{\Fun\bigl(\Fl^{\cG}_{w,\fZ}\bigr)}\otimes \CM,
\end{equation}
for $\cmu\in \cLambda^+$,
such that the diagrams 
$$
\CD
j_{\cmu_1,*}\star j_{\cmu_2,*}\star \CM @>>>
j_{\cmu_1,*}\star (\CL^{\cmu_2}_{\Fl^{\cG}_{w,\fZ}}
\underset{\Fun\bigl(\Fl^{\cG}_{w,\fZ}\bigr)}\otimes \CM)  \\
@V{\sim}VV    @VVV   \\
j_{\cmu_1+\cmu_2,*}\star \CM  & &
\CL^{\cmu_2}_{\Fl^{\cG}_{w,\fZ}}
\underset{\Fun\bigl(\Fl^{\cG}_{w,\fZ}\bigr)}\otimes (j_{\cmu_1,*}
\star \CM) \\ @VVV     @VVV  \\
\CL^{\cmu_1+\cmu_2}_{\Fl^{\cG}_{w,\fZ}}
\underset{\Fun\bigl(\Fl^{\cG}_{w,\fZ}\bigr)}\otimes \CM 
@>{\sim}>>
\CL^{\cmu_2}_{\Fl^{\cG}_{w,\fZ}}
\underset{\Fun\bigl(\Fl^{\cG}_{w,\fZ}\bigr)}\otimes 
\CL^{\cmu_1}_{\Fl^{\cG}_{w,\fZ}}
\underset{\Fun\bigl(\Fl^{\cG}_{w,\fZ}\bigr)}\otimes \CM
\endCD
$$
are commutative for all $\cmu_1,\cmu_2\in \cLambda^+$.

\begin{lem}  \label{map from Wak}
Under the above conditions, specifying a map
$$\Gamma^{\Hecke_\fZ}(\Gr_G,\CW^\fZ_w)\to \CM,$$
which is compatible with the action of
$\Fun\bigl(\Fl^{\cG}_{w,\fZ}\bigr)$,
and which intertwines the morphisms \eqref{lambda inv of Wak} and 
\eqref{lambda inv M}, is equivalent to specifying a map
$\BV_\crit\to \CM$ such that the following diagram is commutative:
$$
\CD
\CF_{V^\cla}\star \BV_\crit @>{\sim}>>
\CV^\cla_\fZ\underset{\fZ^\reg_\fg}\otimes \BV_\crit
@>>> \CV^\cla_\fZ \underset{\fZ^\reg_\fg} \otimes \CM \\
@VVV  & & @V{\kappa^{-,\cla,w}}VV  \\
j_{\cla,*} \star \BV_\crit @>>> j_{\cla,*} \star \CM @>>>
\CL^{\cla}_{\Fl^{\cG}_{w,\fZ}}
\underset{\Fun\bigl(\Fl^{\cG}_{w,\fZ}\bigr)}\otimes \CM.
\endCD
$$
\end{lem}

\ssec{}   
Let us now specialize to the case $w=1$, and let us take
$\CM=\BW_{\crit,0}$.  We define an action of
$\Fun\bigl(\Fl^{\cG}_{1,\fZ}\bigr)$ on it via the tautological action of 
$\fH^\reg_{\crit}$ on $\BW_{\crit,0}$ and the isomorphism
$$\Spec(\fH^\reg_\crit)\simeq \Fl^{\cG}_{1,\fZ}.$$ 

The isomorphism \eqref{lambda inv M} is given
by \propref{twist of Wakimoto, local} and the identification
\begin{equation} \label{line bundle ident, fundamental}
\CL^\cla_{\fH,0}\simeq \CL^{\cla}_{\Fl^{\cG}_{1,\fZ}},
\end{equation}
which follows from the construction of $\on{map}^M_{\geom}$ in 
\secref{constr geom Miura}.

The fact that the conditions of \lemref{map from Wak} are satisfied, 
follows from the definitions.

\ssec{}      \label{sect ind and surj}

Thus, we have constructed a map
\begin{equation} \label{map to Wak}
\Gamma^{\Hecke_\fZ}(\Gr_G,\CW^\fZ_1)\to \BW_{\crit,0},
\end{equation}
and our present goal is to show that it is an isomorphism.

First, let us show that this map is injective. We have a natural map
$$\on{Ind}^{\Hecke_\fZ}(\delta_{1,\Gr_G})\underset{\fZ^\reg_\fg}
\otimes \Fun\bigl(\Fl^{\cG}_{1,\fZ}\bigr)\to \CW^\fZ_1,$$ where
$\on{Ind}^{\Hecke_\fZ}$ is as in \cite{FG4}, Sect. 2.5.  Moreover, by
using the same argument as in \cite{ABBGM}, Proposition 3.2.5, we find
that the object $\CW^\fZ_1\in \on{D}(\Gr_G)^{\Hecke_\fZ}_\crit\mod^I$ does
not have sub-objects that do not intersect the image of
$\on{Ind}^{\Hecke_\fZ}(\delta_{1,\Gr_G})\underset{\fZ^\reg_\fg}
\otimes \Fun\bigl(\Fl^{\cG}_{1,\fZ}\bigr)$.

Therefore, since the functor $\Gamma^{\Hecke_\fZ}$ is an equivalence,
to prove the injectivity of \eqref{map to Wak}, it suffices to show
that the composition
\begin{align*}
&\BV_\crit\underset{\fZ^\reg_\fg}\otimes
\Fun\bigl(\Fl^{\cG}_{1,\fZ}\bigr)\simeq
\Gamma^{\Hecke_\fZ}\Bigl(\Gr_G,\on{Ind}^{\Hecke_\fZ}(\delta_{1,\Gr_G})
\underset{\fZ^\reg_\fg} \otimes
\Fun\bigl(\Fl^{\cG}_{1,\fZ}\bigr)\Bigr)\to \\ &\to
\Gamma^{\Hecke_\fZ}(\Gr_G,\CW^\fZ_1)\to \BW_{\crit,0}.
\end{align*}

However, by construction, this composition is the canonical map
$$\BV_\crit\underset{\fZ^\reg_\fg}\otimes \fH^\reg_\fg\simeq
(\BW_{\crit,0})^{G[[t]]}\hookrightarrow
\BW_{\crit,0},$$ whose injectivity follows from \cite{F:wak}.

\ssec{}

The proof of surjectivity of the map \eqref{map to Wak} is similar to
that of \cite{FG4}, Proposition 4.13. Namely, we will show that for
$\cla\in \cLambda^+$ the map
$$j_{\cla+\crho,*}\star \BV_{\crit} \underset{\fZ^\reg_\fg}\otimes
\fH^\reg_\crit \to j_{\cla+\crho,*}\star \BW_{\crit,0}\simeq
\BW_{\crit,0}\underset{\fH^\reg_\crit} \otimes
\CL^{\cla+\crho}_{\fH,0}$$ is surjective.

The above map can be obtained by convolution with 
$j_{\cla\cdot w_0,*}$ from the map
\begin{equation} \label{map after rho}
j_{w_0\cdot \crho,*}\star \BV_\crit
\underset{\fZ^\reg_\fg}\otimes  \fH^\reg_\crit \to
j_{w_0\cdot \crho,*}\star \BW_{\crit,0}.
\end{equation}

Since the functor $j_{\cla\cdot w_0,*}$ is right-exact and sends
partially integrable representations to
$D^{<0}(\hg_\crit\mod_\reg)^I$, it is sufficient to show that the
cokernel of the map \eqref{map after rho} is partially integrable. To
establish the latter fact, by \cite{FG2}, Theorem 18.2.1, it is
sufficient to show that the map
$$H^\semiinf\bigl(\fn^-\ppart,t\fn[[t]],(j_{w_0\cdot \crho,*} \star
\BV_\crit)\otimes \Psi_{-\crho}\bigr)\underset{\fZ^\reg_\fg}\otimes
\fH^\reg_\crit\to H^\semiinf\bigl(\fn^-\ppart,t\fn[[t]],(j_{w_0\cdot
\crho,*}\star \BW_{\crit,0}) \otimes\Psi_{-\crho}\bigr)$$ is a
surjection. However, by \cite{FG2}, 18.1.1, the above map of
semi-infinite cohomologies identifies with
$$H^\semiinf(\fn\ppart,\fn[[t]],\BV_\crit\otimes \Psi_0)
\underset{\fZ^\reg_\fg}\otimes  \fH^\reg_\crit\to 
H^\semiinf(\fn\ppart,\fn[[t]],\BW_{\crit,0}\otimes \Psi_0).$$

According to \cite{FB}, Theorem 15.1.9, and \cite{FG2},
Theorem 18.2.4, the latter map becomes the identity isomorphism under
the natural identification of both sides with $\fH^\reg_\crit$. This
completes the proof of \thmref{ident prince Wak}.

\ssec{}   \label{another way}

We will now discuss yet another way of obtaining the module
$\BW_{\crit,0}$.  Consider the $\hg_\crit$-modules
\begin{equation} \label{dir sys}
j_{\cla,*}\star \BV_{\crit}\otimes \fl^\cla_{-\rho}
\end{equation}
for $\cla\in \cLambda^+$.

We claim that whenever $\cla_1-\cla_2=\cmu\in \cLambda^+$ there exists
a natural map
$$j_{\cla_1,*}\star \BV_{\crit} \otimes \fl^{\cla_1}_{-\rho}
\to j_{\cla_2,*}\star \BV_{\crit} \otimes \fl^{\cla_2}_{-\rho}.$$

Indeed, we have a map of D-modules
\begin{equation} \label{act ingr}
j_{\cla_1,*}\star \CF_{V^{\cmu}}\to
j_{\cla_1,*}\star j_{\cmu,\Gr_G,*}\simeq j_{\cla_2,\Gr,*},
\end{equation}
and we compose it with the map
$$\BV_\crit\otimes \fl^\cmu_{\rho}
\overset{\on{id}\otimes \kappa^{\cmu}}
\longrightarrow \BV_\crit\underset{\fZ^\reg_\fg}\otimes
\CV^\cmu_\fZ\simeq \CF_{V^{\cmu}}\star \BV_\crit.$$

{}From the fact that the maps $\kappa^{\cmu}$ satisfy the Pl\"ucker
relations, it follows that the maps \eqref{dir sys} form a directed
system.

\begin{thm}  \label{Wak dir lim}
There exists a canonical isomorphism
$$\underset{\longrightarrow}{\lim}\, 
j_{\cla,*}\star \BV_{\crit}\otimes
\fl^\cla_{-\rho} \to \BW_{\crit,0}.$$
\end{thm}

The rest of this section will be devoted to the proof of this
theorem. First, we construct a map from the LHS to the RHS. By
definition, this amounts to a compatible system of maps
$$j_{\cla,*}\star \BV_{\crit}\otimes
\fl^\cla_{-\rho} \to \BW_{\crit,0}.$$
These maps are given by
$$j_{\cla,*}\star \BV_{\crit}\otimes
\fl^\cla_{-\rho} \to
j_{\cla,*}\star \BW_{\crit,0}\otimes
\fl^\cla_{-\rho}\to \BW_{\crit,0},$$
where the first arrow comes by convolution from the canonical map
$\phi:\BV_\crit\to \BW_{\crit,0}$, and the second map is the isomorphism
of \propref{twist of Wakimoto, local}.

The fact that these maps are compatible for different $\cla$ follows
from \propref{crucial comp}.

\medskip

To construct the map from the RHS to the LHS in \thmref{Wak dir lim}
we will use \lemref{map from Wak}. For that we need to endow the LHS
with an action of the algebra $\Fun(\Fl^{\cG}_{1,\fZ})$. Note that
this algebra is canonically isomorphic to the direct limit
$$\underset{\longrightarrow}{\lim}\, \CV^\cla_\fZ\otimes
\fl^\cla_{-\rho},$$ where the transition maps are
defined whenever $\cla_2-\cla_1=\cmu\in \cLambda^+$ and are equal to
$$\CV^{\cla_1}_\fZ\otimes \fl^{\cla_1}_{-\rho}\to
\CV^{\cla_1}_\fZ\underset{\fZ^\reg_\fg} \otimes \CV^\cmu_\fZ  \otimes 
\fl^\cmu_{-\rho}
\otimes \fl^{\cla_1}_{-\rho}\to
\CV_\fZ^{\cla_2}\otimes \fl^{\cla_2}_{-\rho},$$
where the first arrow corresponds to the map 
$\kappa^{\cmu}:\fl^\cmu_{\rho}\to \CV^\cmu_\fZ$.

We define the desired action of $\Fun(\Fl^{\cG}_{1,\fZ})$ on the
direct limit appearing in \thmref{Wak dir lim} by means of
$$\CV^\cmu_\fZ\underset{\fZ^\reg_\fg}\otimes \bigl(j_{\cla,*}\star
\BV_{\crit}\otimes \fl^\cla_{-\rho} \bigr)\otimes
\fl^\cmu_{-\rho} \simeq j_{\cla,*}\star
\CF_{V^\cmu}\underset{G[[t]]}\star \BV_{\crit}\otimes 
\fl^{\cla+\cmu}_{-\rho}\to j_{\cla+\cmu,*}\star \BV_{\crit}\otimes
\fl^{\cla+\cmu}_{-\rho}.$$

The data of morphisms \eqref{lambda inv M} for the RHS of \thmref{Wak
dir lim} is evident from the definitions. It is straightforward to
check that the diagram appearing in \lemref{map from Wak} is
commutative, thereby giving rise to a map from the RHS to the LHS in
\thmref{Wak dir lim}. Moreover, it is easy to check that the two maps
constructed above are mutually inverse.

\section{Wakimoto modules attached to other Schubert cells}
\label{other Wakimoto}

In this section we will prove \thmref{ident w Wak}. We will derive it
from \thmref{ident prince Wak} using the technique of chiral modules
over chiral algebras.

\ssec{Plan of the proof}   \label{Wak w}     \label{w schemes}

We need to establish an isomorphism
\begin{equation} \label{need to show w}
\Gamma^{\Hecke_\fZ}(\CW^\fZ_w) \simeq
\BW_{\crit,w(\rho)-\rho,\reg}\otimes \fl^{w(\crho)-\crho}_{w(\rho)}
\end{equation}
for all $w \in W$. 
To construct a map in one direction (from left to right)
we will use the universal property of $\CW^\fZ_w$, described in
\lemref{map from Wak}.

\medskip

First, we need to make the algebra $\Fun(\Fl^{\cG}_{w,\fZ})$ act on
$\BW_{\crit,w(\rho)-\rho,\reg}$. This results from
the isomorphism \eqref{H RS}. 

\medskip

Next, we need to establish that the isomorphism \eqref{lambda inv M}
holds for $\CM=\BW_{\crit,w(\rho)-\rho,\reg}$. On the one hand, by
\propref{twist of Wakimoto, local},
$$j_{\cla,*}\star \BW_{\crit,w(\rho)-\rho}\simeq
\BW_{\crit,w(\rho)-\rho}\underset{\fH_\crit^{\RS,w_0(w(\rho)-\rho)}}
\otimes \CL^\cla_{\fH,w(\rho)-\rho},$$
where 
\begin{equation}  \label{CLa const}
\CL^\cla_{\fH,w(\rho)-\rho}\simeq \fH_\crit^{\RS,w_0(w(\rho)-\rho)}\otimes
\fl^\cla_{w(\rho)}.
\end{equation}
This implies that 
$$j_{\cla,*}\star \BW_{\crit,w(\rho)-\rho,\reg}\simeq
\BW_{\crit,w(\rho)-\rho,\reg}
\underset{(\fH_\crit^{\RS,w_0(w(\rho)-\rho)}\underset{\fZ^\nilp_\fg}\otimes
\fZ^\reg_\fg)} \otimes \CL^\cla_{\fH,w,\reg},$$ where we denote
$$\CL^\cla_{\fH,w,\reg}:=\CL^\cla_{\fH,w(\rho)-\rho}
\underset{\fZ^\nilp_\fg}\otimes \fZ^\reg_\fg.$$

By \secref{univ ppty}, we need to show that the isomorphism
$\Fl^{\cG}_{w,\fZ}\simeq \Spec\Bigl(\fH_\crit^{\RS,w_0(w(\rho)-\rho)}
\underset{\fZ^\nilp_\fg}\otimes \fZ^\reg_\fg\Bigr)$ of \eqref{H RS}
lifts to an isomorphism of $\cH$-torsors
\begin{equation} \label{ident tors}
\{\cla\mapsto \CL^\cla_{\Fl^{\cG}_{w,\fZ}}\}
\overset{\{\gamma^{\cla}\}}\Longleftrightarrow \{\cla\mapsto
\CL^\cla_{\fH,w,\reg}\}.
\end{equation}

By \eqref{CLa const}, the latter amounts to an isomorphism
of $\cH$-torsors
\begin{equation} \label{ident tors, special}
\{\cla\mapsto \CL^\cla_{\Fl^{\cG}_{w,\fZ}}\} \text{ and }
\{\cla\mapsto \CO_{\Fl^{\cG}_{w,\fZ}}\otimes \fl^\cla_{w(\rho)}\}.
\end{equation}

\ssec{An identification of $\cH$-torsors}   \label{H torsors}

We fill first show that there exists {\it some} isomorphism as in
\eqref{ident tors, special} that respects the $\Aut(\D)$-actions.

Recall the D-scheme 
\begin{equation} \label{conn with pole}
\ConnX^{\RS,\rho-w(\rho)}\underset{\Op^\nilp_{\cg,X}}\times \Op_{\cg,X}
\end{equation}
(see \secref{part cases}). Its restriction to $X-x$ is isomorphic to
$\ConnX$, and its fiber over $x$ is isomorphic to $\MOp^{w,\reg}_\cg$,
once we identify $\D_x$ with $\D$.

For $\cla\in \cLambda$, let us denote by $\CL^\cla_{\MOp,w,X}$ the
line bundle on \eqref{conn with pole} equal to the pull-back from $X$
of the line bundle
$\omega_X^{\langle\rho,\cla\rangle}\bigl(\langle\rho-w(\rho),
\cla\rangle\cdot x)\bigr)$.

We claim that this line bundle identifies with the pull-back of the
line bundle $\CL^\cla_{\MOp_{\cg,X}}$ under the map \eqref{RS
extension}.  This follows from the fact that the corresponding
isomorphism holds tautologically over $X-x$, and that both line
bundles have connections that are regular at $x$.

\medskip

Restricting the above isomorphism of line bundles to $x\in X$, we
obtain that $\CL^\cla_{\MOp_\cg^{w,\reg}}$ identifies with
$\CO_{\MOp_\cg^{w,\reg}}\otimes \omega_x^{\langle
w(\rho),\cla\rangle}$.  The latter isomorphism respects the
$\Aut(\D)$-actions, where $\Aut(\D)$ acts on $\omega_x^{\langle
w(\rho),\cla\rangle}$ via the isomorphism $\Aut(\D)\simeq \Aut(\D_x)$,
corresponding to the above choice of an isomorphism $\D\simeq
\D_x$. Finally, let us observe that the character of $\Aut(\D)$ on
$\omega_x^{\langle w(\rho),\cla\rangle}$ equals that of
$\fl^\cla_{w(\rho)}$.

\medskip

Let us denote by $\fH_{\crit,w,X}$ the D-algebra
$$\fH_{\crit,w,X}:=\fH_{\crit,X}^{\RS,w_0(w(\rho)-\rho)}
\underset{\fz_{\fg,X}^\nilp}\otimes \fz_{\fg,X}.$$
Via $\on{map}^M_{\geom}$ it identifies with the algebra of regular functions
on the D-scheme \eqref{conn with pole}. For $\cla\in \cLambda$, let
$\CL^\cla_{\fH,w,X}$ denote the corresponding line bundle with
a {\it regular} connection along $X$ over it.

The restriction of $\CL^\cla_{\fH,w,X}$ to $X-x$ identifies with
$\CL^\cla_{\fH,X}$.  The restriction to the fiber over $x$,
identifies, as a line bundle over $\Fl^\cG_{w,\fZ}$ via \eqref{H RS},
with $\CL^\cla_{\Fl^\cG_{w,\fZ}}$.  By the above discussion, we obtain
an identification
$$\CL^\cla_{\Fl^\cG_{w,\fZ}} \simeq \CO_{\Fl^\cG_{w,\fZ}}\otimes 
\omega_x^{\langle w(\rho),\cla\rangle},$$
as required in \eqref{ident tors, special}.

\medskip

For future reference note that for $\cla\in \cLambda^+$ we have a
canonical map
$$\kappa^{-,\cla,w}_\fH:\CV_{\fZ,X}\underset{\fz_{\fg,X}}\otimes
\fH_{\crit,w,X}\to \CL^\cla_{\fH,w,X}.$$

Its restriction to $X-x$ equals the map $\kappa^{-,\cla}_\fH$ of
\eqref{des Miura}, and its restriction to $x$ coincides with the
canonical map \eqref{kappa w}.

\ssec{}

According to \secref{univ ppty}, the next step in the definition of
the map in \eqref{need to show w} is construction of a map
\begin{equation} \label{Vak to Wak}
\BV_{\crit}\underset{\fZ^\reg_\fg}\otimes \Fun(\Fl^{\cG}_{w,\fZ})\to
\BW_{\crit,w(\rho)-\rho,\reg}\otimes \omega_x^{\langle
\rho-w(\rho),\crho\rangle}.
\end{equation}
Here $\omega_x^{\langle \rho-w(\rho),\crho\rangle}$ is the
$1$-dimensional $\Aut(\D)$-module as in \secref{H torsors}. It is
non-canonically isomorphic to $\fl^{w(\crho)-\crho}_{w(\rho)}$.
Eventually, we will fix an isomorphism \eqref{ident tors}, thereby
fixing a choice for the above isomorphism of lines as well.

\medskip

Finding a morphism as in \eqref{Vak to Wak} is equivalent to
exhibiting a map 

\begin{equation} \label{av Wak}
\BV_\crit\underset{\fZ^\reg_\fg}\otimes \Fun(\Fl^{\cG}_{w,\fZ})\to 
\on{Av}_{G[[t]]/I}\bigl(\BW_{\crit,w(\rho)-\rho,\reg}\bigr)
\otimes \omega_x^{\langle \rho-w(\rho),\crho\rangle}
\end{equation}
(see \cite{FG2}, Sect. 20.2, for the definition of the averaging
functor $\on{Av}_{G[[t]]/I}$).

By the definition of $\BW_{\crit,w(\rho)-\rho,\reg}$ and using the fact
that the $\hg_\crit$-module $\BW_{\crit,w(\rho)-\rho}$ is flat over
$\fH_\crit^{\RS,w_0(w(\rho)-\crho)}$, the expression in the RHS of
\eqref{av Wak} can be rewritten as
$$\on{Av}_{G[[t]]/I}\bigl(\BW_{\crit,w(\rho)-\rho}\bigr)
\underset{\fH_\crit^{\RS,w_0(w(\rho)-\crho)}}{\overset{L}\otimes}
\Fun(\Fl^{\cG}_{w,\fZ})\otimes 
\omega_x^{\langle\rho-w(\rho),\crho\rangle}.$$

By construction of Wakimoto modules,
\begin{equation} \label{av of diff}
\on{Av}_{G[[t]]/I}(\BW_{\crit,w(\rho)-\rho})\simeq
H^\semiinf(\fn\ppart,\fn[[t]],\fD^{ch}(G)_{\crit,x})
\overset{L}{\underset{\fh[[t]];\fh}\otimes} \BC^{w_0(w(\rho)-\rho)},
\end{equation}
where $\BC^{w_0(w(\rho)-\rho)}$ is the corresponding character of
$\fh$.

\medskip

Consider the component of the expression appearing on the RHS of the above
formula that has degree zero with respect to $\BG_m\subset
\Aut(\D)$. It identifies with the $w_0(w(\rho)-\rho)$-weight space in the 
Lie algebra cohomology $H^\bullet(\fn,\Fun(G))$.

By the Bott-Borel Weil theorem, the latter is concentrated in the
cohomological degree $\ell(w)$, and is canonically isomorphic to
$\BC$. 

\medskip

The corresponding vector in the $\ell(w)$-th cohomology of
$\on{Av}_{G[[t]]/I}(\BW_{\crit,w(\rho)-\rho})$ is easily seen to be
$\fg[[t]]$-invariant, so we obtain a map
$$\BV_\crit\to
h^{\ell(w)}\left(\on{Av}_{G[[t]]/I}(\BW_{\crit,w(\rho)-\rho})\right),$$
and, hence, a map
\begin{multline}   \label{compl map}
\left(\BV_\crit \underset{\fZ_\fg^\reg}\otimes \Fun(\Fl^{\cG}_{w,\fZ})\right)
\underset{\fH_\crit^{\RS,w_0(w(\rho)-\crho)}}{\overset{L}\otimes}
\Fun(\Fl^{\cG}_{w,\fZ})\otimes 
\omega_x^{\langle\rho-w(\rho),\crho\rangle}[-\ell(w)]\to \\
\to \on{Av}_{G[[t]]/I}\bigl(\BW_{\crit,w(\rho)-\rho}\bigr)
\underset{\fH_\crit^{\RS,w_0(w(\rho)-\crho)}}{\overset{L}\otimes}
\Fun(\Fl^{\cG}_{w,\fZ})\otimes 
\omega_x^{\langle\rho-w(\rho),\crho\rangle}.
\end{multline}

We claim that the $0$th cohomology of the 
LHS of the expression in \eqref{compl map} is canonically isomorphic to
$\BV_\crit\underset{\fZ^\reg_\fg}\otimes \Fun(\Fl^{\cG}_{w,\fZ})$,
which would produce the desired (non-zero!) map in \eqref{av Wak}.

To establish this isomorphism, it would be enough to show that
\begin{equation} \label{one Tor}
\on{Tor}^{\ell(w)}_{\fH_\crit^{\RS,w_0(w(\rho)-\rho)}}
\bigl(\Fun(\Fl^{\cG}_{w,\fZ}), \Fun(\Fl^{\cG}_{w,\fZ})\bigr)\simeq 
\Fun(\Fl^{\cG}_{w,\fZ})\otimes 
\omega_x^{\langle w(\rho)-\rho,\crho\rangle}.
\end{equation}

We can identify the LHS of \eqref{one Tor} with
$\Lambda^{\ell(w)}\Bigl(N^*_{\Fl^{\cG}_{w,\fZ}/
\Spec(\fH_\crit^{\RS,w_0(w(\rho)-\rho)})}\Bigr)$, and the RHS
with $\CL^{\crho-w(\crho)}_{\Fl^\cG_{\fZ,w}}$, using 
\secref{H torsors}. Hence, \eqref{one Tor} is
equivalent to an isomorphism of line bundles
$$\Lambda^{\ell(w)}\Bigl(N^*_{\Fl^{\cG}_{w,\fZ}/
\Spec(\fH_\crit^{\RS,w_0(w(\rho)-\rho)})}\Bigr)\simeq
\CL^{\crho-w(\crho)}_{\Fl^\cG_{\fZ,w}}$$
over $\Fl^\cG_{\fZ,w}$. Using the diagram 
$$ 
\CD 
\Fl^{\cG}_{w,\fZ} @>>>
\Spec\left(\fH_\crit^{\RS,w_0(w(\rho)-\rho)}\right) \\ 
@V{\sim}VV @V{\sim}VV \\ 
\MOp_\cg^{w,\reg} @>>> \MOp_\cg^{w,\nilp}, 
\endCD
$$
we can translate the existence of the above isomorphism to
that between the line bundles
$$\Lambda^{\ell(w)}\Bigl(N^*_{\MOp_{\cg}^{w,\reg}/
\MOp_{\cg}^{w,\nilp}}\Bigr)\simeq
\CL^{\crho-w(\crho)}_{\MOp_{\cg}^{w,\reg}}$$ over
$\MOp_{\cg}^{w,\reg}$. The latter isomorphism follows from \cite{FG2},
Corollary 3.6.3.

\medskip

We will need the following property of the map \eqref{Vak to Wak}:

\begin{prop}  \label{Vak to Wak, semiinf} \hfill 

\smallskip

\noindent{\em (1)} The map 
\begin{align*}
&H^\semiinf\left(\fn\ppart,\fn[[t]], \BV_{\crit}\otimes \Psi_0\right)
\underset{\fZ^\reg_\fg}\otimes \Fun(\Fl^{\cG}_{w,\fZ})\to \\
&H^\semiinf\left(\fn\ppart,\fn[[t]],
\BW_{\crit,w(\rho)-\rho,\reg}\otimes \Psi_0\right)\otimes
\omega_x^{\langle \rho-w(\rho),\crho\rangle},
\end{align*}
induced by \eqref{Vak to Wak}, is an isomorphism.

\smallskip

\noindent{\em (2)} 
The map
\begin{align*}
&H^\semiinf\left(\fn^-\ppart,t\fn^-[[t]], j_{w_0\cdot\crho,*}\star
\BV_{\crit}\otimes \Psi_{-\crho}\right)
\underset{\fZ^\reg_\fg}\otimes \Fun(\Fl^{\cG}_{w,\fZ})\to  \\
&H^\semiinf\left(\fn^-\ppart,t\fn^-[[t]], j_{w_0\cdot\crho,*}\star
\BW_{\crit,w(\rho)-\rho,\reg}\otimes \Psi_{-\crho}\right)
\otimes \omega_x^{\langle \rho-w(\rho),\crho\rangle}
\end{align*}
is also an isomorphism.

\end{prop}

The proof of this proposition will be given in 
\secref{proof of Vak to Wak, semiinf}.

\ssec{}

By \lemref{map from Wak}, in order to complete the construction of the map
in \eqref{need to show w} we need to choose an isomorphism
in \eqref{ident tors} so that the following diagram becomes commutative
for every $\cla\in \cLambda^+$:

\begin{equation} \label{CD w}
\CD \BV_{\crit}\underset{\fZ^\reg_\fg}\otimes \CV^\cla_{\fZ}
@>>> \BW_{\crit,w(\rho)-\rho,\reg}\underset{\fZ^\reg_\fg}\otimes
\CV^\cla_{\fZ}\otimes \omega_x^{\langle\rho-w(\rho),\crho\rangle} \\ 
@V{\sim}VV @V{\kappa^{-,\cla,w}}VV \\
\CF_{V^\cla}\underset{G[[t]]}\star \BV_{\crit} & &
\BW_{\crit,w(\rho)-\rho,\reg}\underset{\Fun(\Fl^{\cG}_{w,\fZ})}
\otimes \CL^\cla_{\Fl^{\cG}_{w,\fZ}} \otimes \omega_x^{\langle
\rho-w(\rho),\crho\rangle} \\ @VVV @V{\sim}V{\gamma^\cla}V \\
j_{\cla,\Gr_G,*}\underset{G[[t]]}\star \BV_{\crit} & &
\BW_{\crit,w(\rho)-\rho,\reg}\underset{\Fun(\Fl^{\cG}_{w,\fZ})}
\otimes \CL^\cla_{\fH,w,\reg} \otimes \omega_x^{\langle
\rho-w(\rho),\crho\rangle} \\ 
@V{\sim}VV @V{\sim}VV \\ 
j_{\cla,*}\star \BV_{\crit} @>>> 
j_{\cla,*}\star \BW_{\crit,w(\rho)-\rho,\reg} \otimes
\omega_x^{\langle \rho-w(\rho),\crho\rangle}, 
\endCD
\end{equation}
where $\kappa^{-,\cla,w}$ is as in \eqref{kappa w}.

To construct the isomorphism \eqref{ident tors}, we will realize
$\BW_{\crit,w(\rho)-\rho,\reg} \otimes 
\omega_x^{\langle \rho-w(\rho),\crho\rangle}$ as a fiber at 
$x\in X$ of a certain chiral
$\CA_{\fg,\crit,X}$-module.

\ssec{}  \label{chiral modules}

We will use the following general construction. Let $\CA$ be a chiral
algebra, let $\CM$ be a torsion-free chiral $\CA$-module on $X-x$, and
let $\CN_1,\CN_2$ be two chiral $\CA$-modules, supported at $x$.  Let
\begin{equation}  \label{ch pairing}
\jmath_*(\CM)\otimes \CN_1\to \imath_!(\CN_2)
\end{equation}
be a chiral pairing (see \cite{CHA} and \cite{FG3}, Sect. 2.1). Here
$\jmath$ and $\imath$ denote the embeddings $X-x\hookrightarrow X$ and
$x\hookrightarrow X$, respectively.

Let $\bv\in \CN_1$ be a vector, which is annihilated by the Lie-*
action of $\CA$, and such that the resulting map $\jmath_*(\CM)\to
\imath_!(\CN_2)$ is surjective.  Then
$$\CM':=\on{ker}\bigl(\jmath_*(\CM)\to \imath_!(\CN_2)\bigr)$$ is a
chiral $\CA$-module, whose fiber at $x$ is $\CN_2$.

\ssec{}

We apply the above construction in the following situation. We let
$\CM$ be $\BW_{\crit,X}$. Recall (see \cite{F:wak} and \cite{FG2},
Sect. 10.3), that by the construction of Wakimoto modules,
$\BW_{\crit,X}$ is in fact a chiral algebra that contains
$\CA_{\fg,\crit,X}$ as a chiral subalgebra.

We let $\CN_1=\CN_2:=\BW_{\crit,w(\rho)-\rho,\reg}\otimes
\omega_x^{\langle \rho-w(\rho),\crho\rangle}$,
which is naturally a chiral $\BW_{\crit,X}$-module, supported
at $x\in X$. Finally, we let $\bv$ to be the image of the 
vacuum vector in $\BV_{\crit}$ under the map \eqref{Vak to Wak}.

\begin{lem}
The resulting map
$$\jmath_*(\BW_{\crit,W})\to
\imath_!\Bigl(\BW_{\crit,w(\rho)-\rho,\reg}\otimes \omega_x^{\langle
\rho-w(\rho),\crho\rangle}\Bigr)$$ is surjective.
\end{lem}

\begin{proof}

{}From the construction of Wakimoto modules we obtain that
$\BW_{\crit,X}$-submodules of $\BW_{\crit,w(\rho)-\rho,\reg}\otimes
\omega_x^{\langle \rho-w(\rho),\crho\rangle}$ are in bijection
with $\fH_{\crit,X}$-submodules of $\Fun(\Fl^{\cG}_{w,\fZ})\otimes 
\omega_x^{\langle \rho-w(\rho),\crho\rangle}$, and the correspondence
is given by applying the functor 
$H^\semiinf\left(\fn\ppart,\fn[[t]], ?\otimes \Psi_0\right)$.

Hence, the assertion of the lemma follows from \propref{Vak to Wak,
semiinf}(1).

\end{proof}

Let us denote the resulting chiral $\CA_{\fg,\crit,X}$-module by
$\BW_{\crit,w,X}$.  By construction, it comes equipped with a map of
chiral $\CA_{\fg,\crit,X}$-modules $\BV_{\crit,X}\to \BW_{\crit,w,X}$.

\medskip

By construction, we have the following
commutative diagram:
$$
\CD
\fz_{\fg,X} @>>> \BV_{\crit,X} \\
@VVV  @VVV  \\
\fH_{\crit,w,X} @>>> \BW_{\crit,w,X},
\endCD
$$
whose restriction to $X-x$ is
$$
\CD
\fz_{\fg,X-x} @>>> \BV_{\crit,X-x}  \\
@V{\varphi}VV   @V{\phi}VV  \\
\fH_{\crit,X-x} @>>> \BW_{\crit,X-x},
\endCD
$$
and whose fiber at $x$ is the diagram
$$
\CD
\fZ^\reg_\fg @>>> \BV_{\crit} \\
@VVV  @VV{\text{\eqref{Vak to Wak}}}V  \\
\Fun(\Fl^{\cG}_{w,\fZ}) @>>> \BW_{\crit,w(\rho)-\rho,\reg}\otimes 
\omega_x^{\langle \rho-w(\rho),\crho\rangle}.
\endCD
$$

\ssec{}

For $\cla\in \cLambda^+$ consider the chiral module $j_{\cla,*,X}\star
\BW_{\crit,w,X}$. We claim that this convolution is concentrated in 
the cohomological degree $0$. 

Indeed, {\it a priori}, it is concentrated
in non-positive cohomolgical degrees, since the functor $j_{\cla,*,X}\star ?$
is right-exact. Now, it does not have strictly negative cohomologies,
because this is true for both chiral $\CA_{\fg,\crit,X}$-modules
$\jmath_*(\BW_{\crit,X})$ and 
$\imath_!\bigl(\BW_{\crit,w(\rho)-\rho,\reg} 
\otimes \omega_x^{\langle \rho-w(\rho),\crho\rangle}\bigr)$.

\medskip

Recall the $\cH$-torsor $\{\cla\mapsto \CL^\cla_{\fH,w,X}\}$ on 
$\Spec(\fH_{\crit,w,X})$ (see \secref{H torsors}).

\begin{prop}   \label{extension of maps}  \hfill

\smallskip

\noindent{\em (1)} The map $\CL^\cla_{\fH,X-x}\to j_{\cla,X-x,*}\star
\BW_{\crit,X-x}$, equal to the composition
$$\CL^\cla_{\fH,X-x}\hookrightarrow
\CL^\cla_{\fH,X-x}\underset{\fH_{\crit,X-x}}\otimes
\BW_{\crit,X-x}\simeq j_{\cla,X-x,*}\star\BW_{\crit,X-x},$$ extends to
a map $\CL^\cla_{\fH,w,X}\to j_{\cla,*,X}\star \BW_{\crit,w,X}$.

\smallskip

\noindent{\em (2)}
The resulting map on the level of fibers at $x$
\begin{align*}
&\CL^\cla_{\Fl^{\cG}_{w,\fZ}}\to j_{\cla,*}\star
\BW_{\crit,w(\rho)-\rho,\reg}\otimes \omega_x^{\langle
\rho-w(\rho),\crho\rangle}\simeq \\ &\simeq
\CL^\cla_{\fH,w,\reg}\underset{\Fun(\Fl^{\cG}_{w,\fZ})}\otimes
\BW_{\crit,w(\rho)-\rho,\reg} 
\otimes \omega_x^{\langle \rho-w(\rho),\crho\rangle}
\end{align*}
comes from an {\em isomorphism} onto $\CL^\cla_{\fH,w,\reg}$, followed
by the map 
\begin{multline*}
\CL^\cla_{\fH,w,\reg}\hookrightarrow
 \CL^\cla_{\fH,w,\reg}\underset{\Fun(\Fl^\cG_{\fZ,w})}\otimes
\Bigl(\Fun(\Fl^\cG_{\fZ,w})\underset{\fZ^\reg_\fg}\otimes \BV_\crit\Bigr)
\overset{\text{\eqref{Vak to Wak}}}\longrightarrow \\
\to \CL^\cla_{\fH,w,\reg}\underset{\Fun(\Fl^\cG_{\fZ,w})}\otimes
\BW_{\crit,w(\rho)-\rho,\reg} 
\otimes \omega_x^{\langle \rho-w(\rho),\crho\rangle}.
\end{multline*}
\end{prop}

Let us assume this proposition and finish the construction of the map
in \eqref{need to show w}. First, note that the second assertion of
the proposition defines the sought-after identification $\gamma^\cla$
of \eqref{ident tors}.

\medskip

Next, observe that the first assertion of the proposition implies
that there exists a diagram:
\begin{equation}   \label{CD W(2)}
\CD \CF_{V^\cla,X}\star \BV_{\crit,X} @<{\sim}<<
\BV_{\crit,X}\underset{\fz_{\fg,X}}\otimes \CV_{\fZ,X} \\ @VVV @VVV \\
j_{\cla,*,X}\star \BV_{\crit,X} & &
\BV_{\crit,X}\underset{\fz_{\fg,X}}\otimes \CV_{\fZ,X}
\underset{\fz_{\fg,X}}\otimes \fH_{\crit,w,X} \\ @VVV
@V{\kappa_\fH^{-,\lambda,w}}VV \\ j_{\cla,*,X}\star \BW_{\crit,w,X}
@<<< \BV_{\crit,X}\underset{\fz_{\fg,X}}\otimes \CL^\cla_{\fH,w,X},
\endCD
\end{equation}
which is commutative, since it extends the corresponding commutative
diagram over $X-x$.

Let us consider the fiber of the above diagram over $x$. 
The two resulting maps
$$\BV_{\crit,X}\underset{\fz_{\fg,X}}\otimes
\CV_{\fZ,X}\rightrightarrows
\BW_{\crit,w(\rho)-\rho,\reg}\underset{\Fun(\Fl^{\cG}_{w,\fZ})}\otimes
\CL^\cla_{\fH,w,\reg} \otimes \omega_x^{\langle
\rho-w(\rho),\crho\rangle}$$ are equal to the two circuits in the
diagram \eqref{CD w}.  Thus, we obtain a well-defined map 
\eqref{need to show w}.

\medskip

The proof that this map is an isomorphism uses the same argument as
the one used in the proof of \thmref{ident prince Wak} in
\secref{sect ind and surj}:

The surjectivity assertion follows from \propref{Vak to Wak,
semiinf}(2) as in {\it loc. cit.}  To prove the injectivity, it is
enough to show that the map \eqref{Vak to Wak} itself is injective.

However, according to \cite{FG1}, the functor
$$
\fZ^\reg_\fg\mod \to \hg_\crit\mod^{G[[t]]}_\reg, \qquad \CT \mapsto
\BV_{\crit}\underset{\fZ^\reg_\fg}\otimes \CT
$$
is an equivalence of categories. Therefore for any
$\fZ^\reg_\fg$-module $\CT$, the submodules of
$\BV_{\crit}\underset{\fZ^\reg_\fg}\otimes \CT$ are in bijection with
the $\fZ^\reg_\fg$-submodules of $\CT$ itself, and any such submodule is
determined by its image in
$$\CT\simeq H^\semiinf\Bigl(\fn\ppart,\fn[[t]],
(\BV_{\crit}\underset{\fZ^\reg_\fg}\otimes \CT)\otimes \Psi_0\Bigr).$$
Therefore, the injectivity of \eqref{Vak to Wak} follows from
\propref{Vak to Wak, semiinf}(1).

This completes the proof of \thmref{ident w Wak} modulo \propref{Vak
to Wak, semiinf} and \propref{extension of maps}. The remainder of
this section is devoted to the proof these two propositions.

\ssec{Proof of \propref{Vak to Wak, semiinf}} 
\label{proof of Vak to Wak, semiinf}

First, let us notice that assertions (1) and (2) of the proposition
are equivalent by \cite{FG2}, Proposition 18.1.1.

Secondly, both the LHS and the RHS, appearing in \propref{Vak to Wak,
semiinf}(2), are isomorphic to $\Fun(\Fl^{\cG}_{w,\fZ})\otimes
\omega_x^{\langle\rho,\crho\rangle}$, as graded
$\Fun(\Fl^{\cG}_{w,\fZ})$-modules. Since the degree $0$ component of
$\Fun(\Fl^{\cG}_{w,\fZ})$ consists of scalars, it sufficient to show
that the map in \propref{Vak to Wak, semiinf}(1) is non-zero.

Let us assume the contrary. Then, by \cite{FG2}, Theorem 18.2.1,
the map
$$j_{w_0\cdot\crho,*}\star \BV_\crit\to
j_{w_0\cdot\crho,*}\star \BW_{\crit,w(\rho)-\rho,\reg}\otimes 
\omega_x^{\langle \rho-w(\rho),\crho\rangle}$$
has a partially integrable image. Let us denote it by $\CM$, and
consider the composition
$$j_{\cla\cdot w_0,*}\star j_{w_0\cdot\crho,*}\star \BV_\crit\to
j_{\cla\cdot w_0,*}\star \CM\to j_{\cla\cdot w_0,*}\star
j_{w_0\cdot\crho,*}\star \BW_{\crit,w(\rho)-\rho,\reg}\otimes 
\omega_x^{\langle \rho-w(\rho),\crho\rangle}$$
for $\cla\in \cLambda^+$. On the one hand, 
$j_{\cla\cdot w_0,*}\star j_{w_0\cdot\crho,*}\star \BV_\crit\simeq
j_{\cla,*}\star \BV_\crit$ and 
$$j_{\cla\cdot w_0,*}\star j_{w_0\cdot\crho,*}\star 
\BW_{\crit,w(\rho)-\rho,\reg}\otimes 
\omega_x^{\langle \rho-w(\rho),\crho\rangle}\simeq 
j_{\cla,*}\star \BW_{\crit,w(\rho)-\rho,\reg}\otimes 
\omega_x^{\langle \rho-w(\rho),\crho\rangle},$$
and the above composed map is non-zero, since the functor 
$j_{\cla,*}\star ?$ is invertible on the derived category.

\medskip

On the other hand, since $\CM$ is partially integrable, the
convolution $j_{\cla\cdot w_0,*}\star \CM$ belongs to
$D^{<0}(\hg_\crit\mod)$, and its map to $j_{\cla,*}\star
\BW_{\crit,w(\rho)-\rho,\reg}\otimes \omega_x^{\langle
\rho-w(\rho),\crho\rangle}\in \hg_\crit\mod$ is necessarily $0$.

\ssec{Proof of \propref{extension of maps}}  \label{proof of ext}

Consider the D-submodule of $j_{\cla,*,X}\star \BW_{\crit,w,X}$
consisting of sections which are annihilated by the Lie-* action of
$\CA_{\fg,\crit,X}$. This is a torsion-free chiral module over
$\fH_{\crit,w,X}$. The restriction of this D-module to $X-x$
identifies canonically with $\CL^\cla_{\fH,X}$.  The claim about the
extension follows now from \cite{FG3}, Proposition 3.4.

\medskip

To prove the second assertion of the proposition, we claim that 
it is enough to show that the resulting map
\begin{equation} \label{map of lines}
\CL^\cla_{\Fl^{\cG}_{w,\fZ}}
\to j_{\cla,*}\star\BW_{\crit,w(\rho)-\rho,\reg}\otimes 
\omega_x^{\langle \rho-w(\rho),\crho\rangle}
\end{equation}
is non-zero.

Indeed, dividing by the line $\omega_x^{\langle w(\rho),\cla\rangle}$,
we obtain that the map in \propref{extension of maps}(2) corresponds
to a non-zero $G[[t]]$-invariant vector in
$\BW_{\crit,w(\rho)-\rho,\reg}\otimes \omega_x^{\langle
\rho-w(\rho),\crho\rangle}$, which has degree $0$ with respect to
$\BG_m\subset \Aut(\D)$. Hence it must coincide, up to a scalar, with
the spherical vector used in the construction of the map \eqref{Vak to
Wak}.

\medskip

Suppose that the map \eqref{map of lines} was zero. We would obtain that
in the commutative diagram \eqref{CD W(2)} the composed map
$$\CF_{V^\cla,X}\star \BV_{\crit,X}\to j_{\cla,*,X}\star \BW_{\crit,w,X}$$
is such that its fiber at $x$ is zero.

In other words, we obtain that the map
$$\CF_{V^\cla}\star \BV_\crit\to j_{\cla,*}\star \BV_\crit\to
j_{\cla,*}\star  \BW_{\crit,w(\rho)-\rho,\reg}\otimes 
\omega_x^{\langle \rho-w(\rho),\crho\rangle},$$
obtained from \eqref{Vak to Wak}, vanishes.

\medskip

Consider the map
\begin{equation} \label{rho & la}
j_{w_0\cdot \crho,*}\star j_{\cla,*}\star \BV_\crit\to
j_{w_0\cdot \crho,*}\star j_{\cla,*}\star 
\BW_{\crit,w(\rho)-\rho,\reg}\otimes 
\omega_x^{\langle \rho-w(\rho),\crho\rangle}.
\end{equation}
We obtain that its composition with
$$\IC_{w_0\cdot (\crho+\cla),\Gr_G}\star \BV_\crit\to j_{w_0\cdot
(\crho+\cla),\Gr_G,*} \star \BV_\crit\simeq j_{w_0\cdot
\crho,*}\star j_{\cla,*} \star \BV_\crit$$ vanishes.

However, by \cite{ABBGM}, Sect. 2.3, the cokernel $j_{w_0\cdot
(\crho+\cla),\Gr_G,*}/\IC_{w_0\cdot (\crho+\cla),\Gr_G}$ is partially
integrable. Hence, the image of the map in \eqref{rho & la} is
partially integrable. But this leads to a contradiction as in the
proof of \propref{Vak to Wak, semiinf} given above.

\section{The general case}    \label{more general}

\ssec{}
In this section we will carry out one more step in the proof 
\thmref{full ident}. Namely, we will prove that an isomorphism
\begin{equation} \label{need to prove general}
_w\BW(\CT)\simeq \sG\Bigl(\CL^{\on{twist}}_w
\underset{\CO_{\Fl^\cG_{w,\on{th},\fZ}}}\otimes \CT\Bigr),
\end{equation}
holds functorially in $\CT\in \QCoh(\Fl^\cG_{w,\on{th},\fZ})$ for
{\it some} bi-module $\CL^{\on{twist}}_w$ over the topological
commutative algebra $\Fun(\Fl^\cG_{w,\on{th},\fZ})$.
Moreover, we will show that as a {\it right} module, $\CL^{\on{twist}}_w$
is a line bundle over $\Fl^\cG_{w,\on{th},\fZ}$.

As an immediate corollary of the above isomorphism, combined with
\secref{sE Z}, we obtain the following:
\begin{cor}
For every $\CT\in \QCoh(\Fl^\cG_{w,\on{th},\fZ})$ the Wakimoto module
$_w\BW(\CT)$ is $I$-equivariant. 
\end{cor}

The main step in the proof of \eqref{need to prove general} will be
the following:

\begin{prop}   \label{isom on Ext}
Let $\CT^1$ and $\CT^2$ be two quasi-coherent sheaves on $\Fl^\cG_{\fZ}$.
Then the functor $\sE^\fZ:D^+\bigl(\QCoh(\Fl^\cG_{\fZ})\bigr)\to
D^+\bigl(\on{D}(\Gr_G)^{\Hecke_\fZ}_\crit\mod^{I^0}\bigr)$ induces
a bijection
$$\Ext^i_{\QCoh(\Fl^{\cG}_{\fZ})}(\CT^1,\CT^2)\to
R^i\Hom_{\on{D}(\Gr_G)^{\Hecke_\fZ}_\crit\mod^{I^0}}\bigl(\sE^\fZ(\CT^1),
\sE^\fZ(\CT^2)\bigr)$$ for $i=0,1$.
\end{prop}

Note that the map in the proposition is {\it not}
bijective for $i\geq 2$, for otherwise the functor $\sE^\fZ$
would be an equivalence of categories, which it is not.

\medskip

Let us show how this proposition implies \eqref{need to prove general}.

\medskip

\noindent{\bf Step 1.}
We claim that $_w\BW(\CT)$ is isomorphic to
$\sG(\CT'):=\Gamma^{\Hecke_\fZ}(\Gr_G,\sE^\fZ(\CT'))$ for
{\it some} object $\CT'\in \QCoh(\Fl^\cG_{w,\on{th},\fZ})$.

\medskip

Without loss of generality we can assume that $\CT$ is
supported on $k$-th infinitesimal neighborhood of $\Fl^{\cG}_{w,\fZ}$
in $\Fl^{\cG}_\fZ$. We will argue by induction on $k$.

If $k=0$, i.e., when $\CT$ is scheme-theoretically supported on
$\Fl^{\cG}_{w,\fZ}$, we have:
$$_w\BW(\CT)\simeq
\BW_{\crit,w(\rho)-\rho,\reg}\underset{\Fun(\Fl^{\cG}_{w,\fZ})}\otimes
\CT,$$ and the assertion follows from \thmref{ident w Wak}.

Suppose now that $k>1$. Then we can write $\CT$ as an extension
$$0\to \CT^1\to \CT\to \CT^2\to 0,$$ with $\CT^1$ supported on
on $\Fl^{\cG}_{w,\fZ}$ and $\CT^2$ supported on the $k-1$-st
infinitesimal neighborhood of $\Fl^{\cG}_{w,\fZ}$ in
$\Fl^{\cG}_\fZ$. Then, by the induction hypothesis, there exist
objects $\CT'{}^1,\CT'{}^2\in \QCoh(\Fl^\cG_{w,\on{th},\fZ})$, such that
$$\sG(\CT'{}^i)\simeq {}_w\BW(\CT^i)$$
for $i=1,2$. 

By the main theorem of \cite{FG4}, the extension
$$0\to {}_w\BW(\CT^1)\to {}_w\BW(\CT)\to {}_w\BW(\CT^2)\to 0$$ 
comes from some extension
$$0\to \sE^\fZ(\CT'{}^1)\to \CF\to \sE^\fZ(\CT'{}^2)\to 0$$ in
$\on{D}(\Gr_G)^{\Hecke_\fZ}_\crit\mod^{I^0}$. However, by
\propref{isom on Ext}
$$\CF\simeq \sE^\fZ(\CT')$$
for some $\CT'\in \QCoh(\Fl^\cG_{\fZ})$, as required. Evidently,
$\CT'$, being an extension of two quasi-coherent sheaves 
that belong to $\QCoh(\Fl^\cG_{w,\on{th},\fZ})$, itself belongs
to this subcategory.

\medskip

\noindent{\bf Step 2.}  \propref{isom on Ext} for $i=0$, combined with
the equivalence of \cite{FG4}, implies that the assignment $\CT\mapsto
\CT'$, constructed above, is a functor
$\QCoh(\Fl^\cG_{w,\on{th},\fZ})\to \QCoh(\Fl^\cG_{w,\on{th},\fZ})$;
let us denote it by $_w\sQ$. Since each of the functors $_w\BW$,
$\Gamma^{\Hecke_\fZ}$ and $\sE|_{\QCoh(\Fl^\cG_{w,\on{th},\fZ})}$ is
exact and faithful, we obtain that $_w\sQ$ is also exact and faithful.

Hence, $_w\sQ$ has the form
$$\CT\mapsto \CL^{\on{twist}}_w
\underset{\Fun(\Fl^\cG_{w,\on{th},\fZ})}\otimes\CT,$$ for a certain
$\Fun(\Fl^\cG_{w,\on{th},\fZ})$-bimodule $\CL^{\on{twist}}_w$. We will
show now that as a {\it right} $\Fun(\Fl^\cG_{w,\on{th},\fZ})$-module,
$\CL^{\on{twist}}_w$ is a line bundle. In fact, we will show that it
is non-canonically trivial.

This is equivalent to showing that there exists a functorial isomorphism
between the vector space underlying $_w\sQ(\CT)$ and that of $\CT$.
We will do this by comparing the semi-infinite cohomologies 
$H^\semiinf(\fn\ppart,\fn[[t]],?\otimes \Psi_0)$.

\medskip

\noindent{\bf Step 3.} 

Recall that by \cite{FG2}, Sect. 12.4, there exists an isomorphism
\begin{equation} \label{naive semiinf of Wak}
H^\semiinf(\fn\ppart,\fn[[t]],\BW(\CT)\otimes \Psi_0)\simeq \CT,
\end{equation}
which is functorial in $\CT\in \wh\fH_\crit\mod$, but it is
non-canonical in the sense that it depends on the choice of the
coordinate $t$ on $\D$.

\medskip

Thus, to identify $_w\sQ(\CT)$ and $\CT$ as vector spaces, it suffices
to prove the following:

\begin{prop}  \label{semiinf of convolution}
For any $\CT'\in \QCoh(\Fl^{\cG}_\fZ)$, there exists a canonical
quasi-isomorphism
$$H^\semiinfb\Bigl(\fn\ppart,\fn[[t]],\sG(\CT')\otimes
\Psi_0\Bigr)\simeq R\Gamma(\Fl^\cG_\fZ,\CT').$$
\end{prop}

\begin{proof}

By the definition of convolution, the LHS of the proposition is given
by applying the functor of derived $\cB_\fZ$-invariants to the complex
\begin{equation} \label{semiinf of pull-back}
H^\semiinfb\biggl(\fn\ppart,\fn[[t]],
\Gamma^{\Hecke_\fZ}\Bigl(\Gr_G,\act^*(\CW^\fZ_{w_0})
\underset{\Fun(\cG_\fZ)}\otimes \wt\CT'\Bigr)\otimes \Psi_0\biggr),
\end{equation}
where $\wt\CT'$ denotes the pull-back of $\CT'$ to $\cG_\fZ$.

However, by \cite{FG4}, Sect. 2.12, 
\begin{align*}
&H^\semiinfb\biggl(\fn\ppart,\fn[[t]],
\Gamma^{\Hecke_\fZ}\Bigl(\Gr_G,\act^*(\CW_{w^\fZ_0})\Bigr) \otimes
\Psi_0\biggr)\simeq \\ &\simeq
H^\semiinfb\Bigl(\fn\ppart,\fn[[t]],\Gamma^{\Hecke_\fZ}(\Gr_G,\CW^\fZ_{w_0})
\otimes \Psi_0\Bigr)\underset{\fZ^\reg_\fg}\otimes \Fun(\cG_\fZ)\simeq
\Fun(\cG_\fZ).
\end{align*}

Hence, the expression in \eqref{semiinf of pull-back} is isomorphic to
$\wt\CT'$, as a $\Fun(\cG_\fZ)$-module, endowed with a
$\cB_\fZ$-action. Finally, we have:
$$R\on{Inv}(\cB_\fZ,\wt\CT')\simeq R\Gamma(\Fl^{\cG}_\fZ,\CT'),$$
which is what we had to show.

\end{proof}

\ssec{Proof of \propref{isom on Ext}}   \hfill

\medskip

\noindent{\bf Step 1.}  We claim that it is enough to prove a version
of the proposition for the categories $(\Fl^\cG_\fZ,
\on{D}(\Gr_G)^{\Hecke_\fZ}_\crit\mod^{I^0})$ replaced by $(\Fl^\cG,
\on{D}(\Gr_G)^{\Hecke}_\crit\mod^{I^0})$, i.e., that for
$\CT^1,\CT^2\in \QCoh(\Fl^\cG)$ the map
\begin{equation} \label{isom ext simple}
\Ext^i_{\QCoh(\Fl^{\cG})}(\CT^1,\CT^2)\to
R^i\Hom_{\on{D}(\Gr_G)^{\Hecke}_\crit\mod^{I^0}}\bigl(\sE(\CT^1),
\sE(\CT^2)\bigr)
\end{equation}
is an isomorphism for $i=0,1$.

Indeed, by choosing a trivialization of the $\cG$-torsor
$\CP_{\cG,\fZ}$, we can identify $\Fl^\cG_\fZ$ (resp.,
$\on{D}(\Gr_G)^{\Hecke_\fZ}_\crit\mod^{I^0}$) with the category of
objects of $\QCoh(\Fl^\cG)$ (resp.,
$\on{D}(\Gr_G)^{\Hecke}_\crit\mod^{I^0}$), endowed with an action of
the algebra $\fZ^\reg_\fg$, and the functor $\sE^\fZ$ is obtained from
$\sE$ by extension of scalars.

The fact that \eqref{isom ext simple} is an isomorphism for $i=0,1$ is
a formal corollary of \cite{ABG}. Below we will give an independent
proof.

\medskip

\noindent{\bf Step 2.} 
We claim that it is enough to show that the maps
\begin{equation} \label{map on Ext la}
\Ext^i_{\QCoh(\Fl^{\cG})}\bigl(\CL^{\cla}_{\Fl^{\cG}},\CT\bigr)
\to R^i\Hom_{\on{D}(\Gr_G)^{\Hecke}_\crit\mod^{I^0}}
\bigl(\sE(\CL^{\cla}_{\Fl^{\cG}}),\sE(\CT)\bigr)
\end{equation}
are isomorphisms for $i=0,1$, $\cla\in \cLambda^+$, 
and any $\CT\in \QCoh(\Fl^{\cG})$.

\medskip

Without loss of generality we can assume that $\CT^1$ 
is coherent, and let $$...\to\CQ^2\to \CQ^1\to \CQ^0\to \CT^1\to 0,$$
be a resolution
which each $\CQ^i$ is isomorphic to a direct sum of line bundles
$\CL^{\cla^j_i}_{\Fl^{\cG}}$.

The cohomological dimension of the functor $\sE$ is {\it a priori}
bounded by the cohomological dimension of the category
$\QCoh(\Fl^\cG)$, which is $\dim(\Fl^\cG)$. Let $\cla\in \cLambda^+$
be such that $\cla^j_i+\cla\in \cLambda^+$ for $i\leq
\dim(\Fl^\cG)+1$. Then by \propref{verify full A},
$\sE(\CL^{\cla^j_i+\cla}_{\Fl^{\cG}})\simeq
\on{Ind}^\Hecke(j_{\cla^j_i+\cla,\Gr_G,*})$.  In particular, this
implies that $\sE(\CL^\cla_{\Fl^\cG}\otimes\CT)\in \on{D}(\Gr_G)^{\Hecke}_\crit\mod^{I^0}$,
i.e., $\sE(\CL^\cla_{\Fl^\cG}\otimes\CT)$ does not have higher cohomologies.

\medskip

By \propref{coh and I conv}, for $\CT^1,\CT^2$ as above we have a
commutative diagram
$$
\CD
R\Hom(\CT^1,\CT^2)  @>{\sE}>> 
R\Hom\bigl(\sE(\CT^1),\sE(\CT^2)\bigr) \\
@V{\sim}VV  @VVV  \\
R\Hom(\CL^\cla_{\Fl^\cG}\otimes \CT^1,\CL^\cla_{\Fl^\cG}\otimes \CT^2)
@>{\sE}>>
R\Hom\bigl(j_{\cla,*}\star \sE(\CT^1), j_{\cla,*}\star
\sE(\CT^2)\bigr).  \endCD
$$

Now notice that the functor $\wt{j}_{\cla,*}\underset{I^0}\star ?$ is
a self-equivalence of
$D\bigl(\on{D}(\Gr_G)^{\Hecke}_\crit\mod^{I^0}\bigr)$.  (Its
quasi-inverse is given by the $!$-convolution with
$\wt{j}_{-\cla,*}$.)  Hence, the right vertical arrow in the above
diagram is an isomorphism as well.

\medskip

Hence, we can replace the initial sheaves $\CT^1$ and $\CT^2$ by their
twists with respect to $\CL^\cla_{\Fl^\cG}$. Moreover, we can use the
3-term resolution
$$\CL^\cla_{\Fl^\cG}\otimes \CQ_2\to \CL^\cla_{\Fl^\cG}\otimes \CQ_1\to
\CL^\cla_{\Fl^\cG}\otimes \CQ_0\to \CL^\cla_{\Fl^\cG}\otimes \CT^1$$
to compute both sides of \eqref{isom ext simple}. This performs the
required reduction in Step 2. 

Using \propref{coh and I conv} again, we reduce the assertion further to 
the case $\cla=0$.

\medskip

\noindent{\bf Step 3.}  Let us recall the general set-up of
\secref{coherent convolutions}.  Let $\CM$ (resp., $\CN$) be a
$\cG$-equivariant (resp., $\cB$-equivariant) object of $\CC$, and let
$\CT$ be an object of $\QCoh(\Fl^\cG)$.  In this case
$R\Hom_\CC(\CM,\CN)$ is naturally an object of the derived category of
$\cB$-modules. Let
$\wt{R\Hom}(\CM,\CN)$ denote the associated complex of
$\cG$-equivariant quasi-coherent sheaves on $\Fl^\cG$.

For $\CT\in \QCoh(\Fl^\cG)$ we have:
\begin{equation} \label{convert to sections}
R\Hom_{\CC}\bigl(\CM,\CT\overset{R}\ast \CN\bigr)\simeq
R\Gamma\Bigl(\Fl^\cG,\CT\underset{\CO_{\Fl^\cG}}\otimes 
\wt{R\Hom}(\CM,\CN)\Bigr).
\end{equation}

Applying this to $\{\CC,\CM,\CN\})$ being
$$\{\QCoh(\Fl^\cG),\CO_{\Fl^\cG},\BC_{w_0}\} \text{ and }
\{\on{D}(\Gr_G)^\Hecke_\crit\mod^{I^0},
\sE(\CO_{\Fl^\cG}),\sE(\BC_{w_0})\},$$ we obtain a commutative diagram
$$
\CD
R\Hom(\CO_{\Fl^\cG},\CT)   @>{\text{\eqref{isom ext simple}}}>> 
R\Hom_{\on{D}(\Gr_G)^{\Hecke}_\crit\mod^{I^0}}
\bigl(\sE(\CO_{\Fl^{\cG}}),\sE(\CT)\bigr)  \\
@V{\sim}VV     @V{\sim}VV   \\
R\Gamma \Bigl(\Fl^\cG,\CT\underset{\CO_{\Fl^\cG}}\otimes 
\wt{R\Hom}(\CO_{\Fl^\cG},\BC_{w_0})\Bigr)
@>>> 
R\Gamma \Bigl(\Fl^\cG,\CT\underset{\CO_{\Fl^\cG}}\otimes 
\wt{R\Hom}(\sE(\CO_{\Fl^\cG}),\sE(\BC_{w_0}))\Bigr) \endCD
$$

Thus, we obtain that it is enough to show that \eqref{map on Ext la}
is an isomorphism for $i=0,1$, $\cla=0$ and $\CT=\BC_{w_0}$.

\medskip

Let $\Fl^\cG_1\subset \Fl^\cG$ be the open Schubert cell, and let us
apply the above commutative diagram again with
$\CT=\CO_{\Fl^\cG_1}$. Since $\Fl^\cG_1$ is affine, we obtain that the
isomorphism of \eqref{map on Ext la} for $\CT=\BC_{w_0}$ will follow
once we establish it for $\CT=\CO_{\Fl^\cG_1}$.

\medskip

\noindent{\bf Step 4.}
Let us denote $\sE(\CO_{\Fl^\cG_1})$ by $\CW_1$, and recall that
$\sE(\CO_{\Fl^\cG})\simeq \on{Ind}^\Hecke(\delta_{1,\Gr_G})$.
It remains to show that 
\begin{equation}  \label{last hom}
\Fun(\Fl^\cG_1)\to 
\Hom_{\on{D}(\Gr_G)_\crit\mod^{I^0}}\Bigl(\delta_{1,\Gr_G},\CW_1\Bigr)
\end{equation}
is an isomorphism and
\begin{equation} \label{last ext}
\Ext^1_{\on{D}(\Gr_G)_\crit\mod^{I^0}}
\Bigl(\delta_{1,\Gr_G},\CW_1\Bigr)=0,
\end{equation}
where in both formulas $\CW_1$ appears as an object
of $\on{D}(\Gr_G)_\crit\mod^{I^0}$ via the tautological forgetful functor.

\medskip

Recall (see \secref{another way}) that 
$\CO_{\Fl^\cG_1}$ can be written as a filtered direct limit
$$\underset{\cmu\in \cLambda^+}{\underset{\longrightarrow}{\lim}}\,
\CL^\cmu_{\Fl^\cG}\otimes \fl^{-\cmu},$$
where $\fl^\cmu$ denotes the $\cB$-stable line in $V^\cmu$, and
$\fl^{-\cmu}\subset (V^\cmu)^*$ is its dual.

Hence, $\CW_1$ can also be written down as a filtered direct limit
$$\underset{\cmu\in \cLambda^+}{\underset{\longrightarrow}{\lim}}\,
\on{Ind}^\Hecke(j_{\cmu,\Gr_G,*})\otimes \fl^{-\cmu}.$$

Recall that $\on{Ind}^\Hecke(\CF)\simeq \CF\underset{G[[t]]}\star
\CF_{R(\cG)}$.  We have $\CF_{R(\cG)}\simeq \underset{\cnu}\oplus\,
\CF_{(V^\cnu)^*}\otimes V^\cnu$.  Now the isomorphism in \eqref{last
hom} follows from the fact that
$$
\begin{cases}
& \Hom(\delta_{1,\Gr_G}, j_{\cmu,\Gr_G,*}\underset{G[[t]]}\star
\CF_{(V^\cnu)^*})=0,\,\, \cmu\neq \cnu, \\
& \Hom(\delta_{1,\Gr_G}, j_{\cmu,\Gr_G,*}\underset{G[[t]]}\star
\CF_{(V^\cnu)^*}\simeq \BC,\,\, \cmu=\cnu.
\end{cases}
$$

The vanishing of \eqref{last ext} follows  from the fact that 
$$\Ext^1(\delta_{1,\Gr_G},j_{\cmu,\Gr_G,*}\underset{G[[t]]}\star
\CF_{(V^\cnu)^*})\simeq \Ext^1(\IC_{\Grb^\cnu},j_{\cmu,\Gr_G,*})=0,$$
by the parity vanishing of IC-stalks on $\Gr_G$.

\section{Lie algebroids and the renormalized enveloping algebra}
\label{action}

In this section we will study the interaction between
the isomorphism $\on{map}_{\geom}$ and certain canonical
Lie algebroids defined on the two sides of this isomorphism.

\ssec{}  \label{algbrds}

Let $N^*_{\fZ^\reg_\fg/\fZ_\fg}$ denote the conormal to
$\Spec(\fZ^\reg_\fg)$ inside $\Spec(\fZ_\fg)$, equipped with a natural
topology. Recall (see \cite{BD}, Sect. 3.6 or \cite{FG2}, Sect. 7.4)
that $N^*_{\fZ^\reg_\fg/\fZ_\fg}$ has a natural structure of Lie
algebroid over $\Spec(\fZ_\fg)$.

On the other hand, let us recall the groupoid
$\Isom_{\cG,\fZ^\reg_\fg}$ over $\Spec(\fZ^\reg_\fg)$, see
\eqref{consider groupoid}, and let $\isom_{\cG,\fZ^\reg_\fg}$ be the
corresponding Lie algebroid.

According to \cite{BD}, Theorem 3.6.7, there exists a canonical
morphism (in fact, an isomorphism) of Lie algebroids:
\begin{equation} \label{algebroid map}
\upsilon_{\geom}:N^*_{\fZ^\reg_\fg/\fZ_\fg}\to \isom_{\cG,\fZ_\fg^\reg}.
\end{equation}

Below we will recall the definition of this map. Let us note that both
the Lie algebroid structure on $N^*_{\fZ^\reg_\fg/\fZ_\fg}$ and the
morphism $\upsilon_{\geom}$ depend on an additional choice of a
one-parameter deformation $\kappa_\hslash$ of the level away from the
critical value.

\ssec{The renormalized enveloping algebra}   \label{rnrm}

Let $U^{\ren,\reg}(\hg_\crit)$ be the renormalized universal
enveloping algebra at the critical level, which is defined in
\cite{BD}, Sect. 5.6 (see also \cite{FG2}, Sect. 7.4 for a review). It
admits a natural filtration, with the $0$-th term
$U^{\ren,\reg}(\hg_\crit)_0$ isomorphic to the topological algebra
$$U^\reg(\hg_\crit):=\wt{U}_\crit(\hg)\underset{\fZ_\fg}\otimes
\fZ^\reg_\fg,$$ responsible for the category $\hg_\crit\mod_\reg$ with
its tautological forgetful functor to $\Vect$. The first associated
graded quotient
$U^{\ren,\reg}(\hg_\crit)_1/U^{\ren,\reg}(\hg_\crit)_0$ is isomorphic
to
$$U^\reg(\hg_\crit)\underset{\fZ^\reg_\fg}\hattimes
N^*_{\fZ^\reg_\fg/\fZ_\fg}.$$ 
Note that $\fZ^\reg_\fg$ is a subalgebra in $U^{\ren,\reg}(\hg_\crit)$,
but it is no longer central.

\medskip

Let us recall some basic constructions related to
$U^{\ren,\reg}(\hg_\crit)$.

\medskip

\noindent{\bf 1.} Let $\CM_1,\CM_2$ be two objects of
$\hg_\crit\mod_\reg$, on which the action of $U^\reg(\hg_\crit)$ has
been extended to an action of $U^{\ren,\reg}(\hg_\crit)$. Then
$\Hom_{\hg_\crit}(\CM_1,\CM_2)$ acquires a natural action of
$N^*_{\fZ^\reg_\fg/\fZ_\fg}$ via
$N^*_{\fZ^\reg_\fg/\fZ_\fg}\hookrightarrow
U^{\ren,\reg}(\hg_\crit)_1/U^{\ren,\reg}(\hg_\crit)_0$.

\medskip

\noindent{\bf 2.} Let $\CM$ be a $U^{\ren,\reg}(\hg_\crit)$-module,
and let $\CL$ be a $\fZ^\reg_\fg$ module that carries a compatible
action of $N^*_{\fZ^\reg_\fg/\fZ_\fg}$. Then
$\CM\underset{\fZ^\reg_\fg}\otimes \CL$ is naturally a
$U^{\ren,\reg}(\hg_\crit)$-module.

Indeed, we define the action of
$$U^{\ren,\reg}(\hg_\crit)_1\underset{U^\reg(\hg_\crit)
\underset{\fZ^\reg_\fg}\hattimes N^*_{\fZ^\reg_\fg/\fZ_\fg}}\times
N^*_{\fZ^\reg_\fg/\fZ_\fg}\subset U^{\ren,\reg}(\hg_\crit)$$ to be the
sum of the given $U^{\ren,\reg}(\hg_\crit)_1$-action on $\CM$ and the
$N^*_{\fZ^\reg_\fg/\fZ_\fg}$-action on $\CL$. From the relations that
realize $U^{\ren,\reg}(\hg_\crit)$ as a quotient of the universal
enveloping algebra of $U^{\ren,\reg}(\hg_\crit)_1$, it follows that
the above definition extends to a well-defined
$U^{\ren,\reg}(\hg_\crit)$-action on
$\CM\underset{\fZ^\reg_\fg}\otimes \CL$.

\medskip

\noindent{\bf 3.} Let $\CM_\hslash$ be a flat
$\BC[\hslash]$-family of $\hg_{\kappa_\hslash}$-modules, such that
$\CM_\crit:=\CM_\hslash/\hslash\cdot \CM$ belongs to $\hg_\crit\mod_\reg$, 
then $\CM_\crit$ acquires an action of $U^{\ren,\reg}(\hg_\crit)$.  

The prime example of this situation is $\CM_\hslash=\BV_{\hslash}$.
In this case we obtain the canonical $U^{\ren,\reg}(\hg_\crit)$-action
on $\BV_\crit$.

\medskip

\noindent{\bf 4.} The adjoint action of $G\ppart$ on
$U^\reg(\hg_\crit)$ extends naturally to an action on
$U^{\ren,\reg}(\hg_\crit)$. Hence, the category
$U^{\ren,\reg}(\hg_\crit)\mod$ is acted on by critically twisted
D-modules on $G\ppart$ by convolutions (see \cite{FG2}, Sect. 22).

In particular, if $\CM$ is a $U^{\ren,\reg}(\hg_\crit)$-module, which
is $\fg[[t]]$-integrable as a $\hg_\crit$-module, and $\CF\in
\on{D}(\Gr_G)_\crit\mod$, we obtain a well-defined object
$\CF\underset{G[[t]]}\star \CM$ in the derived category of
$U^{\ren,\reg}(\hg_\crit)$-modules. Applying this to $\CM=\BV_\crit$,
we recover the canonical $U^{\ren,\reg}(\hg_\crit)$-action on
$\CF\underset{G[[t]]}\star \BV_\crit\simeq \Gamma(\Gr_G,\CF)$.


 

\ssec{}

We are now ready to construct the map $\upsilon_{\geom}$.  By
definition, $\isom_{\cG,\fZ_\fg^\reg}$ is the Atiyah algebroid of the
$\cG$-torsor $\CP_{\cG,\fZ}$. Hence, to specify $\upsilon_{\geom}$ we
need to make $N^*_{\fZ^\reg_\fg/\fZ_\fg}$ act on $\CV_\fZ$ for every
$V\in \Rep(\cG)$ in a way compatible with tensor products.

However, by construction, $\CV_\fZ\simeq
\Hom_{\hg_\crit}(\BV_\crit,\Gamma(\Gr_G,\CF_V))$, and the required
action follows from {\bf 1} and {\bf 4} above. In other words, the
isomorphism $\Gamma(\Gr_G,\CF_V)\simeq
\BV_\crit\underset{\fZ^\reg_\fg}\otimes \CV_\fZ$ is compatible with
the $U^{\ren,\reg}(\hg_\crit)$-actions, where on the RHS it is given
via the construction of {\bf 2} above. The compatibility with tensor
products follows from the compatibility between the constructions in
{\bf 2} and {\bf 4}.

\medskip

We will now recall the definition of another map
\begin{equation} \label{algebroid map, FF}
\upsilon_{\alg}:N^*_{\fZ^\reg_\fg/\fZ_\fg}\to
\isom_{\cG,\fZ_\fg^\reg},
\end{equation}
following \cite{BD}, Sect. 3.7.13.

We recall that the isomorphism of topological algebras
\begin{equation} \label{full FF}
\on{map}_{\alg}:\fZ_\fg\simeq \Fun(\Op_\cg(\D^\times))
\end{equation}
respects the Poisson structures, where on the LHS the Poisson bracket
is defined via the Drinfeld-Sokolov reduction (see \cite{FG2}, Sect. 4
for a review), using the $\cG$-invariant form $\check\kappa$ on $\cg$,
corresponds to the deformation $\kappa_\hslash$ of the level off the
critical value (see \secref{review Miura} for the precise
formulation).

Thus, the map $\on{map}_{\alg}$ induces an isomorphism of algebroids
$$N^*_{\fZ^\reg_\fg/\fZ_\fg}\simeq
\on{map}^*_{\alg}\bigl(N^*_{\Op^\reg_\cg/\Op_\cg(\D^\times)}\bigr).$$
Combining this with the isomorphism
\begin{equation} \label{isom alg opers}
N^*_{\Op^\reg_\cg/\Op_\cg(\D^\times)}\simeq \isom_{\cG,\Op^\reg_\cg}
\end{equation}
(see, e.g., \cite{FG2}, Sect. 4.4, where the latter is explained), we
obtain the map of \eqref{algebroid map, FF}. By construction, the map
$\upsilon_{\alg}$ is an isomorphism.

\medskip

Our present goal of this section is to prove the following:

\begin{thm}  \label{two maps of algebroids}
The maps $\upsilon_{\geom}$ and $\upsilon_{\alg}$ coincide.
\end{thm}

The above theorem has been proved in \cite{BD},
Proposition 3.5.13, simultaneously with \thmref{two maps
coincide}. Namely, in {\it loc. cit.} it was shown that the algebroid
$\isom_{\cG,\fZ_\fg^\reg}$ over $\Spec(\fZ^\reg_\fg)$ does not admit
non-trivial automorphisms. We will give a constructive proof of
this result, which will occupy the rest of this section.

\ssec{}   \label{sect two actions on H}

The algebroid $\isom_{\cG,\fZ_\fg^\reg}$ acts
naturally on the scheme $\Fl^{\cG}_{1,\fZ}$, and it is easy to see
that the corresponding action on $\Fun(\Fl^{\cG}_{1,\fZ})$ is
faithful. Hence, to prove \thmref{two maps of algebroids}, it suffices
to see that the two resulting actions of $N^*_{\fZ^\reg_\fg/\fZ_\fg}$
on $\Fun(\Fl^{\cG}_{1,\fZ})$-- one via $\upsilon_{\geom}$ and another
via $\upsilon_{\alg}$-- coincide.

\medskip

Consider the $U^{\ren,\reg}(\hg_\crit)$-action on $\BW_{\crit,0}$,
corresponding to the $\hslash$-family $\BW_{\hslash,0}$.  Recall
that we have an isomorphism
$$\fH^\reg_\crit\simeq \End_{\hg_\crit}(\BW_{\crit,0}).$$
Hence, from \secref{rnrm}({\bf 1}) we obtain an action
of $N^*_{\fZ^\reg_\fg/\fZ_\fg}$ on $\fH^\reg_\crit$.

Recall now isomorphism
$$\fH^\reg_\crit\simeq \Fun(\Fl^\cG_{1,\fZ}).$$

Thus, to prove \thmref{two maps of algebroids} it suffices to prove
the following two assertions:

\begin{prop} \label{two actions on H}
The above action of $N^*_{\fZ^\reg_\fg/\fZ_\fg}$ on $\fH^\reg_\crit$
goes over under the map $\upsilon_{\geom}$ to the natural action of
$\isom_{\cG,\fZ_\fg^\reg}$ on $\Fun(\Fl^\cG_{1,\fZ})$.
\end{prop}

\begin{prop} \label{char ups alg}
The above action of $N^*_{\fZ^\reg_\fg/\fZ_\fg}$ on $\fH^\reg_\crit$
goes over under the map $\upsilon_{\alg}$ to the natural action of
$\isom_{\cG,\fZ_\fg^\reg}$ on $\Fun(\Fl^\cG_{1,\fZ})$.
\end{prop}

\ssec{Proof of \propref{two actions on H}} \label{proof two actions on
H}

Let $\Spec(\CB)$ be an affine scheme over $\Spec(\fZ^\reg_\fg)$,
endowed with a map to $\Fl^\cG_\fZ$. Note that such a data is
specified by a $\cH$-torsor $\{\cla\mapsto \CL^\cla_\CB\}$ and a
collection of maps
$$\kappa_\CB^{-,\cla}:\CV^\cla_\fZ\to \CL^\cla_\CB,\,\,\cla\in
\cLambda^+,$$ satisfying the Pl\"ucker equations.

Let $\ff$ be a Lie algebroid over $\Spec(\fZ^\reg_\fg)$, endowed with a map
$\upsilon:\ff\to \isom_{\cG,\fZ_\fg^\reg}$. In particular, for every
$V\in \Rep(\cG)$ we obtain an $\ff$-action on $\CV_\fZ$.

\medskip

Suppose, in addition, that we are given an action of $\ff$ on $\CB$.
Then these data are compatible with the natural
$\isom_{\cG,\fZ_\fg^\reg}$-action on $\Fl^\cG_\fZ$ if and only if the
following holds:

\noindent For every $\cla\in \cLambda^+$ there exists an $\ff$-action
on the line bundle $\CL^\cla_\CB$, such that the map
$\kappa_\CB^{-,\cla}$ is compatible with the $\ff$-actions.  (Note
that such an action on $\CL^\cla_\CB$ is {\it a priori} unique, since
the induced maps $\CV^\cla_\fZ\underset{\fZ^\reg_\fg}\otimes \CB\to
\CL^\cla_\CB$ are surjective.)

\medskip

Let us apply this for $\CB=\fH^\reg_\crit$ and
$\ff=N^*_{\fZ^\reg_\fg/\fZ_\fg}$.  Recall that the map
$\Spec(\fH^\reg_\crit)\to \Fl^\cG_\fZ$ corresponds to the collection
of line bundles $\CL^\cla_\fH$ and the maps $\kappa^{-,\cla}_\fH$ of
\eqref{des Miura, local}.

\medskip

Consider the Wakimoto module $\BW_{\crit,0}$, endowed with the natural
$U^{\ren,\reg}(\hg_\crit)$-action. By construction, the map
$\phi:\BV_\crit\to \BW_{\crit,0}$ deforms off the critical level to a
map $\phi_\hslash:\BV_{\hslash}\to \BW_{\hslash,0}$; hence the map
$\phi$ is compatible with the $U^{\ren,\reg}(\hg_\crit)$-actions.

Consider now $j_{\cla,*}\star \BW_{\crit,0}$. From \secref{rnrm}({\bf
4}), we obtain that $j_{\cla,*}\star \BW_{\crit,0}$ also carries an
action of $U^{\ren,\reg}(\hg_\crit)$, and the map
$$\Gamma(\Gr_G,\CF_{V^\cla})\to j_{\cla,*}\star \BW_{\crit,0}$$
of \eqref{imp comp} is compatible with the
$U^{\ren,\reg}(\hg_\crit)$-actions.

{}From \secref{rnrm}({\bf 1}) We obtain that $\CL^\cla_\fH\simeq
\Hom_{\hg_\crit}(\BW_{\crit,0}, j_{\cla,*}\star \BW_{\crit,0})$
acquires an action of $N^*_{\fZ^\reg_\fg/\fZ_\fg}$, and the map
$\kappa^{-,\cla}_\fH$, which equals
\begin{multline*}
\CV^\cla_\fZ\simeq \Hom_{\hg_\crit}\Bigl(\BV_\crit,
\Gamma(\Gr_G,\CF_{V^\cla})\Bigr)\to \\
\to\Hom_{\hg_\crit}\Bigl(\BV_\crit,j_{\cla,*}\star \BW_{\crit,0}\Bigr)
\overset{\sim}\leftarrow
\Hom_{\hg_\crit}\Bigl(\BW_{\crit,0},j_{\cla,*}\star \BW_{\crit,0}\Bigr)
\simeq \CL^\cla_\fH,
\end{multline*}
is compatible with the $N^*_{\fZ^\reg_\fg/\fZ_\fg}$-actions, which is
what we had to show.

\ssec{The Poisson structure on Miura opers}   \label{review Miura}

In order to prove \propref{char ups alg} we need to make a digression
and discuss the Poisson structure on the space of Miura opers. 

\medskip

Let $\frac{d\kappa_\hslash}{d\hslash}|_{\kappa_\crit}$ be the
derivative of $\kappa_\hslash$ at $\kappa_\crit$, which is, by
definition, a non-degenerate $G$-invariant bilinear form on $\fg$. We
can view it as a non-degenerate $W$-invariant bilinear form on $\fh$,
which, in turn, can be thought of as a non-degenerate $W$-invariant
bilinear form on $\check\fh$, or as a $\cG$-invariant bilinear form
on $\cg$. We will denote the latter by $\check\kappa$.

\medskip

Recall the scheme $\ConnDt$. Proceeding as in \cite{FG2}, Sect. 4.3,
(with $\cG$ replaced by $\cH$), using the form $\check\kappa$, we
define a Poisson structure on $\ConnDt$.

Let $\Isom_{\cH,\ConnDt}$ denote the natural groupoid acting on
$\ConnDt$, whose fiber over two points $\wh\chi_1,\wt\chi_2\in
\ConnDt$ is the ind-scheme of automorphisms of the $\cH$-torsor
$\omega^\crho_{\D^\times}$ over $\D^\times$, that transforms 
connection $\wh\chi_1$ to $\wh\chi_2$. Let
$\isom_{\cH,\ConnDt}$ denote the corresponding Lie algebroid on
$\ConnDt$.

As in \cite{FG2}, Sect. 4.3, one easily shows that there exists a 
canonical isomorphism of Lie algebroids
\begin{equation} \label{Miura algebroids}
\Omega^1(\ConnDt)\simeq \isom_{\cH,\ConnDt}.
\end{equation}

\medskip

Recall that the Poisson structure on $\Op_\cg(\D^\times)$ also
depended on the form $\check\kappa$ on $\cg$. It is
a straightforward calculation that the map 
$\on{MT}: \ConnDt\to \Op_\cg(\D^\times)$ is compatible with
the Poisson structures. Therefore, the pull-back
$\on{MT}^*\bigl(\Omega^1(\Op_\cg(\D^\times))\bigr)$ acquires
a structure of Lie algebroid and we have a homomorphism
of Lie algebroids over $\ConnDt$:
$$\on{MT}^*\bigl(\Omega^1(\Op_\cg(\D^\times))\bigr)\to
\Omega^1(\ConnDt).$$

\medskip

Recall now the Lie algebroid $\isom_{\cG,\Op_\cg(\D^\times)}$ over
$\Op_\cg(D^\times)$.  By the definition of generic Miura opers, this
algebroid acts on the ind-scheme $\MOp_{\cg,\gen}(\D^\times)$, 
which, as we know, identifies with $\ConnDt$.

Hence, the pull-back
$\on{MT}^*\bigl(\isom_{\cG,\Op_\cg(\D^\times)}\bigr)$ acquires the
structure of a Lie algebroid on the ind-scheme $\ConnDt$.  In fact,
the anchor map $\on{MT}^*\bigl(\isom_{\cG,\Op_\cg(\D^\times)}\bigr)\to
T(\ConnDt)$ naturally factors through a map
$$\on{MT}^*\bigl(\isom_{\cG,\Op_\cg(\D^\times)}\bigr)\to
\isom_{\cH,\ConnDt}.$$

Finally, recall that we have a canonical isomorphism of algebroids on
$\Op_\cg(\D^\times)$ (see \cite{FG2}, Sect. 4.3):
\begin{equation} \label{DS}
\Omega^1(\Op_\cg(\D^\times))\simeq \isom_{\cG,\Op_\cg(\D^\times)}.
\end{equation}

The following results from the constructions:
\begin{lem}  \label{mkdv flows}
The following diagram is commutative
$$ 
\CD 
\on{MT}^*\bigl(\Omega^1(\Op_\cg(\D^\times))\bigr)
@>{\on{MT}^*\text{\eqref{DS}}}>>
\on{MT}^*\bigl(\isom_{\cG,\Op_\cg(\D^\times)}\bigr) \\ @VVV @VVV \\
\Omega^1(\ConnDt) @>{\text{\eqref{Miura algebroids}}}>>
\isom_{\cH,\ConnDt}.  \endCD
 $$
 \end{lem}

As the result, we obtain:

\begin{cor} \label{mkdv reg}
The action of $N^*_{\Op^\reg_\cg/\Op_\cg(\D^\times)}$ on $\ConnD$,
resulting from the Poisson structure on on $\ConnDt$ and the map
$\on{MT}:\ConnD\to \Op^\reg_\cg$, identifies via
$$N^*_{\Op^\reg_\cg/\Op_\cg(\D^\times)}
\overset{\text{\eqref{isom alg opers}}}\simeq 
\isom_{\cG,\Op^\reg_\cg} \text{ and }
\ConnD\overset{\text{\eqref{Miura as Con}}}\simeq
\MOp_{\cg,\gen}^\reg$$ with the canonical action of
$\isom_{\cG,\Op^\reg_\cg}$ on $\MOp_{\cg,\gen}^\reg$.
\end{cor}

\ssec{Proof of \propref{char ups alg}}

The deformation $\kappa_\hslash$ of the level defines a
non-commutative deformation $\wh\fH_{\kappa_\hslash}$ of
$\wh\fH_\crit$. This endows $\wh\fH_\crit$ with a Poisson
structure. Moreover, as we shall see in \secref{act of alg on Wak},
the map $\varphi:\fZ_\fg\to \wh\fH_\crit$ is Poisson and has the
following property:

\medskip

The resulting action of $N^*_{\fZ^\reg_\fg/\fZ_\fg}$ on
$\fH^\reg_\crit$ coincides with the one coming from the isomorphism
$\fH^\reg_\crit\simeq \End_{\hg_\crit}(\BW_{\crit,0})$ via the
construction of \secref{rnrm}({\bf 1}).

\medskip

In addition, it is easy to see that the isomorphism, induced by
$\on{map}^M_{\alg}$,
$$\wh\fH_\crit\to \Fun(\ConnDt)$$ respects the Poisson
structures. Hence, the same is true for $\on{map}^M_{\geom}=\tau\circ
\on{map}^M_{\alg}$.

Thus, we obtain a commutative diagram of Poisson ind-schemes:
$$ \CD \Spec(\wh\fH_\crit) @>{\on{map}^M_{\geom}=\tau\circ
\on{map}^M_{\alg}}>> \ConnDt \\ @V{\varphi}VV @V{\on{MT}}VV \\
\Spec(\fZ_\fg) @>{\on{map}_{\geom}=\on{map}_{\alg}}>>
\Op_\cg(\D^\times), \endCD
$$
and the assertion of the proposition follows from \corref{mkdv reg}.

\section{Lie algebroids and Wakimoto modules} \label{renorm Wak}

\ssec{Algebroids acting on categories}     \label{act of alg on cat}

Let $\bA$ be a commutative algebra, and let $\ff$ be a (topological)
Lie algebroid over $\bA$. We shall assume that as an $\bA$-module,
$\ff$ is the dual of a discrete projective $\bA$-module, denoted
$\ff^*$.

Let $\CC$ be a $\bA$-linear category, i.e., $\bA$ acts by
endomorphisms on every object of $\CC$ in a functorial way.  In this
case one can introduce the notion of {\it action of $\ff$ on $\CC$}.
This is, by definition, the same as an action on $\CC$ of the formal
groupoid $\fF$, corresponding to $\ff$. (The latter notion is spelled
out explicitly for groupoids in \cite{Ga}, and for group
ind-schemes in \cite{FG2}, Sect. 22; the generalization to the case of
arbitrary formal groupoids is straightforward.)

A basic example of this situation is when $\CC$ is taken to be the
category of $\bA$-modules.

\medskip

An action of $\ff$ on $\CC$ is specified by the following data. For
every object $\CM\in \CC$ there must be a functorially assigned
extension
$$0\to \CM\underset{\bA}\otimes \ff^*\to \act_\ff^*(\CM)\to \CM\to
0,$$ such that for $\ba\in \bA$, the difference between its action on
$\CM_\ff$ as an object of $\CC$, and the action coming from the
functoriality of $\act^*_\ff$ and the action of $\ba$ on $\CM$, is the
map
$$\CM\mapsto \CM\underset{\bA}\otimes \ff^*$$
given by the image of $d(\ba)$ under the dual of the anchor map
$\Omega^1(\bA)\to \ff^*$.

The functor $\act^*_\ff$ must, in addition, be equipped with a {\it
Lie constraint}, which is a natural transformation that relates the
iteration $\act_\ff^*\circ \act_\ff^*$ with the Lie bracket on
$\ff$. This natural transformation must satisfy an identity for the
3-fold iteration.  We will not spell this out explicitly.

\medskip

Given an action of $\ff$ on $\CC$, we say that an object $\CM\in \CC$
if $\ff$-equivariant if we are given a splitting $\act_\ff^*(\CM)
\leftarrow \CM$, compatible with the Lie constraint.

If $\CC_1$ and $\CC_2$ are two categories, acted on by $\ff$, there
is an evident notion of functor between them, compatible with the
$\ff$-actions.

\ssec{}     \label{gen renorm}

Let us give a typical example of how actions of Lie algebroids 
arise in practice. Let $\bB$
be a topological $\bA$-algebra. Suppose we are given 
a topological Lie algebroid $\ff'$ over $\bA$ that fits into a 
diagram of Lie algebras
$$
\CD
& & \bB \\
& & @AAA \\
0 @>>>  \ff''  @>>>  \ff' @>>>  \ff @>>>  0,
\endCD
$$
and we are given a continuous action
of $\ff'$ on $\bB$ by derivations, which extends the action of $\ff$
on $\bA$ and $\ff''$ on $\bB$.

Note that in this case we can form a topological associative
enveloping algebra, call it $\bB^\ren$, which is universal with
respect to the property that $\bB$ maps to it, as an associative
subalgebra, and $\ff'$, as a Lie algebra, in a compatible
way.

Let $\bB\mod$ denote the category of (discrete, continuous)
$\bB$-modules.

\begin{lem} \label{arise actions}  \hfill

\smallskip

\noindent{\em (1)}
Under the above circumstances we have a canonical action 
of $\ff$ on the category $\bB\mod$.

\smallskip

\noindent{\em (2)} Specifying an $\ff$-equivariant structure on
$\CM\in \bB\mod$ is equivalent to extending the $\bB$-action on it to
a $\bB^\ren$-action.

\end{lem}

We will use the following general construction. 

\medskip

Let $\bB_\hslash$ be a flat $\BC[\hslash]$ family of topological
associative algebras; set $\bB_0=\bB_\hslash/\hslash\cdot \bB_\hslash$,
and let $\bZ$ be the center of $\bB_0$. As in \cite{BD}, Sect. 5.6,
$\bZ$ acquires a natural Poisson bracket. 

Let $\bA\subset \bZ$ be a (closed) subalgebra, closed under the
Poisson bracket, and let $\bI\subset \bA$ be an open Poisson ideal; in
particular $\bA^\reg:=\bA/\bI$ is a discrete Poisson algebra.  We will
assume, in addition, that the multiplication map
$\bI\underset{\bA}{\wh\otimes}\bI\to \bI$ is a closed embedding, and
that $\bI/\bI^2$ is the dual of a projective $\bA^\reg$-module.  Then
$\bI/\bI^2$ is naturally a topological Lie algebroid over
$\bA^\reg$. Denote $\bB^\reg:=\bB_0/\bB_0\cdot \bI$.

\begin{propconstr}  \label{algebroid general}  \hfill

\smallskip

\noindent{\em (1)}
Under the above circumstances, the $\bA^\reg$-linear 
category $\bB^\reg\mod$ carries a natural action of $\bI/\bI^2$.

\smallskip

\noindent{\em (2)}
If $\CM_\hslash$ is a flat family of $\bB_\hslash$-modules, such that
the action of $\bI$ on $\CM_0:=\CM/\hslash\cdot \CM$ is zero,
then, as an object of $\bB^\reg\mod$, $\CM_0$ is naturally
$\bI/\bI^2$-equivariant.

\end{propconstr}

\begin{proof}

Let $\bB^\sharp_\hslash$ be the $\BC[\hslash]$-submodule in the
localization of $\bB_\hslash$ with respect to $\hslash$, formed by
expressions $\frac{\bb_\hslash}{\hslash}$, where
$\bb_0:=\bb_\hslash\on{mod} \hslash\in \bI$. Set $\bB^\sharp=
\bB^\sharp_\hslash/\hslash \cdot \bB^\sharp_\hslash$. It fits into a
short exact sequence
$$0\to \bB_0/\bI\to \bB^\sharp \to \bI\to 0.$$
The algebra $\bA$ acts naturally on $\bB^\sharp$ by left 
multiplication. Set $\bB^\flat:=\bB^\sharp/\bI\cdot \bB^\sharp$.

By the assumption on $\bI$, we have a short exact sequence:
\begin{equation} \label{B flat}
0\to \bB^\reg\to \bB^\flat\to \bI/\bI^2\to 0.
\end{equation}

This brings us to the context of \lemref{arise actions} with
$\ff=\bI/\bI^2$ and $\ff'=\bB^\flat$.

\end{proof}

We will apply the above Proposition-Construction to
topological associative algebras, attached to chiral algebras
at a given point of a curve. We should remark that one
can avoid dealing with topological associative algebras
and instead of using \propconstrref{algebroid general},
one can work directly with chiral algebras,
as was done in \cite{FG1}, Sect. 4:

\medskip

Let $\CB_\hslash$ be a $\BC[\hslash]$-flat family of chiral algebras.
Denote $\CB_0=\CB_\hslash/\hslash\cdot \CB_\hslash$, and
let $\fz(\CB_0)$ be its center. It is known that $\fz(\CB_0)$ acquires a 
natural coisson structure. 

Let $\CA$ be a (commutative) chiral subalgebra of $\fz(\CB_0)$,
and we will assume that $\CA$ is stable under the coisson bracket.
We will assume that $\CA$ is smooth; in particular, the
module $\Omega^1(\CA)$ is projective as a $\CA\otimes \on{D}_X$-module.

Let $\CB_0\mod_\reg$ denote the full subcategory in the
category of chiral $\CB_0$-modules, supported at $x$,
on which the Lie-* action of $\CA$ is zero. By definition,
this is a $\CA_x$-linear category, where $\CA_x$ denotes
the fiber of $\CA$ at $x$. 

We claim that under the above circumstances, we have a 
naturally defined action of 
$$N^*_{\CA_x/\wh\CA_x}\simeq H^0_{DR}(\D_x,\Omega^1(\CA))$$ 
on the category $\CB_0\mod_\reg$.

\medskip

Indeed, by proceeding along the lines of \propconstrref{algebroid
general}, one can produce a Lie-* algebroid over $\CA$, denoted
$\CB^\flat$, that fits into a short exact sequence
$$0\to \CB_0/\CA\to \CB^\flat\to \Omega^1(\CA)\to 0.$$
We define the Lie-* algebroid $\ff'$ to be 
$H^0_{DR}(\D_x,\CB^\flat)$, and then apply \lemref{arise actions}.

\ssec{}   \label{algebroids in context}

Let us return to the setting of the affine algebra at the critical
level.  We will consider various categories over the commutative
algebra $\fZ^\reg_\fg$.

Let us consider several examples of $\fZ^\reg_\fg$-linear categories,
equipped with an action of the algebroid $\isom_{\cG,\fZ_\fg^\reg}$.

\medskip

\noindent{\bf 1.}
Let $\CC$ be $\QCoh(\Fl^\cG_\fZ)$. The action of the
groupoid $\Isom_{\cG,\fZ_\fg^\reg}$ on $\Fl^\cG_\fZ$ defines a natural
action of $\Isom_{\cG,\fZ_\fg^\reg}$ on $\QCoh(\Fl^\cG_\fZ)$.
Hence, we obtain the action on $\QCoh(\Fl^\cG_\fZ)$ of its
algebroid, which is by definition, $\isom_{\cG,\fZ_\fg^\reg}$ .

Evidently, for every $w\in W$ the subcategory 
$\QCoh(\Fl^\cG_{w,\on{th},\fZ})$ also carries an action of 
$\isom_{\cG,\fZ_\fg^\reg}$. (But, of course, this action
does not come from an action of $\Isom_{\cG,\fZ_\fg^\reg}$.)

\medskip

\noindent{\bf 2.}
Let us now take $\CC=\on{D}(\Gr_G)^{\Hecke_\fZ}_\crit\mod$.
It carries a natural action of the groupoid $\Isom_{\cG,\fZ_\fg^\reg}$,
and hence also of the algebroid $\isom_{\cG,\fZ_\fg^\reg}$.

Let us write down this action more explicitly. For
$(\CF,\{\alpha_V\})\in \on{D}(\Gr_G)^{\Hecke_\fZ}_\crit\mod$
we define the corresponding
extension $\act^*_{\isom_{\cG,\fZ_\fg^\reg}}(\CF,\{\alpha_V\})$ as follows.

The underlying D-module is the trivial extension
$$\CF\underset{\fz^\reg_\fg}\otimes (\isom_{\cG,\fZ_\fg^\reg})^*\oplus
\CF,$$ with the $\fZ^\reg_\fg$-action twisted by map
$$\fZ^\reg_\fg\overset{d}\to \Omega^1(\fZ^\reg_\fg)\to
(\isom_{\cG,\fZ_\fg^\reg})^*.$$ The corresponding isomorphisms
$\alpha_V$ are obtained from the original ones by adding the term
involving the $\isom_{\cG,\fZ_\fg^\reg}$-action on $\CV_\fZ$.

\medskip

By construction, the functor $\sE^\fZ:\QCoh(\Fl^\cG_\fZ)\to
\on{D}(\Gr_G)^{\Hecke_\fZ}_\crit\mod$, defined in \secref{sE Z},
respects the action of the groupoid $\Isom_{\cG,\fZ_\fg^\reg}$,
and hence of the algebroid $\isom_{\cG,\fZ_\fg^\reg}$.

One can see the latter compatibility explicitly for quasi-coherent 
sheaves $\CT$ that are direct images from affine open subsets
of $\Fl^\cG_\fZ$, using the explicit description of 
$\sE^\fZ(\CT)$ in this case as a co-equalizer, see \secref{sE Z}.

\medskip

\noindent{\bf 3.}  Finally, let us take $\CC=\hg_\crit\mod_\reg$. This
example fits into the framework of (the chiral algebra version of)
\propconstrref{algebroid general} for the $\BC[\hslash]$-family
$\CB_\hslash:=\CA_{\fg,\kappa_\hslash,X}$, and
$\CA=\fz(\CA_{\fg,\crit,X})$.

We obtain an action of $N^*_{\fZ^\reg_\fg/\fZ_\fg}$ on
$\hg_\crit\mod_\reg$, and hence also of $\isom_{\cG,\fZ_\fg^\reg}$, via
the identification $\upsilon_{\geom}$.

Note that the topological associative algebra $\bB^\ren$ identifies,
by definition, with the renormalized algebra
$U^{\ren,\reg}(\hg_\crit)$.

\medskip

As was shown in \cite{FG4}, Sect. 2.9., the functor
$$\Gamma^{\Hecke_\fZ}: \on{D}(\Gr_G)^{\Hecke_\fZ}_\crit\mod\to
\hg_\crit\mod_\reg$$ is compatible with the
$\isom_{\cG,\fZ_\fg^\reg}$-actions. From this we obtain that the
functor
$$\sG=\Gamma^{\Hecke_\fZ}\circ \sE:\QCoh(\Fl^\cG_\fZ)\to
\hg_\crit\mod_\reg$$ is also compatible with the
$\isom_{\cG,\fZ_\fg^\reg}$-actions.

\ssec{}   \label{act of alg on Wak funct}

Recall the topological algebra $\wh\fH_\crit$ and consider the
category $\Bigl(\wh\fH_\crit\underset{\fZ_\fg}\otimes
\fZ^\reg_\fg\Bigr)\mod$, from which, we recall, we had the Wakimoto
functor to the category $\hg_\crit\mod_\reg$.  We will establish the
following:
\begin{propconstr}   \label{action of renorm on Wak}  \hfill

\smallskip

\noindent{\em (1)} The category
$\Bigl(\wh\fH_\crit\underset{\fZ_\fg}\otimes \fZ^\reg_\fg\Bigr)\mod$
carries a natural action of $N^*_{\fZ^\reg_\fg/\fZ_\fg}$.

\smallskip

\noindent{\em (2)} The functor
$\BW:\Bigl(\wh\fH_\crit\underset{\fZ_\fg}\otimes
\fZ^\reg_\fg\Bigr)\mod\to \hg_\crit\mod_\reg$ is compatible with the
$N^*_{\fZ^\reg_\fg/\fZ_\fg}$-actions.
\end{propconstr}

The rest of this subsection is devoted to the proof of this
\propref{action of renorm on Wak}. To carry it out 
we need to review the framework in which Wakimoto modules were defined.
We will follow the conventions of \cite{FG2}, Sect. 10.

\medskip

For every level $\kappa$ we have the Heisenberg chiral algebra
$\fH_{\kappa,X}$, and a chiral algebra, denoted
$\fD^{ch}(\oGB)_{\kappa,X}$, which is isomorphic to the tensor product
$$\fD^{ch}(N)_X\otimes \fH_{\kappa,X},$$ where $\fD^{ch}(N)_X$ is the
chiral algebra of differential operators on $N$, which is independent
of the level.

In particular, for every $\fH_{\kappa,X}$-module $\CT$, supported at
$x\in X$, we can consider the $\fD^{ch}(\oGB)_{\kappa,X}$-module
$\fD^{ch}(N)_x\otimes \CT$, where $\fD^{ch}(N)_x$ denotes the vacuum
module of $\fD^{ch}(N)_X$ at $x$.

\medskip

We have a canonical bosonization map $\fl_\kappa:\CA_{\fg,\kappa,X}\to
\fD^{ch}(\oGB)_{\kappa,X}$. The Wakimoto functor associates at $\CT$
as above the module $\fD^{ch}(N)_x\otimes \CT$, regarded as a
$\CA_{\fg,\kappa,X}$-module via $\fl_\kappa$.

\medskip

Let us take $\kappa=\kappa_\crit$. In this case $\fH_{\crit,X}$ equals
the center of $\fD^{ch}(\oGB)_{\crit,X}$. Moreover, we have the
following commutative diagram:

\begin{equation} \label{two maps of center}
\CD
\fz_{\fg,X}   @>>> \CA_{\fg,\crit,X}\\
@V{\varphi}VV    @V{\fl_\crit}VV  \\
\fH_{\crit,X}  @>>>  \fD^{ch}(\oGB)_{\crit,X}.
\endCD
\end{equation}

We apply the construction of \secref{gen renorm} for $\CB_{\hslash}$
being each of the following three chiral algebras:
$\CA_{\fg,\kappa_\hslash,X}$, $\fD^{ch}(\oGB)_{\kappa_\hslash,X}$ and
$\fH_{\kappa_\hslash,X}$.  (In what follows, we will replace the
subscript $\kappa_\hslash$ by just $\hslash$, for brevity.)

In each of these cases we take $\CA\subset \fz(\CB_0)$ to be the image
of $\fz_{\fg,X}:=\fz(\CA_{\fg,\crit,X})$. The functoriality of the
construction in \secref{gen renorm} makes the assertion of
\propconstrref{action of renorm on Wak} manifest.

\ssec{}  \label{act of alg on Wak}

Let us make several additional remarks on the above construction. 

First, we see explicitly that the homomorphism $\varphi:\fZ_\fg\to
\wh\fH_\crit$ is indeed a Poisson map, and the action of
$N^*_{\fZ^\reg_\fg/\fZ_\fg}$ on the category
$$\Bigl(\wh\fH_\crit\underset{\fZ_\fg}\otimes
\fZ^\reg_\fg\Bigr)\mod\simeq \fH_{\crit,X}\mod_\reg$$ corresponds to
the action of the algebroid $N^*_{\fZ^\reg_\fg/\fZ_\fg}$ on the
ind-scheme $$\Spec(\wh\fH_\crit)\underset{\Spec(\fZ_\fg)}\times
\Spec(\fZ^\reg_\fg).$$

\medskip

Secondly, assume that $\CT_\hslash$ was a flat family of 
$\fH_{\hslash,X}$-modules, such that the action of $\fZ_\fg$ on
$\CT_0$ factors through $\fZ^\reg_\fg$. 

On the one hand, by \propconstrref{algebroid general}(2), $\CT_\crit$
is naturally $N^*_{\fZ^\reg_\fg/\fZ_\fg}$-equivariant, as an object of
$\Bigl(\wh\fH_\crit\underset{\fZ_\fg}\otimes \fZ^\reg_\fg\Bigr)\mod$.
Therefore, $\BW(\CT_\crit)$ is
$N^*_{\fZ^\reg_\fg/\fZ_\fg}$-equivariant as an object of
$\hg_\crit\mod_\reg$, and hence carries an action of
$U^{\ren,\reg}(\hg_\crit)$.

On other hand, we can consider the family
$\CM_\hslash:=\BW(\CT_\hslash)$ as in \secref{rnrm}({\bf 3}), and
hence $\CM_\crit\simeq \BW(\CT_\crit)$ acquires a
$U^{\ren,\reg}(\hg_\crit)$-action.

The following assertion results from the definitions:

\begin{lem}   \label{two actions of renorm}
Under the above circumstances, the two
$U^{\ren,\reg}(\hg_\crit)$-actions on $\BW(\CT_\crit)$ coincide.
\end{lem}

\ssec{}   

Let us recall the isomorphism of ind-schemes \eqref{induction parameters}
$$\Spec(\wh\fH_\crit)\underset{\Spec(\fZ_\fg)}\times
\Spec(\fZ^\reg_\fg)\simeq \underset{w\in W}\sqcup\,
\Fl^\cG_{w,\on{th},\fZ}.$$

On the one hand, we have the action of $N^*_{\fZ^\reg_\fg/\fZ_\fg}$
on the LHS, described in \secref{act of alg on Wak}. On the other hand,
we have a natural action of $\isom_{\cG.\fZ^\reg_\fg}$ on the RHS.

\begin{prop}   \label{reg mkdv}
Under the identification
$\upsilon_{\geom}:N^*_{\fZ^\reg_\fg/\fZ_\fg}\simeq
\isom_{\cG.\fZ^\reg_\fg}$ the above actions on the two sides of
\eqref{induction parameters} coincide.
\end{prop}

\begin{proof}

Using \thmref{two maps of algebroids}, we obtain that the assertion of
the proposition is equivalent to the following generalization of
\lemref{mkdv reg}:

The action of $N^*_{\Op^\reg_\cg/\Op_\cg(\D^\times)}$ on
$\ConnDt\underset{\Op_\cg(\D^\times)}\times \Op_\cg^\reg$, resulting
from the Poisson structure on $\ConnDt$ and the map
$\on{MT}:\ConnD\to \Op^\reg_\cg$, identifies via
$$N^*_{\Op^\reg_\cg/\Op_\cg(\D^\times)}\overset{\text{\eqref{isom alg
opers}}} \simeq \isom_{\cG,\Op^\reg_\cg} \text{ and }
\ConnDt\underset{\Op_{\cg}(\D^\times)}\times \Op^\reg_\cg
\overset{\text{\eqref{Miura as Con, general}}}\simeq \underset{w\in
W}\sqcup\, \MOp^{w,\on{th},\reg}_{\cg}$$ with the natural
$\isom_{\cG,\Op^\reg_\cg}$-action on $\MOp^\reg_\cg$.

\medskip

Using \lemref{mkdv flows}, this is, in turn, equivalent to the fact
that the canonical action of $\isom_{\cG,\Op_\cG(\D^\times)}$ on
$\MOp_{\cg,\gen}(\D^\times)$ is such that the induced action of
$\isom_{\cG,\Op^\reg_\cg}$ on
$\MOp_{\cg,\gen}(\D^\times)\underset{\Op_{\cg}(\D^\times)}\times
\Op^\reg_\cg$ is compatible with the action of
$\isom_{\cG,\Op_\cG(\D^\times)}$ on $\MOp^\reg_\cg$ under the map
$$\MOp_{\cg,\gen}(\D^\times)\underset{\Op_{\cg}(\D^\times)}\times
\Op^\reg_\cg \hookrightarrow \MOp^\reg_\cg.$$

The required compatibility follows from the construction of the latter
map in \cite{FG2}, Sect. 3.6.

\end{proof}

\ssec{The functor $\sG$ and algebroids}

Let us now return to the setting of \thmref{full ident}. Recall that
by \secref{more general} for every $w\in W$ there exists a functor
$_w\sQ:\QCoh(\Fl^\cG_{w,\on{th},\fZ})\to \QCoh(\Fl^\cG_{w,\on{th},\fZ})$
and an isomorphism
$$_w\BW\simeq \sG\circ {}_w\sQ.$$

We claim now that the functor $_w\sQ$ is compatible with the action of
$\isom_{\cG.\fZ^\reg_\fg}$. Indeed, this follows from the
corresponding properties of the functors $_w\BW$ and $\sG$,
established above, and \propref{isom on Ext}.

\medskip

Recall now that the functor $_w\sQ$ has the form
$$\CT\mapsto
\CL^{\on{twist}}_w\underset{\Fun(\Fl^\cG_{w,\on{th},\fZ})}\otimes\CT,$$
for some $\Fun(\Fl^\cG_{w,\on{th},\fZ})$-bimodule
$\CL^{\on{twist}}_w$.  Moreover, we have shown that
$\CL^{\on{twist}}_w$, regarded as a right
$\Fun(\Fl^\cG_{w,\on{th},\fZ})$-module, is a line bundle.

\begin{cor}  \hfill

\smallskip

\noindent{\em (1)} The left and right actions of
$\Fun(\Fl^\cG_{w,\on{th},\fZ})$ on $\CL^{\on{twist}}_w$ coincide.

\smallskip

\noindent{\em (2)} There exists a (unique)
$\isom_{\cG.\fZ^\reg_\fg}$-equivariant structure on
$\CL^{\on{twist}}_w$, which induces the above
$\isom_{\cG.\fZ^\reg_\fg}$-equivariant structure on $_w\sQ$.

\end{cor}

\begin{proof}

Let us choose a trivialization of $\CL^{\on{twist}}_w$ as a right
$\Fun(\Fl^\cG_{w,\on{th},\fZ})$-module. Then the left action of
$\Fun(\Fl^\cG_{w,\on{th},\fZ})$ defines a homomorphism ${\mathbf
\gamma}:\Fun(\Fl^\cG_{w,\on{th},\fZ})\to
\Fun(\Fl^\cG_{w,\on{th},\fZ})$.  The assertion of the first point of
the corollary is equivalent to the fact that ${\mathbf \gamma}$ equals
the identity map.

However, the fact that the functor 
${\mathbf \gamma}^*:\QCoh(\Fl^\cG_{w,\on{th},\fZ})\to
\QCoh(\Fl^\cG_{w,\on{th},\fZ})$ is compatible with the
$\isom_{\cG.\fZ^\reg_\fg}$-action on this category implies that
${\mathbf \gamma}$ commutes with the action of this algebroid. In
addition, from \thmref{ident w Wak}, we know that ${\mathbf
\gamma}|_{\Fl^\cG_{w,\fZ}}$ equals the identity map.

Hence, ${\mathbf \gamma}$ is the identity map, restricted to the
minimal $\isom_{\cG.\fZ^\reg_\fg}$-stable subscheme of 
$\Fl^\cG_{w,\on{th},\fZ}$, that contains
$\Fl^\cG_{w,\fZ}$. But since this minimal subscheme is
$\Fl^\cG_{w,\on{th},\fZ}$ itself, assertion (1) of the lemma follows,
while assertion (2) is evident.

\end{proof}

Thus, to finish the proof of \thmref{full ident}, it is sufficient to show
the following:

\begin{thm} \label{determine line bundle}
There exists a canonical isomorphism of line bundles
$$\CL^{\on{twist}}_w\simeq
\CL^{\crho-w(\crho)}_{\Fl^\cG_{w,\on{th},\fZ}},$$ that respects the
$\isom_{\cG,\fZ^\reg_\fg}$-equivariant structures.
\end{thm}

We recall that the line bundle on the RHS is by definition the
restriction of the line bundle $\CL^{\crho-w(\crho)}_{\Fl^\cG_{\fZ}}$
to $\Fl^\cG_{w,\on{th},\fZ}$; the
$\isom_{\cG.\fZ^\reg_\fg}$-equivariant structure on it comes from the
equivariant structure on $\CL^{\crho-w(\crho)}_{\Fl^\cG_{\fZ}}$ with
respect to the groupoid $\Isom_{\cG.\fZ^\reg_\fg}$.

\ssec{Some corollaries} Let $\CT\in \QCoh(\Fl^\cG_{w,\on{th},\fZ})$ be
$\isom_{\cG,\fZ^\reg_\fg}$- equivariant. Then by \propconstrref{action
of renorm on Wak}, the corresponding Wakimoto module $_w\BW(\CT)$ will
carry an action of the renormalized algebra
$U^{\ren,\reg}(\hg_\crit)$.

Considering $\CL^{\crho-w(\crho)}_{\Fl^\cG_{w,\on{th},\fZ}}
\underset{\Fun(\Fl^\cG_{w,\on{th},\fZ})}\otimes\CT$ with its 
$\isom_{\cG,\fZ^\reg_\fg}$-equivariant structure, from
\secref{act of alg on cat} we obtain that 
$\sG\Bigl(\CL^{\crho-w(\crho)}_{\Fl^\cG_{w,\on{th},\fZ}}
\underset{\Fun(\Fl^\cG_{w,\on{th},\fZ})}\otimes \CT\Bigr)$ also
carries an action of $U^{\ren,\reg}(\hg_\crit)$.

Now, the fact that the isomorphism of \thmref{determine line bundle}
is $\isom_{\cG,\fZ^\reg_\fg}$-equivariant implies that the isomorphism
$$_w\BW(\CT)\simeq
\sG\Bigl(\CL^{\crho-w(\crho)}_{\Fl^\cG_{w,\on{th},\fZ}}
\underset{\Fun(\Fl^\cG_{w,\on{th},\fZ})}\otimes\CT\Bigr)$$ of 
\thmref{full ident} respects the $U^{\ren,\reg}(\hg_\crit)$-actions.

\medskip

As an example let us take $\CT=\CO_{\Fl^\cG_{1,\fZ}}$. We obtain:

\begin{cor}
The isomorphism
$$\Gamma^{\Hecke_\fZ}(\CW^\fZ_1)\simeq \BW_{\crit,0}$$ of
\thmref{ident prince Wak} respects the
$U^{\ren,\reg}(\hg_\crit)$-actions.
\end{cor}

Let us give a direct proof of this corollary.

\begin{proof}

By construction, $\CW^\fZ_1$, as an
$\isom_{\cG,\fZ^\reg_\fg}$-equivariant object of
$\on{D}(\Gr_G)^{\Hecke_\fZ}_\crit\mod$, is a quotient of the direct
sum of objects of the form
$\on{Ind}^{\Hecke_\fZ}(j_{\cla,*,\Gr})\underset{\fZ^\reg_\fg}\otimes
\CL^{-\cla}_{\Fl^{\cG}_{1,\fZ}}$. Hence, it is enough to show that
each of the morphisms
$$\Gamma(\Gr_G,j_{\cla,\Gr_G,*})\underset{\fZ^\reg_\fg}\otimes
\CL^{-\cla}_{\Fl^{\cG}_{1,\fZ}}\to \BW_{\crit,0}$$ is compatible with
the $U^{\ren,\reg}(\hg_\crit)$-actions.

However, the latter morphism is by definition the composition
\begin{align*}
&\Gamma(\Gr_G,j_{\cla,\Gr_G,*})\underset{\fZ^\reg_\fg}\otimes
\CL^{-\cla}_{\Fl^{\cG}_{1,\fZ}}\simeq \Bigl(j_{\cla,*}\star
\BV_\crit\Bigr) \underset{\fZ^\reg_\fg}\otimes
\CL^{-\cla}_{\Fl^{\cG}_{1,\fZ}}\overset{\phi} \to \Bigl(
j_{\cla,*}\star \BW_{\crit,0}\Bigr)\underset{\fZ^\reg_\fg}\otimes
\CL^{-\cla}_{\Fl^{\cG}_{1,\fZ}}\simeq \\ &\simeq
\Bigl(\BW_{\crit,0}\underset{\fH^\reg_\crit}\otimes
\CL^{\cla}_{\fH}\Bigr) \underset{\fZ^\reg_\fg}\otimes
\CL^{-\cla}_{\Fl^{\cG}_{1,\fZ}}\twoheadrightarrow \BW_{\crit,0},
\end{align*}
where all the arrows are compatible with the
$U^{\ren,\reg}(\hg_\crit)$-actions (see \secref{proof two actions on
H}).

\end{proof}

Let us return to the general setting of Theorems \ref{full ident} and
\ref{determine line bundle}, and consider the semi-infinite cohomology
functor
$$\CM\mapsto H^\semiinf(\fn((t)),\fn[[t]],\CM\otimes \Psi_0):
\hg_\crit\mod_\reg\to \fZ^\reg_\fg\mod.$$ Recall (see Section 18.3 of
\cite{FG2}) that this functor is naturally compatible with the
$\isom_{\cG,\fZ^\reg_\fg}$-actions on both categories.

Note, that on the other hand, by \propref{semiinf of convolution}, the
composite functor
$$\QCoh(\Fl^\cG_{w,\on{th},\fZ}) \overset{\sG}\to
\hg_\crit\mod_\reg\to \fZ^\reg_\fg\mod$$ is isomorphic to the identity
functor. Moreover, from the construction, it is easy to see that this
isomorphism respects the $\isom_{\cG,\fZ^\reg_\fg}$-equivariant
structures.

On the other hand, by \eqref{naive semiinf of Wak}, we have a
functorial (but dependent on some choices) isomorphism
$$H^\semiinf(\fn((t)),\fn[[t]],{}_w\BW(\CT)\otimes \Psi_0)\simeq \CT$$
for $\CT\in \QCoh(\Fl^\cG_{w,\on{th},\fZ})$.

Thus, combining this with the assertion of \thmref{full ident},
we obtain a functorial isomorphism
$$\CT\simeq \CL^{\crho-w(\crho)}_{\Fl^\cG_{w,\on{th},\fZ}}
\underset{\Fun(\Fl^\cG_{w,\on{th},\fZ})}\otimes\CT,$$ i.e., a
trivialization of the line bundle
$\CL^{\crho-w(\crho)}_{\Fl^\cG_{w,\on{th},\fZ}}$ over
$\Fl^\cG_{w,\on{th},\fZ}$. Let us explain what this trivialization is.

\medskip

First, let us observe that the choice of
a character $\Psi_0$ amounts to a trivialization of the $H$-torsor
$\omega_x^\crho$. However, since different choices of $\Psi_0$ are
$H[[t]]$-conjugate, we obtain that when we apply the functor
$H^\semiinf(\fn((t)),\fn[[t]],?\otimes \Psi_0)$ to $\hg_\crit$-modules
that are $I$-equivariant, the spaces that we obtain for different
choices of $\Psi_0$ are canonically isomorphic. In other words, on the
category $\hg_\crit\mod^I$, the functor
$H^\semiinf(\fn\ppart,\fn[[t]],?\otimes \Psi_0)$ is well-defined.

For the same reason, if $\CM\in \hg_\crit\mod^I$ is in addition 
equivariant with respect to $\Aut(\D)$, then the cohomology
$H^\semiinf(\fn\ppart,\fn[[t]],?\otimes \Psi_0)$ acquires a 
$\Aut(\D)$-action. 

\medskip

Canonically, we have an isomorphism
\begin{equation} \label{less naive semiinf of Wak}
H^\semiinf(\fn\ppart,\fn[[t]],{}_w\BW(\CT)\otimes \Psi_0)\simeq
\CT\otimes \omega_x^{\langle w(\rho)-\rho,\crho\rangle}.
\end{equation}

Thus, rather than trivializing the line bundle
$\CL^{\crho-w(\crho)}_{\Fl^\cG_{w,\on{th},\fZ}}$ we need to identify
it with $\Fun(\Fl^\cG_{w,\on{th},\fZ})\otimes \omega_x^{\langle
w(\rho)-\rho,\crho\rangle}$ in a $\Aut(\D)$-equivariant way. It is
easy to see that there exists a unique such identification, with
induces the isomorphism \eqref{ident tors, special} over
$\Fl^\cG_{w,\fZ}$.

\ssec{Proof of \thmref{determine line bundle}}

Recall the algebroid $\isom^\nilp_{\Op_\cg}$ over the scheme
$\Op_\cg^\nilp$ (see \cite{FG2}, Sect. 4.5), and let
$\isom_{\fZ^\nilp_\fg}$ denote the corresponding algebroid over
$\Spec(\fZ^\nilp_\fg)$.

By {\it loc. cit.}, it preserves the subscheme
$\Spec(\fZ^\reg_\fg)\subset \Spec(\fZ^\nilp_\fg)$. Hence, the
restriction $\isom_{\fZ^\nilp_\fg}|_{\Spec(\fZ^\reg_\fg)}$ has a
natural structure of Lie algebroid; we will denote it by
$\isom_{\cB,\fZ^\reg_\fg}$.  By \cite{FG2}, Sect. 4.5,
$\isom_{\cB,\fZ^\reg_\fg}$ identifies with the sub-algebroid of
$\isom_{\cG,\fZ^\reg_\fg}$, corresponding to the reduction of
$\CP_{\cG,\fZ}$ to $\cB$.

\medskip

Since $\cB$ acts transitively along the fibers of
$\Fl^{\cG}_{w,\fZ}\to \Spec(\fZ^\reg_\fg)$, and since the stabilizers
of any point of $\Fl^{\cG}_{w,\fZ}$ in $\isom_{\cG,\fZ^\reg_\fg}$ and
$\isom_{\cB,\fZ^\reg_\fg}$ coincide, it suffices to show that there
exists a $\isom_{\cB,\fZ^\reg_\fg}$-equivariant isomorphism
\begin{equation} \label{two line bundles}
\CL^{\on{twist}}_w|_{\Fl^{\cG}_{w,\fZ}}\simeq 
\CL^{\crho-w(\crho)}_{\Fl^{\cG}_{w,\fZ}}
\end{equation}
of line bundles on $\Fl^{\cG}_{w,\fZ}$.

\medskip

Recall that along with the renormalized algebra
$U^{\ren,\reg}(\hg_\crit)$, there exists its version
$U^{\ren,\nilp}(\hg_\crit)$, corresponding to the quotient
$\fZ^\nilp_\fg$ of $\fZ_\fg$ (see \cite{FG2}, Sect. 7.4).

As in \propconstrref{action of renorm on Wak}, we obtain that if
$$\CT\in \QCoh\Bigl(\Spec\bigl(\wh\fH_{\crit}\underset{\fZ_\fg}\otimes
\fZ^\nilp_\fg\bigr)\Bigr),$$ is
$\isom_{\cB,\fZ^\nilp_\fg}$-equivariant, then the Wakimoto module
$\BW(\CT)$ acquires an action of the renormalized algebra
$U^{\ren,\nilp}(\hg_\crit)$.

Similarly, if $\CT\in \QCoh(\Fl^{\cG}_\fZ)$ is
$\isom_{\cB,\fZ^\reg_\fg}$-equivariant, then $\sG(\CT)$ also acquires
an action of $U^{\ren,\nilp}(\hg_\crit)\underset{\fZ^\nilp_\fg}\otimes
\fZ^\reg_\fg$.

\medskip

Recall now the isomorphism
\begin{equation} \label{w again}
\Gamma^{\Hecke_\fZ}(\Gr_G,\CW^\fZ_w)\simeq
\BW_{\crit,w(\rho)-\rho,\reg}\otimes \omega_x^{\langle
\rho-w(\rho),\crho\rangle} \overset{\text{\eqref{ident tors,
special}}}\simeq
\BW\bigl(\CL^{w(\crho)-\crho}_{\Fl^{\cG}_{w,\fZ}}\bigr)
\end{equation}
given by \thmref{ident w Wak}. By the definition of
$\CL^{\on{twist}}_w$, we obtain an isomorphism of line bundles
appearing in \eqref{two line bundles}.  The fact that this isomorphism
respects the $\isom_{\cB,\fZ^\reg_\fg}$-action is equivalent to the
fact that the isomorphism of \eqref{w again} respects the action of
$U^{\ren,\nilp}(\hg_\crit)$.

This is, in turn, equivalent to the map
$$\BV_{\crit}\to
\BW\bigl(\CL^{w(\crho)-\crho}_{\Fl^{\cG}_{w,\fZ}}\bigr),$$
of \eqref{Vak to Wak} being compatible with the
$U^{\ren,\nilp}(\hg_\crit)$-action.  Recall that the latter morphism
comes by adjunction from a map
\begin{equation} \label{renorm action on Av}
\BV_{\crit}\to \on{Av}_{G[[t]]/I}
\Bigl(\BW\bigl(\CL^{w(\crho)-\crho}_{\Fl^{\cG}_{w,\fZ}}
\bigr)\Bigr).
\end{equation} 
Thus, we need to show that the map in \eqref{renorm action on Av} is
compatible with the action of $U^{\ren,\nilp}(\hg_\crit)$, where on
the RHS it is defined by functoriality as in \secref{rnrm}({\bf 4}).
The latter essentially follows from the definition of the 
map \eqref{Vak to Wak}:

\medskip

Consider the family of Wakimoto modules 
$\BW_{\hslash,w(\rho)-\rho}$,
which is, by definition, induced from the
$\fH_{\hslash,X}$-module
$$\pi_{\hslash,w_0(w(\rho)-\rho)} :=
\on{Ind}^{\wh\fh_\hslash}_{\fh[[t]]}
(\BC^{w_0(w(\rho)-\rho)}).$$

Its fiber at the critical level
$$\pi_{\crit,w_0(w(\rho)-\rho)}\simeq
\fH_\crit^{\RS,w_0(w(\rho)-\rho)}$$ is supported over
$\Spec\bigl(\wh\fH_{\crit}\underset{\fZ_\fg}\otimes
\fZ^\nilp_\fg\bigr)$, hence, it acquires an action of
$\isom_{\fZ^\nilp_\fg}$. By the $\nilp$-version of
\propconstrref{action of renorm on Wak}, the module
$\BW_{\crit,w(\rho)-\rho}$ is endowed an action of
$U^{\ren,\nilp}(\hg_\crit)$, and by \lemref{two actions of renorm},
this action coincides with the one coming from the family
$\BW_{\hslash,w(\rho)-\rho}$.

\medskip

We have an $\hslash$-family of maps
$$\BV_{\hslash}\to
h^{\ell(w)}\Bigl(\on{Av}_{G[[t]]/I}(\BW_{\hslash,w(\rho)-\rho})\Bigr)$$
and, hence, the correspoding map at the critical level
\begin{equation} \label{av one}
\BV_\crit\to
h^{\ell(w)}\Bigl(\on{Av}_{G[[t]]/I}(\BW_{\crit,w(\rho)-\rho})\Bigr)
\end{equation}
is compatible with the $U^{\ren,\nilp}(\hg_\crit)$-action. A
compatibility of constructions in points {\bf 3} and {\bf 4}
of \secref{rnrm} shows that the $U^{\ren,\nilp}(\hg_\crit)$-action
on the RHS of \eqref{av one} coincides
with the one coming by functoriality from its action on
$\BW_{\crit,w(\rho)-\rho}$ via 
\propconstrref{action of renorm on Wak}.

\medskip

Let $\imath_w$ denote the closed embedding
$$\Fl^\cG_{w,\fZ}\hookrightarrow
\Spec(\fH_\crit^{\RS,w_0(w(\rho)-\rho)}).$$ We have a canonical
$\isom_{\cB,\fZ^\nilp_\fg}$-equivariant isomorphism:
$$\on{Av}_{G[[t]]/I}
\Bigl(\BW\bigl(\CL^{w(\crho)-\crho}_{\Fl^{\cG}_{w,\fZ}}\bigr)\Bigr)
\simeq\on{Av}_{G[[t]]/I}(\BW_{\crit,w(\rho)-\rho})
\underset{\fH_\crit^{\RS,w_0(w(\rho)-\rho)}}{\overset{L}\otimes}
\imath_w{}_*\Bigl(\CL^{w(\crho)-\crho}_{\Fl^{\cG}_{w,\fZ}}\Bigr).$$

The map \eqref{renorm action on Av} is, by construction, obtained from
the map \eqref{av one} and the identification
$$L^{\ell(w)}\imath^*_w\Bigl(\imath_w{}_*\Bigl(\Fun(\Fl^\cG_{w,\fZ})\Bigr)
\Bigr)\simeq \CL^{\crho-w(\crho)}_{\Fl^{\cG}_{w,\fZ}}.$$ Thus, we need
to see that the latter is compatible with the
$\isom_{\cB,\fZ^\nilp_\fg}$-actions.

\medskip

Translating it using the isomorphism $\on{map}^M_{\geom}$, we obtain
that the LHS of the above expression identifies with the line bundle
$$\Lambda^{\ell(w)}\Bigl(N^*_{\MOp^{w,\reg}_\cg/\MOp^{w,\nilp}_\cg}\Bigr)$$
over $\MOp_\cg^{w,\reg}$. By Theorem 3.6.2 of \cite{FG2}, the latter
does indeed identify with $\CL^{\crho-w(\crho)}_{\MOp_\cg^{w,\reg}}$
in a $\isom^\nilp_{\Op_\cg}$-equivariant way, as required.

\section{Renormalized Wakimoto modules}    \label{resolution}

In this section we will study a particular family of 
$\hg_\crit\mod$-modules, obtained by a certain 
renormalization procedure. We will show that they 
coincide with Wakimoto modules corresponding to some 
particular quasi-coherent sheaves on the formal
neighborhoods of the Schubert strata in $\Fl^{\cG}_\fZ$.

\ssec{}

Recall the context of \secref{gen renorm}. Let 
$\bB_\hslash\mod$ be the category of all (discrete)
$\bB_\hslash$-modules, and let
$\bB_\hslash\mod^{fl}\subset \bB_\hslash\mod$ be the full
subcategory consisting of $\BC[\hslash]$-flat modules. Let 
\begin{equation} \label{embed cat}
\bB_\hslash\mod^{fl,\reg}\subset \bB_\hslash\mod^{fl}
\end{equation}
be the full subcategory, consisting of modules $\CN_\hslash$, for
which the action of $\bA$ on $\CN_0:=\CN_\hslash/\hslash\cdot
\CN_\hslash$ factors through $\bA^\reg$. (Recall that in this case
$\CN_0$ acquires a natural $\bB^\ren$-action.)

\begin{propconstr}   \label{renorm module}
The tautological embedding \eqref{embed cat} admits a left
adjoint.
\end{propconstr}

\begin{proof}

For $\CM_\hslash\in \bB_\hslash\mod^{fl}$ let $\CM^\sharp_\hslash$
denote its modification spanned by the symbols
$\frac{m_\hslash}{\hslash}$ where $m_\hslash$ is such that $m_0\in
\bI\cdot \CM_0$. We have a canonical map $\CM_\hslash\to
\CM^\sharp_\hslash$ and a short exact sequence
$$0\to \CM_0/\bI\cdot \CM_0\to \CM^\sharp_0\to \bI\cdot \CM_0\to 0.$$

Let us denote by $\CM_\hslash\mapsto \CM^{k \cdot \sharp}_\hslash$ the
$k$-th iteration of the functor $\CM_\hslash\mapsto
\CM^{\sharp}_\hslash$. The fiber $\CM^{k\cdot \sharp}_0$ admits a
$k+1$-term filtration $\left(\CM^{k\cdot \sharp}_0\right)_i$, such
that
$$
\begin{cases}
& \on{gr}^j\left(\CM^{k\cdot \sharp}_0\right)\simeq \bI^{j-1}\cdot
\CM_0/\bI^j\cdot \CM_0,\,\, 1\leq j\leq k \text{ and }\\ &
\on{gr}^{k+1}\left(\CM^{k\cdot \sharp}_0\right)\simeq \bI^k\cdot \CM_0.
\end{cases}
$$ Moreover, the submodule 
$\left(\CM^{k\cdot\sharp}_0\right)_k\subset \CM^{k\cdot \sharp}_0$ 
is annihilated by $\bI$.

\medskip

Thus, we obtain a sequence of maps
$$...\to \CM_\hslash^{(k-1)\cdot \sharp} \to \CM^{k\cdot
  \sharp}_\hslash\to \CM^{(k+1)\cdot \sharp}_\hslash\to...$$ and we
  set
$$\CM_\hslash^\ren:=\underset{k}{\underset{\longrightarrow}{\lim}}\,
\CM^{k\cdot \sharp}_\hslash.$$ The above computation of fibers of
$\CM^{k\cdot \sharp}_\hslash$ implies that $\CM_\hslash^\ren$ belongs
to $\bB_\hslash\mod^{\reg,fl}$. It satisfies the required adjunction
property by construction.

\end{proof}

For $\CM_\hslash$ as above, let $\CM^\ren_0$ denote the $\bB^\reg$-module
$\CM_\hslash^\ren/\hslash\cdot \CM_\hslash$. By construction, it carries an
action of $B^\ren_0$. As a $\bB^\reg$-module it is equipped with an
increasing filtration, labeled by positive integers, with
\begin{equation} \label{gr renorm}
\on{gr}^j(\CM_0^\ren)\simeq \bI^{j}\cdot \CM_0/\bI^{j+1}\cdot \CM_0.
\end{equation}

\ssec{}

Let us consider the family of chiral algebras $\fH_{\hslash,X}$; let
$\wh\fh_{\hslash}$ be the central extension of $\fh\ppart$,
corresponding to the point $x\in X$, and let
$\bB_\hslash:=\wh\fH_{\hslash}$ be the corresponding family of
associative topological algebras. As in \secref{act of alg on Wak funct},
we let $\bI\subset \bB_0:=\wh\fH_\crit$ be the image of 
$\on{ker}(\fZ_\fg\to \fZ^\reg_\fg)$ under the map
$\varphi:\fZ_\fg\to \wh\fH_\crit$.

\medskip

For $w\in W$ let $\pi_{\hslash,w}$ denote the
$\wh\fH_{\hslash}$-module
$\on{Ind}_{\fh[[t]]}^{\wh\fh_{\hslash}}(\BC^{w_0(w(\rho)-\rho)})$.
We let $\pi_{\crit,w}$ denote the fiber of $\pi_{\hslash,w}$ at the
critical level.
Let $\pi^\ren_{\hslash,w}$ be the $\wh\fH_{\hslash}$-module, 
corresponding to $\pi_{\hslash,w}$ via \propconstrref{renorm module}. 

Let $\pi^\ren_{\crit,w}$ denote the fiber of $\pi^\ren_{\hslash,w}$ at
the critical level.  By \propconstrref{algebroid general},
$\pi^\ren_{\crit,w}$ is equivariant with respect to the algebroid
$N^*_{\fZ^\reg_\fg/\fZ_\fg}\simeq \isom_{\cG,\fZ^\reg_\fg}$, which
amounts to an action of $N^*_{\fZ^\reg_\fg/\fZ_\fg}$ on
$\pi^\ren_{\crit,w}$ compatible with its action on
$\wh\fH_\crit\underset{\fZ_\fg}\otimes \fZ^\reg_\fg$.

\medskip

The main result of this section is the following explicit description of 
$\pi^\ren_{\crit,w}$:

\medskip

Recall that $\on{Dist}_w$ denotes the quasi-coherent sheaf
on $\Fl^{\cG}_\fZ$, underlying the left D-module of distributions
on the subscheme $\Fl^{\cG}_{w,\fZ}$. Consider the object
$$\CL^{w(\crho)-\crho}_{\Fl^{\cG}_\fZ}\underset{\CO_{\Fl^{\cG}_\fZ}}\otimes
\on{Dist}_w\in\QCoh(\Fl^\cG_{w,\on{th},\fZ}),$$ which we think of as a
$(\wh\fH_\crit\underset{\fZ_\fg}\otimes \fZ^\reg_\fg)$-module via the
identification \eqref{induction parameters}.  It is naturally
equivariant with respect to $\isom_{\cG,\fZ^\reg_\fg}\simeq
N^*_{\fZ^\reg_\fg/\fZ_\fg}$.

\begin{thm}    \label{iso}
There exists a canonical isomorphism
$$\pi^\ren_{\crit,w}\simeq
\CL^{w(\crho)-\crho}_{\Fl^{\cG}_\fZ}\underset{\CO_{\Fl^{\cG}_\fZ}}\otimes
\on{Dist}_w\in\QCoh(\Fl^\cG_{w,\on{th},\fZ}),$$
compatible with the $N^*_{\fZ^\reg_\fg/\fZ_\fg}$-action.
\end{thm}

\ssec{Application to Wakimoto modules and BGG type resolution} Let
$\BW_{\hslash,w}$ be the $\hslash$-family of Wakimoto modules, induced
from $\pi_{\hslash,w}$. Applying the above renormalization
construction to $\bB_\hslash=\wh\CA_{\fg,\hslash,x}$ and the module
$\BW_{\hslash,w}$ with respect to the ideal $\on{ker}\bigl(\fZ_\fg\to
\fZ^\reg_\fg\bigr)$ we obtain the renormalized family of Wakimoto
modules, denoted $\BW^\ren_{\hslash,w}$.

\medskip

Let $\BW^\ren_{\crit,w}$ denote the fiber of
$\BW^\ren_{\hslash,w}$ at the critical level.  
As in \secref{act of alg on Wak funct} we have:

\begin{lem}
The lattices $\BW^\ren_{\hslash,w}$ and
$\BW_{\hslash}(\pi^\ren_{\hslash,w})$ in the localization of
$\BW_{\hslash,w}$ with respect to $\hslash$ coincide.
\end{lem}

Hence, in particular, $\BW^\ren_{\crit,w}$ is isomorphic to
$\BW_\crit(\pi^\ren_{\crit,w})$. Moreover, by 
\lemref{two actions of renorm}, the above isomorphism respects 
the $U^{\ren,\reg}(\hg_\crit)$-actions.

Combining this with \thmref{full ident}, we obtain an isomorphism
of $U^{\ren,\reg}(\hg_\crit)$-modules
\begin{equation} \label{ren Wak as G}
\BW^\ren_{\crit,w}\simeq \sG(\on{Dist}_w).
\end{equation}

\medskip

Combining \thmref{iso} with \corref{screening BGG} we obtain a right
resolution of $\BV_\crit$:
$$\BV_\crit \to C^0 \overset{\delta^0}\longrightarrow C^1\to...,$$
whose $k$-th term is
$$\underset{w\in W,\ell(w)=k}\oplus\, \BW^\ren_{\crit,w}.$$

Note that for $w=1$ the module $\on{Dist}_w$ is just
$\Fun(\Fl^{\cG}_{1,\fZ})$ and the corresponding Wakimoto module is
$\BW_{\crit,0}$.  In the case when $w$ is a simple reflection, the
modules $\on{Dist}_w$, and the corresponding Wakimoto modules were
constructed in \cite{FF}, and the differential $\delta^0$ was
described there explicitly as the sum of certain degenerations of the
"screening operators" at the critical level. It was conjectured in
\cite{FF} that this complex may be extended to a resolution of
$\BV_\crit$, which is a particular degeneration of the BGG type
resolution of $\BV_\crit$ for generic levels.

Thus, we have obtained a proof of this conjecture. However, it
would be interesting to find explicit formulas for
the higher differentials $\delta^k: C^k \to C^{k+1}$ of this
resolution in terms of the screening operators.

\ssec{}

The rest of this section is devoted to the proof of
\thmref{iso}. Observe that the pull-back
$\isom_{\cG,\fZ^\reg_\fg}|_{\Fl^\cG_{w,\on{th},\fZ}}$ has a natural
structure of algebroid on $\Fl^\cG_{w,\on{th},\fZ}$; we will denote it
by $\isom_{\cG,w,\on{th}}$. Let $\isom_{\cG,w}\subset
\isom_{\cG,w,\on{th}}|_{\Fl^{\cG}_{w,\fZ}}$ be the corresponding
algebroid on $\Fl^{\cG}_{w,\fZ}$ it contains
$\isom_{\cB,w}:=\isom_{\cB,\fZ^\nilp_\fg}|_{\Fl^{\cG}_{w,\fZ}}$ as a
sub-algebroid.

Since the action of $\isom_{\cG,w,\on{th}}$ on
$\Fl^\cG_{w,\on{th},\fZ}$ is transitive, from Kashiwara's theorem
(see, e.g., \cite{FG1}, Sect. 7) we obtain the following:

\begin{lem}
Let $\CT\in \QCoh(\Fl^\cG_{w,\on{th},\fZ})$ be a module over
$\isom_{\cG,w,\on{th}}$.  Assume that the $\CO$-submodule
$\CT^0\subset \CT$ of sections supported over $\Fl^{\cG}_{w,\fZ}$ is
isomorphic as a $\isom_{\cB,w}$-module to
$\CL^\cla_{\Fl^{\cG}_{w,\fZ}}$ for some $\cla\in \cLambda$. Then $\CT$
is isomorphic to 
$\CL^{w(\crho)-\crho+\cla}_{\Fl^{\cG}_\fZ}
\underset{\CO_{\Fl^{\cG}_\fZ}}\otimes\on{Dist}_w$ as a
$\isom_{\cG,w,\on{th}}$-module.
\end{lem}

We claim that the conditions of this lemma are satisfied for
$\CT=\pi^\ren_{\crit,w}$ and $\cla=0$.  Consider the canonical
filtration on $\pi^\ren_{\crit,w}$. For every $k$ we have a map
\begin{align*}
&\bigl(\isom_{\cG,w,\on{th}}|_{\Fl^{\cG}_{w,\fZ}}/\isom_{\cG,w}\bigr)
\underset{\Fun(\Fl^{\cG}_{w,\fZ})}\otimes
\Bigl(\bigl(\pi^\ren_{\crit,w}\bigl)_k/
\bigl(\pi^\ren_{\crit,w}\bigl)_{k-1}\Bigr)\to \\ &\to
\bigl(\pi^\ren_{\crit,w}\bigr)_{k+1}/\bigl(\pi^\ren_{\crit,w}\bigr)_{k},
\end{align*}
which by \eqref{gr renorm} induces an isomorphism
$$\on{gr}^k\bigl(\pi^\ren_{\crit,w}\bigl)\simeq \pi_{\crit,w}
\underset{\fH^{\RS,w_0(w(\rho)-\rho)}_\crit}\otimes
\Sym^k_{\Fun(\Fl^{\cG}_{w,\fZ})}
\bigl(\isom_{\cG,w,\on{th}}|_{\Fl^{\cG}_{w,\fZ}}/\isom_{\cG,w}\bigr).$$

In particular, we obtain that
$\pi_{\crit,w}\underset{\fH^{\RS,w_0(w(\rho)-\rho)}_\crit}\otimes
\Fun(\Fl^{\cG}_{w,\fZ}) \subset \pi^\ren_{\crit,w}$ equals the
subspace consisting of sections supported scheme-theoretically on
$\Fl^{\cG}_{w,\fZ}\subset \Fl^\cG_{w,\on{th},\fZ}$.

\medskip

Hence, it remains to see that
$$\pi_{\crit,w}\underset {\fZ^\nilp_\fg}\otimes \fZ^\reg_\fg\simeq
\pi_{\crit,w}\underset{\fH^{\RS,w_0(w(\rho)-\rho)}_\crit}\otimes
\Fun(\Fl^{\cG}_{w,\fZ})$$ is isomorphic to $\Fun(\Fl^{\cG}_{w,\fZ})$,
as a module over $\isom_{\cB,w}\simeq
\isom_{\cB,\fZ^\nilp_\fg}|_{\Fl^{\cG}_{w,\fZ}}$. The latter
fact follows from the $\nilp$-version of \propref{reg mkdv}.

\end{document}